\documentclass[11pt,reqno]{article}

\usepackage{geometry}
\usepackage{comment}
\usepackage{enumitem} 
\usepackage{lipsum}
\usepackage{changepage}

\usepackage{amsmath}
\usepackage{amsthm}
\usepackage{amssymb}
\usepackage{mathrsfs}

\usepackage{graphicx}
\usepackage{subcaption}
\usepackage{xcolor}
\usepackage{tikz}
\usepackage{tkz-euclide}
\graphicspath{{./figure/}}
\DeclareGraphicsExtensions{.pdf,.jpeg,.png,.jpg}

\usepackage{thmtools}
\usepackage{etoolbox}

\usepackage{hyperref}
\usepackage[capitalize]{cleveref}

\newcommand{\HH}{\mathbb{H}}
\newcommand{\RR}{\mathbb{R}}
\newcommand{\CC}{\mathbb{C}}
\renewcommand{\SS}{\mathbb{S}} 
\newcommand{\ZZ}{\mathbb{Z}}
\newcommand{\FF}{\mathbb{F}}
\newcommand{\BB}{\mathbb{B}}
\newcommand{\EE}{\mathbb{E}}
\newcommand{\NN}{\mathbb{N}}

\newcommand{\mm}{\mathfrak{m}}

\newcommand{\dd}{\,\mathrm{d}}

\newcommand{\ol}[1]{\overline{#1}}

\DeclareMathOperator{\Int}{Int}
\DeclareMathOperator{\diam}{diam}
\DeclareMathOperator{\dist}{dist}

\DeclareMathOperator{\id}{id}

\DeclareMathOperator{\Area}{Area}

\DeclareMathOperator{\Lip}{Lip}
\DeclareMathOperator{\EL}{EL}

\DeclareMathOperator{\CBB}{CBB}

\newcommand{\refcond}[2]{%
  \hyperref[cond:#1]{(#2)}%
}

\numberwithin{equation}{section}

\newtheoremstyle{noparen}
  {3pt}{3pt}        
  {\itshape}        
  {}                
  {\bfseries}       
  {}                
  { }               
  {\thmname{#1}\thmnumber{ #2} \thmnote{\normalfont #3}}

\theoremstyle{noparen}
\newtheorem{theorem}{Theorem}[section]
\newtheorem{lemma}[theorem]{Lemma}
\newtheorem{proposition}[theorem]{Proposition}
\newtheorem{corollary}[theorem]{Corollary}
\newtheorem{conjecture}{Conjecture}

\newtheoremstyle{nopunct}
  {3pt}{3pt}
  {}                
  {}
  {\bfseries}
  {}
  { }
  {\thmname{#1}\thmnumber{ #2} \thmnote{ #3}}

\theoremstyle{nopunct}

\newtheorem{definition}[theorem]{Definition}

\newtheoremstyle{withcitation}
  {\topsep}{\topsep}
  {\itshape}
  {}
  {\bfseries}
  {.}
  { }
  {\thmname{#1}\thmnumber{ #2}\thmnote{ (#3)}}

\theoremstyle{withcitation}

\crefname{lemma}{Lemma}{Lemma}

\definecolor{checkorange}{RGB}{210,105,0}

\newif\ifdraft
\drafttrue

\usepackage{xcolor}
\usepackage{todonotes}

\ifdraft
    \newcommand{\TODO}[1]{\todo[inline,color=yellow!20]{#1}}
    \newcommand{\CHECK}[1]{%
        \textcolor{checkorange}{\textbf{(CHECK: #1)}}%
    }
    \newcommand{\REF}[1]{\textcolor{red}{\textbf{[REF: #1]}}}
    \newcommand{\FIXME}[1]{%
        \par\noindent
        {\color{red}\textbf{PROOF GAP:} #1}
        \par
    }
\else
    \newcommand{\TODO}[1]{}
    \newcommand{\CHECK}[1]{}
    \newcommand{\REF}[1]{}
    \newcommand{\FIXME}[1]{}
\fi

\makeatletter

\renewcommand*\l@subsection{%
  \@dottedtocline{2}{2em}{3em}%
}

\let\oldl@part\l@part
\renewcommand{\l@part}[2]{%
  \oldl@part{\partname\space #1}{#2}%
}

\makeatother
  
\makeatletter
\let\oldthebibliography\thebibliography
\renewcommand{\thebibliography}[1]{%
  \oldthebibliography{#1}%
  \renewcommand{\@biblabel}[1]{\makebox[\labelwidth][r]{[##1]}}%
}
\makeatother

\begin{document}
\title{The Weyl Problem for Planar CBB$(-1)$ Surfaces in Hyperbolic $3$-Space}
\author{Xinrong Zhao}
\date{}
	\maketitle
\begin{abstract}
    In this paper, we study the Weyl problem for noncompact convex surfaces in hyperbolic $3$-space $\mathbb{H}^3$. We prove that any planar $\operatorname{CBB}(-1)$ surface with countably many ends can be isometrically embedded as the boundary of a closed convex subset $X\subset\mathbb{H}^3$ such that each connected component of the ideal boundary of $X$ is either a round disk or a single point. By known rigidity results, such an embedding is unique up to isometries of $\mathbb{H}^3$. Moreover, we show that the assumption of countably many ends can be dropped if Koebe's circle domain conjecture holds.
\end{abstract}

\tableofcontents
\bigskip

\section{Introduction}
Let $\HH^3_P$ denote the Poincaré ball model of hyperbolic $3$-space. A compact subset $K$ in the ideal sphere $\partial \HH^3_P=\widehat{\CC}$ is called a \emph{circle-type closed set} if each of its connected components is either a round disk or a singleton, and a domain (i.e., an open connected subset) $U\subset\widehat{\CC}$ is called a \emph{circle domain} if its complement $\widehat{\CC}\setminus U$ is a circle-type closed set. Throughout this paper, all surfaces are assumed to be connected and without boundary unless explicitly stated otherwise. A surface is called \emph{planar} if it is noncompact and homeomorphic to a domain in the $2$-sphere $\SS^2$. By a $\CBB (-1)$ surface, possibly with boundary, we mean a complete two-dimensional Alexandrov space with curvature bounded below by $-1$.

Let $(S,d_0)$ be a planar $\CBB(-1)$ surface. A \emph{convex realization} of $(S,d_0)$ is a pair $(F,X)$, where $X\subset\HH^3_P$ is a closed convex subset of dimension at least two and $F$ is an isometry from $(S,d_0)$ onto the intrinsic boundary $\partial X$ of $X$. Specifically, when $\dim X=3$, $\partial X$ is the standard topological boundary of $X$ endowed with the intrinsic path metric induced from $\HH^3_P$; in the degenerate case $\dim X=2$, $\partial X$ is understood as the Alexandrov double $D(X)$ of $X$. Moreover, let $\overline{X}$ denote the closure of $X$ in $\RR^3$. If the ideal boundary $\partial_\infty X:=\overline{X}\cap\partial\HH^3_P$ is a circle-type closed set, then $(F,X)$ is called a \emph{circle-type convex realization} of $(S,d_0)$.

The classical Weyl problem asks whether every positively curved $2$-sphere can be isometrically embedded as a convex surface in Euclidean $3$-space. In the hyperbolic setting, an analogous problem for planar $\CBB(-1)$ surfaces was posed by Luo and Wu \cite{LuoWu}, which can be formulated as follows:
\begin{conjecture}[({\cite{LuoWu}})]\label{conj:1}
	Every planar $\CBB(-1)$ surface $(S,d_0)$ admits a circle-type convex realization in $\HH^3_P$, which is unique up to isometries of $\HH^3_P$.
\end{conjecture}

The main result of this paper is the following
general existence theorem.

\begin{theorem}\label{thm:MainThm0}
Let $(S,d_0)$ be a planar $\CBB(-1)$ surface. If $(S,d_0)$ is conformally equivalent to a circle domain, then it admits a circle-type convex realization in $\HH^3_P$.
\end{theorem}

Here, the underlying conformal structure on $(S,d_0)$ is induced by Reshetnyak's isothermal coordinates; see \cite[Sections~7.1--7.3]{Reshetnyak}. By Koebe's general uniformization theorem (see \cref{thm:GUT-Koebe}), the planar surface $S$ is conformally equivalent to a domain in $\widehat{\CC}$. In 1908, Koebe \cite{Koebe1908} conjectured that every domain in $\widehat\CC$ is conformally equivalent to a circle domain, and he verified this conjecture for finitely connected domains \cite{Koebe1920}. He and Schramm \cite{HeSchramm1993} confirmed this conjecture for domains with countably many boundary components. Combining this uniformization theorem with \cref{thm:MainThm0} and the rigidity results in \cite{luo2026rigidity,Zhao-rigidity} establishes \cref{conj:1} for surfaces with countably many ends:
\begin{theorem}\label{thm:MainThm1}
	Every planar $\CBB(-1)$ surface $(S,d_0)$ with countably many ends admits a circle-type convex realization in $\HH^3_P$. Furthermore, such a realization is unique up to isometries of $\HH^3_P$.
\end{theorem}

By Pogorelov's regularity theorem \cite[Chapter~V, Section~8, Theorem~1]{PogorelovBook}, a convex realization is smooth and strictly convex wherever the intrinsic metric on $S$ is smooth with Gaussian curvature strictly greater than $-1$. As an immediate consequence in the smooth category, we obtain:
\begin{theorem}\label{thm:MainThm1-smooth}
	Let $(S, g)$ be a complete planar Riemannian surface with countably many ends and Gaussian curvature $K_g > -1$ everywhere. Then $(S, g)$ is isometric to a smooth, strictly convex surface $Y \subset \HH^3_P$ that bounds a convex body $X \subset \HH^3_P$ whose ideal boundary $\partial_{\infty} X$ is a circle-type closed set with countably many connected components. Moreover, $Y$ is unique up to isometries of $\HH^3_P$.
\end{theorem}

In view of \cref{thm:MainThm0}, the existence part of \cref{conj:1} is a consequence of Koebe's circle domain conjecture. Conversely, Luo and Wu \cite{LuoWu} showed that \cref{conj:1} implies Koebe's circle domain conjecture. We thus arrive at the following equivalence:
\begin{corollary}
	The existence part of \cref{conj:1} is equivalent to Koebe's circle domain conjecture.
\end{corollary}

OpenAI has recently announced a proof of Koebe's circle domain
conjecture \cite{OpenAIKoebe2026}. Assuming the validity of this proof, \cref{thm:MainThm0} completely settles the existence part of the hyperbolic Weyl problem (\cref{conj:1}) for all planar $\operatorname{CBB}(-1)$ surfaces:
\begin{theorem}
    Every planar $\CBB(-1)$ surface $(S,d_0)$ admits a circle-type convex realization in $\HH^3_P$.
\end{theorem}

\subsection{Historical background}

Convex surfaces are classical objects in global differential geometry and Alexandrov geometry, where a central question is the extent to which the intrinsic metric of a surface determines its extrinsic shape. In $\mathbb{R}^3$, the classical Weyl problem asks whether any Riemannian metric of positive Gaussian curvature on $\SS^2$ can be isometrically embedded as a closed, strictly convex surface \cite{Weyl1916}. This was resolved in the real-analytic case by Lewy \cite{Lewy1938} and in the smooth category by Nirenberg \cite{Nirenberg1953}, while global rigidity for smooth strictly convex closed surfaces was established earlier by Cohn-Vossen \cite{CohnVossen1927}. Beyond the smooth setting, Alexandrov \cite{MR13540} proved that every $\CBB(0)$ metric on $\SS^2$ can be realized as the boundary of a compact convex body in $\mathbb{R}^3$, and Pogorelov \cite{Pogorelov} established the corresponding global rigidity. Together, their foundational theorems completely settle the Weyl problem for closed convex surfaces in $\mathbb{R}^3$ within the Alexandrov category.
For complete (possibly noncompact) convex surfaces in $3$-dimensional space forms, Alexandrov \cite[Chapter~XII, Section~2]{Alexandrov2006} showed that a surface with a complete intrinsic path metric can be realized as the boundary of a closed convex set in the simply connected space form of constant curvature $K$ if and only if it is a complete $\CBB(K)$ surface homeomorphic to a domain in $\SS^2$. Nevertheless, this general existence theorem neither guarantees uniqueness nor prescribes the asymptotic behavior at infinity. 

Whereas complete noncompact convex surfaces in $\mathbb{R}^3$ are topologically limited to planes or cylinders and are well understood, their counterparts in $\HH^3_P$ behave much more subtly: such surfaces admit arbitrary planar topologies and rich end structures, while their ideal boundaries introduce additional flexibility at infinity. As observed in \cite{luo2026rigidity}, the following counterexample by Thurston \cite{Thurston} shows that the condition of vanishing one-dimensional Hausdorff measure is optimal. Let $X_i$ be the hyperbolic convex hull of a non-circular Jordan curve $J_i \subset \partial \mathbb{H}^3_P$ ($i=1,2$). Then each $\partial X_i$ is intrinsically isometric to two copies of $\mathbb{H}^2$, which implies that $\partial X_1$ and $\partial X_2$ are isometric. Note that $X_1$ and $X_2$ are related by an ambient hyperbolic isometry if and only if $J_1$ and $J_2$ are related by a M\"obius or anti-M\"obius transformation. Thus, taking $J_1$ to be a square and $J_2$ to be a non-square rectangle yields the desired counterexample. This demonstrates that the intrinsic metric alone cannot detect the extrinsic geometry at infinity.

To restore rigidity and single out canonical realizations, Luo and Wu \cite{LuoWu} formulated the \emph{strong Weyl problem} in hyperbolic $3$-space (stated as \cref{conj:1} above), which asks whether every planar $\CBB(-1)$ surface admits a unique circle-type convex realization up to isometries of $\HH^3_P$. In the case where the intrinsic metric is hyperbolic, Luo and Wu \cite{LuoWu} showed that the existence of a circle-type convex realization for every planar hyperbolic surface is equivalent to Koebe's circle domain conjecture. Combining this equivalence with the theorem of He and Schramm \cite{HeSchramm1993}, they established the existence of such realizations for all planar hyperbolic surfaces with countably many ends. On the rigidity side, recent works by Luo, Luo, and Rao \cite{luo2026rigidity} and by the author \cite{Zhao-rigidity} proved that for noncompact closed convex subsets $X_1, X_2 \subset \HH^3_P$ of dimension at least $2$ whose ideal boundaries are circle-type closed sets with countably many connected components, every intrinsic isometry $\partial X_1 \to \partial X_2$ extends to a global isometry of $\HH^3_P$. Together with the existence result of \cite{LuoWu}, this settles both the existence and uniqueness of circle-type convex realizations for planar hyperbolic surfaces with countably many ends.

Prior to the present work, however, the existence of circle-type convex realizations remained largely open for general planar $\CBB(-1)$ surfaces. In the smooth disk setting, a recent partial result was obtained by Mesbah \cite[Theorem~1.4]{MesbahEmbeddedConvex}, who proved that every complete smooth conformal metric on the unit disk whose curvature lies in a compact subset of $(-1,0)$ and has uniformly bounded covariant derivatives of all orders with respect to the hyperbolic metric admits a convex realization in $\mathbb{H}^3$ asymptotic to any prescribed Jordan
curve at infinity (in particular, a round circle). For further developments concerning noncompact Weyl problems, Fuchsian rigidity, and convex boundaries of hyperbolic manifolds, we refer the reader to \cite{SchlenkerUnboundedWeyl,ProsanovPolyhedralType,ProsanovFuchsianRigidity,ChenSchlenkerBoundaryData}.
    
\subsection{Outline of the proof and organization of the paper}
The proof of \cref{thm:MainThm0} proceeds in three main stages, progressing from compact surfaces with boundary to complete surfaces of finite funnel-type, and finally to general planar $\CBB(-1)$ surfaces.

\smallskip
\noindent
{\bf Step 1.} Convex cap realizations for genus zero compact surfaces. 
Let $k > -1$ and let $(S, d_0)$ be a compact genus zero $\CBB(k)$ surface with nonempty boundary. We prove that $(S,d_0)$ admits an intrinsic isometric embedding $F$ into the boundary of a compact convex body in $\HH^3_P$ such that the image of each boundary component $\Gamma_i \subset \partial S$ lies in a totally geodesic plane $P_i \subset \HH^3_P$, and the orthogonal projection of $F(S)$ onto $P_i$ is contained in the hyperbolic convex hull of $F(\Gamma_i)$ within $P_i$. We refer to such embeddings as \emph{convex cap realizations}.

\smallskip
\noindent
{\bf Step 2.} Planar $\CBB(k)$ surfaces of finite funnel-type. 
Next, we consider a planar $\CBB(k)$ surface $(S,d_0)$ ($k > -1$) of \emph{finite funnel-type}, that is, $S$ has finitely many ends and admits a convex realization $(F,X)$ in $\HH^3_P$ under which each component of $\partial_{\infty} X$ in the ideal sphere is nondegenerate. We construct an exhaustion $\{S_j\}_{j=1}^\infty$ of $S$ by compact $\CBB(k)$ subsurfaces with boundary. By Step~1, each $S_j$ admits a convex cap realization $F_j$. Passing to a convergent subsequence of $F_j(S_j)$, we obtain a limiting convex realization $F$ of $(S,d_0)$. Combining the convergence theory of convex surfaces with Schramm's transboundary extremal length, we prove that $F$ is indeed a circle-type convex realization.

\smallskip
\noindent
{\bf Step 3.} General planar $\CBB(-1)$ surfaces.
Finally, let $(S,d_0)$ be an arbitrary planar $\CBB(-1)$ surface conformally equivalent to a circle domain. We exhaust $S$ by finite-type subdomains and conformally approximate the metric $d_0$ on these subdomains $V_j$ by complete metrics of finite funnel-type with curvature bounded below by $k_j \downarrow -1$. Applying Step~2 yields a sequence of circle-type convex realizations for the approximating metrics. By the convergence theorems established in \cref{subsec:convergence}, after passing to a subsequence, these realizations converge to a convex realization of $(S,d_0)$. Adapting the transboundary extremal length argument by Luo and Wu in \cite{LuoWu}, we conclude that the limiting ideal boundary is a circle-type closed set, thereby completing the proofs of \cref{thm:MainThm0} and \cref{thm:MainThm1}.

\smallskip

The paper is organized as follows. In \cref{sec:General-Prelimiaries}, we review background material on $\CBB(-1)$ surfaces, transboundary extremal length, and the convergence theory for convex realizations in $\HH^3_P$, all of which will be used throughout this paper. Steps~1, 2, and~3 are then carried out in \cref{sec:Convex-Cap-Realizations,sec:finite-type-Realizations,sec:general-Realizations}, respectively. Finally, in \cref{sec:appendix-convergence,sec:appendix}, we collect several auxiliary results concerning the convergence of convex surfaces (used in \cref{subsec:convergence}) and the complete conformal blow-up construction (used in \cref{subsec:complete-blowup}).

\medskip
\noindent
\textbf{Acknowledgments.}  
I am grateful to Professor Feng Luo, Yanwen Luo, and Zhenghao Rao for helpful discussions on related topics during our earlier collaboration. I would also like to express my sincere gratitude to my advisor, Professor Bobo Hua, for his valuable suggestions and continuous support.

\medskip
\noindent
{\bf Use of AI.} 
During the preparation of this manuscript, the author used generative AI tools for language polishing and literature exploration. All arguments in the main text (Sections~1--5) were developed independently by the author without AI assistance. Some auxiliary results, developed by the author based on initial ideas suggested by AI tools, are presented in two appendices. The author takes full responsibility for the content of this article.

\section{General Preliminaries}\label{sec:General-Prelimiaries}
\subsection{$\CBB(-1)$ surfaces}
In this paper, we adopt the classical triangle-comparison formulation for $\CBB(-1)$ spaces; see \cite[Theorem~8.14\textup{(b)}]{AKP24}. Note that under this convention, a $\CBB(-1)$ space is complete by definition.

\begin{definition}
A complete geodesic metric space $(\mathcal{L},d)$ satisfies the
\emph{$\CBB(-1)$ comparison} if for every geodesic triangle $\triangle pxy$ in $\mathcal{L}$ with minimizing sides $[px]$, $[xy]$, and $[yp]$, and any comparison triangle $\triangle\bar{p}\bar{x}\bar{y}\subset\HH^2$ with the same side lengths, the inequality
$$
 d(p,z)\ge d_{\HH^2}(\bar{p},\bar{z})
$$
holds for every point $z\in[xy]$, where $\bar{z}\in[\bar{x}\bar{y}]$ is the unique point satisfying $d_{\HH^2}(\bar{x},\bar{z})=d(x,z)$.
The space is \emph{locally $\CBB(-1)$} if every point admits a neighborhood in which this comparison holds for all geodesic triangles contained therein.

A \emph{$\CBB(-1)$ surface} is a connected, complete geodesic metric space satisfying the $\CBB(-1)$ comparison whose metric topology is that of a topological $2$-manifold, possibly with boundary. Unless stated otherwise, surfaces are assumed to be connected and without boundary. For an arbitrary $k\in\mathbb{R}$, we define $\CBB(k)$ analogously using the complete, simply connected surface of constant Gaussian curvature $k$ as the model space. When $k>0$, the comparison is required only for geodesic triangles of perimeter less than $2\pi/\sqrt{k}$.
\end{definition}

The globalization theorem for $\CBB(-1)$ spaces takes the following form.

\begin{theorem}[{(\cite[Theorem~8.31]{AKP24})}]\label{thm:Globalization-theorem}
Any complete locally $\CBB(-1)$ length space is $\CBB(-1)$.
\end{theorem}

\begin{definition}\label{def:metric-double}
Let $S$ be a length surface with nonempty boundary. Its
\emph{metric double} $D(S)$ is obtained by gluing two copies of
$S$ pointwise along their boundaries and equipping the quotient with the induced
length metric. The canonical involution interchanging the two copies is denoted
by $\iota:D(S)\to D(S)$. We call $D(S)$ a \emph{planar metric
double} if $S$ is isometric to a two-dimensional closed convex
subset of $\HH^2$ with nonempty boundary.
\end{definition}

For a $\CBB(-1)$ surface, the Alexandrov boundary coincides with
the manifold boundary \cite[Corollary~2.12]{Alexandrov-boundary}.
Consequently, Perelman's doubling theorem (\cite[Theorem~2.1]{Petrunin-Gluing}) yields the following result.

\begin{corollary}\label{cor:perelman-double}
Let $(S,d_0)$ be a compact $\CBB(-1)$ surface with nonempty
boundary. Then $D(S)$ is a closed, connected $\CBB(-1)$ surface.
\end{corollary}

Throughout this paper, an \emph{isometric embedding} $F:(X,d_X)\to (Y,d_Y)$ refers to a map that is a distance-preserving homeomorphism onto its image $(F(X),d_Y|_{F(X)\times F(X)})$. By contrast, the standard differential-geometric notion of an isometric embedding, which preserves path lengths rather than extrinsic distances, is called an \emph{intrinsic isometric embedding}:

\begin{definition}
Let $(X,d_X)$ be a length space and $(Y,d_Y)$ a metric space.
A topological embedding $f:X\to Y$ is called an
\emph{intrinsic isometric embedding} if it is an isometry
from $(X,d_X)$ onto $(f(X),d_{f(X)})$, where $d_{f(X)}$ denotes the induced length metric on $f(X)$ given by
$$
 d_{f(X)}(p,q):= \inf_{\gamma}l_{d_Y}(\gamma),
 \qquad p,q\in f(X),
$$
where $l_{d_Y}(\gamma)$ denotes the length of $\gamma$ with respect to $d_Y$ and the infimum is taken over all rectifiable curves $\gamma$ in $f(X)$ connecting $p$ to $q$.
\end{definition}

\subsection{Extremal length and the finite-type classification}
\label{sec:prelim-finite-type}
We fix an orientation on each planar surface. Let $(S,d)$ be an oriented $\CBB(-1)$ surface, possibly with boundary. Equipped with the induced length metric, its interior $\operatorname{Int} S$ has locally bounded integral curvature in the sense of Alexandrov, and Reshetnyak's isothermal-coordinate theorem endows $\operatorname{Int} S$ with a canonical conformal structure compatible with the orientation, which is defined on all of $S$ when $\partial S=\varnothing$; see \cite[Sections~7.1--7.3]{Reshetnyak} and \cref{sec:appendix} for the isothermal representation. In particular, if $S$ is planar, this conformal structure makes $\operatorname{Int} S$ into a planar Riemann surface, to which Koebe's general uniformization theorem below applies \cite{GUTofKoebe}.

\begin{theorem}\label{thm:GUT-Koebe}
Every planar Riemannian surface is biholomorphic to a domain in
the Riemann sphere $\widehat{\CC}$.
\end{theorem}

\begin{definition}
\label{def:finite-funnel-conformal}
A planar $\CBB(-1)$ surface is said to be of \emph{finite
funnel-type} if, when equipped with its canonical conformal structure, it is conformally equivalent
to a finitely connected circle domain $U\subset\widehat{\CC}$ such that every connected component of $\widehat{\CC}\setminus U$ is a closed round disk.
\end{definition}

Fix a smooth conformal background metric $g$ on a Riemann
surface $S$. For a nonnegative Borel measurable function $\rho$ and a curve family $\Gamma$ in $S$, we set
$$l_\rho(\gamma)=\int_\gamma\rho\,\dd s_g,\qquad A_\rho=\int_S\rho^2\,\dd A_g,\qquad l_\rho(\Gamma)=\inf_{\gamma\in\Gamma}l_\rho(\gamma).$$
The \emph{extremal length} of $\Gamma$ in $S$ is defined by
$$\EL_S(\Gamma) =\sup_{0<A_\rho<\infty} \frac{l_\rho(\Gamma)^2}{A_\rho}.$$
Clearly, $\EL_S(\Gamma)$ depends only on the conformal class of $g$. Moreover, via the canonical isothermal representation, the estimate of extremal length on an oriented $\CBB(-1)$ surface can be derived by computing directly in terms of its intrinsic length element and area measure; see \cite[Section~2.3 and Appendix]{LuoWu} for details.

By an argument similar to that of \cite[Proposition~4.2]{Zhao-rigidity}, we have the following criterion:

\begin{proposition}\label{prop:finite-end-extremal-length}
Let $X \subset\HH^3_P$ be a closed convex set whose ideal boundary $\partial_{\infty} X$ has finitely many connected components, and let $D$ be a connected component of $\partial_{\infty} X$. Fix a closed topological disk $A \subset \partial X$, and let $\Gamma$ be the family of simple closed curves in $\partial X$ separating $D$ from $A$. Then $\EL_{\partial X}(\Gamma)=0$ if and only if $D$ is a singleton.
\end{proposition}

For a finitely connected circle domain $U\subset\widehat{\CC}$, the analogous
criterion follows immediately from the standard modulus of a round annulus: the
extremal length of the family of simple closed curves separating a complementary component of $U$ from a fixed closed disk $A\subset U$ vanishes if and only if that
component is a single point. Combining this with Koebe's uniformization theorem for finitely connected domains and \cref{prop:finite-end-extremal-length}, we arrive at the following characterization:

\begin{corollary}\label{cor:finite-funnel-type}
Let $(S,d_0)$ be a planar $\CBB(-1)$ surface. The
following are equivalent:
\begin{enumerate}
\item[(1)] $S$ is of finite funnel-type.
\item[(2)] For every convex realization $(F,X)$ of $S$, the ideal
boundary $\partial_\infty X$ consists of finitely many connected components,
none of which is a singleton.
\end{enumerate}
\end{corollary}

\subsection{Transboundary extremal length}

In this subsection, we review the notion of transboundary extremal length, which was introduced by O.~Schramm \cite{TransboundaryLength} and plays a central role in this work. 

For a noncompact planar surface $S$, let $E(S)$ denote its space
of ends and
$$
 \mathcal{E}(S)=S\cup E(S)
$$
its end compactification, which is homeomorphic to the $2$-sphere. When $S$
is realized as a domain $U\subset\widehat{\CC}$, its ends are in natural bijection with the connected components of $\widehat{\CC}\setminus U$, and
$\mathcal{E}(U)$ can be identified with the quotient space of $\widehat{\CC}$ obtained by collapsing each complementary component to a single point.
Notably, every homeomorphism between two domains in $\widehat{\CC}$ extends uniquely to a homeomorphism between their
end compactifications; see the opening paragraphs of
\cite[Section~2]{NtalampekosRajala2025}.

Let $S$ be a planar Riemannian surface equipped with a smooth conformal background metric $g$. An \emph{extended metric} on $S$ is a pair $\mm=(\rho\, \dd s_g, \mu)$, where $\rho:S\to\mathbb{R}_{\geq 0}$ is a Borel measurable function and $\mu: E(S)\to\mathbb{R}_{\geq 0}$ is a nonnegative function on the space of ends. 
The \emph{area} of the extended metric $\mm$ is defined as 
$$
A(\mm)=\int_S \rho^2\, \dd A_g +\sum_{e\in E(S)} \mu(e)^2,
$$ 
where $\dd A_g$ is the area element of $g$. 
A \emph{transboundary curve} in $S$ is a continuous map $\gamma$ from an interval into $\mathcal{E}(S)$. The \emph{length} of $\gamma$ with respect to the extended metric $\mm$ is defined by 
$$
l_{\mm}(\gamma)=\int_{\gamma\cap S}\rho\,\dd s_g+\sum_{e\in E(S)\cap \gamma} \mu(e),
$$
where $\dd s_g$ is the line element of $g$ and $\int_{\gamma\cap S}\rho\,\dd s_g$ is understood as $+\infty$ if $\gamma\cap S$ is not locally rectifiable. For a family $\Gamma$ of transboundary curves in $S$, its length with respect to $\mm$ is given by 
$$
l_\mm(\Gamma):=\inf_{\gamma\in\Gamma} l_\mm(\gamma),
$$
and Schramm's \emph{transboundary extremal length} of $\Gamma$ in $S$ is defined as
$$
\EL_S(\Gamma):=\sup_{0<A(\mm)<\infty} \frac{l_\mm(\Gamma)^2}{A(\mm)},
$$
where the supremum is taken over all extended metrics $\mm$ on $S$ with finite positive area. The following two lemmas follow immediately from the definition:

\begin{lemma}\label{lem:transboundary-include-relation}
Let $\Gamma_1$ and $\Gamma_2$ be two curve families in $\mathcal{E}(S)$. If every $\gamma_1 \in \Gamma_1$ contains a subcurve $\gamma_2 \in \Gamma_2$, then $\EL_S(\Gamma_1) \ge \EL_S(\Gamma_2)$. In particular, if $\Gamma_1 \subset \Gamma_2$, then $\EL_S(\Gamma_1) \ge \EL_S(\Gamma_2)$.
\end{lemma}

\begin{lemma}\label{lem:conformal-invariance}
If $\varphi : S \to S'$ is a conformal or anti-conformal homeomorphism between planar Riemannian surfaces, then for any curve family $\Gamma$ in $\mathcal{E}(S)$,
$$
\EL_S(\Gamma) = \EL_{S'}(\varphi(\Gamma)),
$$
where $\varphi(\Gamma) := \{\varphi(\gamma) : \gamma \in \Gamma\}$ via the homeomorphic extension of $\varphi$ to $\mathcal{E}(S)\to\mathcal{E}(S')$.
\end{lemma}

For convenience, we also introduce a restricted variant of transboundary extremal length. Given a distinguished end $e_0\in E(S)$ and a curve family $\Gamma$ in $\mathcal{E}(S)$, we define the \emph{restricted transboundary extremal length} of $\Gamma$ with respect to $e_0$ by
$$
\EL_{S;e_0}(\Gamma):=\sup_{\substack{0<A(\mm)<\infty,\ \mu(e_0)=0}} \frac{l_\mm(\Gamma)^2}{A(\mm)},
$$ 
where the supremum is taken over all extended metrics $\mm=(\rho\,\dd s_g,\mu)$ on $S$ with finite positive area satisfying $\mu(e_0)=0$. Clearly, \cref{lem:transboundary-include-relation,lem:conformal-invariance} remain valid for this restricted version with $\varphi(e_0)$ as the distinguished end of $S'$ in \cref{lem:conformal-invariance}, and both $\EL_S(\Gamma)$ and $\EL_{S;e_0}(\Gamma)$ coincide with the classical extremal length when every curve in $\Gamma$ lies entirely in $S$.

The following extremal length estimate due to Luo and Wu will be useful in our subsequent arguments.

\begin{proposition}[{\cite[Proposition~6.5]{LuoWu}}]\label{prop:fin-multiring}
Let $S$ be a topological annulus in $\widehat{\mathbb{C}} \setminus \{0\}$, $W \subset S$ be a finite disjoint union of closed round disks and points, and $\Gamma$ be the family of all simple closed curves in $\mathcal{E}(S\setminus W)$ separating the two ends of $S$. If $S$ contains $N$ disjoint round annuli $R_i = \{z\in\mathbb{C} : r_i < |z| < 2r_i\}$ ($i = 1, 2, \dots, N$) such that
\begin{enumerate}
    \item[(1)] each annulus $R_i$ separates the two ends of $S$, and
    \item[(2)] no connected component of $W$ intersects two distinct annuli $R_i$ and $R_j$,
\end{enumerate}
then
$$
\EL_{\widehat{\CC}\setminus W}(\Gamma) \leq \frac{216\pi}{N}.
$$
\end{proposition}

Let $\{U_n\}$ be a sequence of finitely connected circle domains in $\widehat{\mathbb{C}}$, and let $U \subset \widehat{\mathbb{C}}$ be a circle domain such that every compact subset of $U$ is contained in $U_n$ for all sufficiently large $n$. For each $n$, let $D_n$ be a complementary component of $U_n$ and $z_n$ be a point in $U_n$. Suppose that $D_n$ converges in the spherical Hausdorff metric to a complementary component $D$ of $U$, and that $z_n$ converges spherically to a point $z \in D$. After applying a M\"obius transformation, we may assume that $z=0$ and $\infty\in U$. 

\begin{lemma}\label{lem:disjoint-rings}
Fix $r_* > 0$. There exists a sequence of positive radii $\{r_i\}_{i=1}^\infty$ with $r_1 < r_*$ and $r_{i+1} < r_i/2$ for all $i \ge 1$ such that for each $j \in \mathbb{N}$, whenever $n$ is sufficiently large, every connected component of $\widehat{\CC}\setminus U_n$ other than $D_n$ intersects at most one of the round annuli
$$
A_{n,i} = \{\xi \in \mathbb{C} : r_i < |\xi - z_n| < 2r_i\}, \qquad 1 \leq i \leq j.
$$
Moreover, if $D=\{0\}$, then for each $j\in\mathbb{N}$, $D_n$ is disjoint from $A_{n,1},\dots,A_{n,j}$ for all sufficiently large $n$.
\end{lemma}

\begin{proof}
We first claim that for every $r>0$, there exist
$\delta\in (0,r/2)$ and $N\in\NN$ such that for every $n\ge N$ and
every connected component $R_n\ne D_n$ of $\widehat{\CC}\setminus U_n$,
$$
R_n\cap B(z_n,\delta)\ne\varnothing
\quad\implies\quad
R_n\subset B(z_n,r).
$$
Suppose, for contradiction, that there exist a subsequence $n_k\to\infty$ and connected components $R_{n_k}\ne D_{n_k}$ of $\widehat{\CC}\setminus U_{n_k}$ such that both $R_{n_k}\cap B(z_{n_k},1/k)$ and $R_{n_k}\setminus B(z_{n_k},r)$ are nonempty for all $k\in\NN$. After passing to a further subsequence if necessary, we may assume that $R_{n_k}$ converges in the spherical Hausdorff metric to a nondegenerate closed round disk $R$ containing $0$. Since every compact subset of $U$ is contained in $U_n$ for all sufficiently large $n$, we must have $R\subset\widehat{\CC}\setminus U$. Because $R$ is connected and $0\in R\cap D$, it follows that $R\subset D$. This is impossible if $D=\{0\}$. If instead $D$ is a nondegenerate disk, then $\operatorname{Int}R$ is a nonempty open subset of $\operatorname{Int}D$. Consequently, any point in $\operatorname{Int}R$ belongs to both $R_{n_k}$ and $D_{n_k}$ for all sufficiently large $k$, contradicting the fact that $R_{n_k}\ne D_{n_k}$. This proves the claim.

Now choose an arbitrary $r_1\in (0,r_*)$.
We construct the radii $r_j$ and integers $N_j$ inductively
such that for every $n\ge N_j$, each complementary
component of $U_n$ other than $D_n$ meets at most one
of the annuli $A_{n,1},\ldots,A_{n,j}$. For $j=1$, the condition is vacuous, so we may take
$N_1=1$. Suppose that $r_1,\ldots,r_j$ and $N_j$
have been constructed with the desired property.
Applying the claim to $r=r_j$, we obtain
$\delta_j\in (0,r_j/2)$ and $M_j\in\NN$ such that
for every $n\ge M_j$ and every complementary component
$R_n\ne D_n$ of $U_n$,
$$
R_n\cap B(z_n,\delta_j)\ne\varnothing
\quad\implies\quad
R_n\subset B(z_n,r_j).
$$
We then choose
$$
r_{j+1}\in (0,\delta_j/2)
\qquad\text{and}\qquad
N_{j+1}=\max\{N_j,M_j\}.
$$
Fix $n\ge N_{j+1}$ and a complementary component
$R_n\ne D_n$ of $U_n$. Since $2r_{j+1}<\delta_j$, we have $A_{n,j+1}\subset B(z_n,\delta_j)$. Thus, if $R_n$ meets $A_{n,j+1}$, the choice of $\delta_j$ yields $R_n\subset B(z_n,r_j)$. Since $r_i\ge r_j$
for all $1\le i\leq j$, this implies that $R_n\cap A_{n,i}=\varnothing$ for all $1\leq i\leq j$.
On the other hand, if $R_n$ does not meet $A_{n,j+1}$, then the induction
hypothesis ensures that $R_n$ meets at most one of $A_{n,1},\ldots,A_{n,j}$. Hence the assertion holds for
$j+1$, completing the induction. By construction, $r_{i+1}<\delta_i/2<r_i/4<r_i/2$. Finally, if $D=\{0\}$, then since both $D_n\to \{0\}$ and $z_n\to 0$, we have $D_n\subset B(z_n,r_j)$ for all sufficiently large $n$, so $D_n$ is disjoint from all the annuli $A_{n,1},\ldots,A_{n,j}$. This completes the proof.
\end{proof}

For any circle domain $U \subset \widehat{\mathbb{C}}$, each connected component $P$ of $\widehat{\mathbb{C}} \setminus U$ corresponds canonically to an end $e_P \in E(U)$. Under this identification, we simply write $\EL_{U;P}$ in place of the restricted transboundary extremal length $\EL_{U;e_P}$.

\begin{proposition}\label{prop:sourse_vanishing}
Let $U_n,U,D_n$, and $D$ be as above. Suppose that $\{z_n\}$ and $\{z_n'\}$ are two sequences in $\widehat{\CC}$ such that $z_n,z_n'\in U_n$ for every $n\in\NN$ and both converge to a point $z\in D$ as $n\to\infty$.
Fix a nondegenerate closed round disk $B_*\subset U$, and let $\Gamma_n$ be the family
of simple closed curves in $\mathcal{E}(U_n)$ separating
$\{z_n,z_n'\}$ from $B_*$. Then
$$
\lim_{n\to\infty}\EL_{U_n;D_n}(\Gamma_n)=0.
$$
\end{proposition}

\begin{proof}
By applying a suitable M\"obius transformation, we may assume without loss of generality that $z = 0$ and $\infty \in \operatorname{Int} B_*$. First, we show that for each $n$, there exists a simple arc $\alpha_n \subset U_n$ joining $z_n$ to $z_n'$ such that $\alpha_n \subset \overline{B}\big(z_n, 8|z_n - z_n'|\big)$. Indeed, when $z_n\neq z_n'$, we start with the Euclidean line segment $[z_n, z_n']$. For each complementary component of $U_n$ that intersects $[z_n,z_n']$, we replace the chord inside the disk by the shorter arc of its boundary circle. By slightly pushing these circular arcs into $U_n$, we obtain a simple arc $\alpha_n \subset U_n$ satisfying the desired inclusion.

Let 
$$
A_{n,i} = \{\xi \in \mathbb{C} : r_i < |\xi - z_n| < 2r_i\}, \qquad n,i\in\NN,
$$
be the round annuli given by \cref{lem:disjoint-rings}. Shrinking $r_* > 0$ if necessary, we may assume that $B_*\subset\{\xi\in\widehat{\CC} : |\xi - z_n| > 2r_1\}$ for all sufficiently large $n$. Fix an integer $j\geq 1$. Since $8|z_n-z_n'|\to 0$ as $n\to\infty$, the inclusion above implies
that 
$$
\alpha_n\subset\{\xi\in\CC : |\xi-z_n|<r_j\}\subset\{\xi\in\CC : |\xi-z_n|<r_i\}
$$ 
for all $1\leq i\leq j$ and all sufficiently large $n$. Consequently, for all sufficiently large $n$, the topological annulus $\mathcal{A}_n:=\widehat{\CC}\setminus (B_*\cup\alpha_n)$
contains the first $j$ annuli $A_{n,1},\dots,A_{n,j}$, each of which separates the two ends of $\mathcal{A}_n$. We now distinguish two cases.

\medskip
\noindent\textbf{Case 1.} $D$ is a single point. By \cref{lem:disjoint-rings}, for all sufficiently large $n$, every connected component of $\widehat{\CC}\setminus U_n$ (including $D_n$) intersects at most one of the annuli $A_{n,1},\dots,A_{n,j}$. Let $\Lambda_n$ be the family of simple closed curves in
$\mathcal{E}(\mathcal{A}_n\cap U_n)$ separating the two ends of $\mathcal{A}_n$. Since every curve in $\Lambda_n$ belongs to $\Gamma_n$, \cref{lem:transboundary-include-relation} and \cref{prop:fin-multiring} yield
$$
\EL_{U_n;D_n}(\Gamma_n)\leq \EL_{U_n;D_n}(\Lambda_n)\leq \EL_{U_n}(\Lambda_n)\leq \frac{216\pi}{j}
$$
for all sufficiently large $n$.

\medskip
\noindent\textbf{Case 2.} $D$ is a nondegenerate closed round disk. Then
$D_n$ is also a nondegenerate closed round disk for all sufficiently large $n$. Set $V_n:=U_n\cup D_n$, and let $\Lambda_n'$ be the family of simple closed curves in
$\mathcal{E}(\mathcal{A}_n\cap V_n)$ separating the two ends of $\mathcal{A}_n$. By \cref{lem:disjoint-rings}, every connected component of $\widehat{\CC}\setminus V_n$ intersects at most one of the first $j$ annuli $A_{n,1},\dots,A_{n,j}$ when $n$ is sufficiently large. Applying \cref{prop:fin-multiring} to $V_n$ gives
$$
\EL_{V_n}(\Lambda_n') \leq \frac{216\pi}{j}.
$$

We now bound $\EL_{U_n;D_n}(\Gamma_n)$ by $\EL_{V_n}(\Lambda_n')$. Fix a smooth conformal background metric $g_{V_n}$ on $V_n$, and let $g_{U_n}$ be its restriction to $U_n$.
Take any extended metric $\mm=(\rho\,\dd s_{g_{U_n}},\mu)$ on
$U_n$ with $0<A(\mm)<\infty$ and $\mu(e_{D_n})=0$, where $e_{D_n}\in E(U_n)$ is the end of $U_n$ corresponding to $D_n$.
Identifying $E(V_n)$ naturally with
$E(U_n)\setminus\{e_{D_n}\}$, we define an extended metric $\widetilde{\mm}=(\widetilde{\rho}\,\dd s_{g_{V_n}},\widetilde{\mu})$
on $V_n$ by setting
$$
\widetilde{\rho}=
\begin{cases}
\rho & \text{on }U_n,\\
0 & \text{on }D_n,
\end{cases}
\qquad \text{and} \qquad
\widetilde{\mu}(e)=\mu(e)\quad\text{for all } e\in E(V_n).
$$
Since $\mu(e_{D_n})=0$, the area of $\widetilde{\mm}$ satisfies
$$
A(\widetilde{\mm})=
\int_{V_n}\widetilde{\rho}^2\,\dd A_{g_{V_n}}
+\sum_{e\in E(V_n)}\widetilde{\mu}(e)^2=
\int_{U_n}\rho^2\,\dd A_{g_{U_n}}
+\sum_{e\in E(U_n)}\mu(e)^2
=A(\mm).
$$
Next, we show that for every $\gamma\in\Lambda_n'$, there exists a curve $\eta\in\Gamma_n$ such that $l_\mm(\eta)\leq l_{\widetilde{\mm}}(\gamma)$.
Consider the natural quotient map $\pi_n:\mathcal{E}(V_n)\to\mathcal{E}(U_n)$ obtained by collapsing $D_n$ to the end $e_{D_n}$. If $\gamma \cap D_n=\varnothing$, then $\eta:=\pi_n\circ\gamma$ is already a simple closed curve in
$\mathcal{E}(\mathcal{A}_n\cap U_n)$ separating the two ends of $\mathcal{A}_n$, so $\eta\in\Gamma_n$ and $l_\mm(\eta)=l_{\widetilde{\mm}}(\gamma)$. Otherwise, $\pi_n(\gamma)$ is a union of simple closed curves in $\mathcal{E}(\mathcal{A}_n\cap U_n)$ based at $e_{D_n}$. Since $\gamma$ separates the two ends of $\mathcal{A}_n$ in $\mathcal{E}(\mathcal{A}_n \cap V_n)$, at least one of these simple closed curves, which we denote by $\eta$, must separate the two ends of $\mathcal{A}_n$ in $\mathcal{E}(\mathcal{A}_n \cap U_n)$. Clearly, $\eta\in\Gamma_n$. Since $\widetilde{\rho}=0$ on $D_n$ and $\mu(e_{D_n})=0$, we have $l_\mm(\eta)\leq l_{\widetilde{\mm}}(\gamma)$. 
Consequently, for any such extended metric $\mm$ on $U_n$,
$$
\frac{l_\mm(\Gamma_n)^2}{A(\mm)}
=\frac{\inf_{\eta\in\Gamma_n} l_\mm(\eta)^2}{A(\widetilde{\mm})}
\leq\frac{\inf_{\gamma\in\Lambda_n'} l_{\widetilde{\mm}}(\gamma)^2}{A(\widetilde{\mm})}
\leq\EL_{V_n}(\Lambda_n').
$$
Taking the supremum over all admissible extended metrics $\mm$ yields
$$
\EL_{U_n;D_n}(\Gamma_n)\leq
\EL_{V_n}(\Lambda_n')
\leq \frac{216\pi}{j}
$$
for all sufficiently large $n$.

\medskip
In both cases, we conclude that $\limsup_{n\to\infty}\EL_{U_n;D_n}(\Gamma_n)
\leq \frac{216\pi}{j}$. Since $j\geq 1$ is arbitrary, letting $j\to\infty$ gives $\lim_{n\to\infty}\EL_{U_n;D_n}(\Gamma_n)=0$. This completes the proof.
\end{proof}

\subsection{Convergence of Convex Realizations in hyperbolic space}\label{subsec:convergence}

In this subsection, we establish several results on the convergence of convex surfaces in $\HH^3$, which will be extensively used in the later argument.

First, we recall the definition of the Hausdorff distance and Gromov--Hausdorff convergence. For nonempty compact subsets $A,B$ of a metric space $(X,d)$, their \emph{Hausdorff distance} in $X$ is defined by
$$
d_H^X(A,B)
=
\max\left\{
\sup_{a\in A}d(a,B),
\sup_{b\in B}d(b,A)
\right\}.
$$
The following compactness theorem is well known; see, e.g., \cite[Theorem 7.3.8]{BBI} for a proof.

\begin{theorem}[(Blaschke)]
\label{thm:blashke-rewrite}
If $Z$ is a compact metric space, then every sequence of
closed subsets of $Z$ has a subsequence that converges in the Hausdorff distance in $Z$ to a closed subset.
\end{theorem}
For any two nonempty compact metric spaces $(X,d_X)$ and $(Y,d_Y)$,
their \emph{Gromov--Hausdorff distance} is defined as
$$
 d_{\mathrm{GH}}\bigl((X,d_X),(Y,d_Y)\bigr)
 :=
 \inf_{Z,\iota_X,\iota_Y}
 d_H^Z\bigl(\iota_X(X),\iota_Y(Y)\bigr),
$$
where the infimum is taken over all metric spaces $(Z,d_Z)$ and all isometric embeddings $\iota_X:(X,d_X)\to(Z,d_Z)$ and $\iota_Y:(Y,d_Y)\to(Z,d_Z)$. A sequence of nonempty compact metric spaces $(X_n,d_{X_n})$ is said to converge to $(X,d_X)$ in the Gromov--Hausdorff sense, denoted by $(X_n,d_{X_n})\xrightarrow{\mathrm{GH}}(X,d_X)$, if $d_{\mathrm{GH}}\bigl((X_n,d_{X_n}),(X,d_X)\bigr)\to0$.

We work in the Klein model $\HH^3_K=\BB^3$, where $\BB^3$ is the open unit ball in $\RR^3$. Unless otherwise specified, the Hausdorff convergence of compact convex sets is taken with respect to the Euclidean metric. For a nonempty compact convex set $K\subset\RR^3$, we adopt the following boundary convention:
\begin{equation}\label{ccv:boundary-convention}
 \mathcal{B}(K):=
 \begin{cases}
  \text{the boundary $\partial K$ of $K$ in $\RR^3$},&\text{if }\dim K=3,\\
  \text{the metric double }D(K),&\text{if }\dim K=2,\\
  K,&\text{if }\dim K\leq1.
 \end{cases}
\end{equation}
We equip $\mathcal{B}(K)$ with its induced Euclidean intrinsic metric $d_{\mathcal{B}(K)}$. More precisely, $d_{\mathcal{B}(K)}$ is the path metric on $\mathcal{B}(K)$ induced from $\RR^3$ when $\dim K=3$, the double metric on $D(K)$ when $\dim K=2$, and the restricted Euclidean metric on $K$ when $\dim K\leq 1$. Thus, $(\mathcal{B}(K),d_{\mathcal{B}(K)})$ is a compact geodesic space. The coordinate map $\pi_K:\mathcal{B}(K)\to\RR^3$ is defined as the inclusion map when $\dim K=3$, the canonical folding map when $\dim K=2$, and the identity map when $\dim K\leq 1$. When $\dim K=2$, we write $x^+$ and $x^-$ for the two copies of $x\in K$ in $K^+$ and $K^-$, respectively, with the identification $x^+=x^-$ for every $x$ in the relative boundary $\partial_{\mathrm{rel}}K$ of $K$. The double metric is given explicitly by
\begin{equation}\label{ccv:double-formulas}
 \begin{aligned}
 d_{\mathcal{B}(K)}(x^+,y^+)&=d_{\mathcal{B}(K)}(x^-,y^-)=\|x-y\|,\\
 d_{\mathcal{B}(K)}(x^+,y^-)
  &=\min_{z\in\partial_{\mathrm{rel}}K}
       \bigl(\|x-z\|+\|z-y\|\bigr).
 \end{aligned}
\end{equation}
In all dimensions, $\pi_K$ is $1$-Lipschitz and preserves the lengths of rectifiable curves.

For $K\subset\overline{\BB^3}$, set
$$
 \mathcal{B}_{\HH}(K):=\pi_K^{-1}(\BB^3).
$$
We denote by $\delta_{\mathcal{B}(K)}$ the hyperbolic intrinsic metric on $\mathcal{B}_{\HH}(K)$, defined by measuring the length of each curve $\gamma$ in $\mathcal{B}_{\HH}(K)$ as the hyperbolic length of its image $\pi_K(\gamma)$ in the Klein model $\HH^3_K$. Accordingly, $\delta_{\mathcal{B}(K)}$ coincides with the induced path metric from $\HH^3_K$ when $\dim K=3$, the hyperbolic double metric on $\mathcal{B}_{\HH}(K)$ when $\dim K=2$, and the restricted hyperbolic metric on $K\cap \HH^3_K$ when $\dim K\leq1$.

In the Klein model $\HH_K^3=\BB^3$, the Riemannian metric tensor at a point $x \in \BB^3$ is given by
$$
g_{\HH}(x)(v, v) = \frac{\|v\|^2}{1 - \|x\|^2} + \frac{\langle x, v \rangle^2}{(1 - \|x\|^2)^2},\qquad\forall v\in T_x\HH_K^3.
$$
Since $\|x\|< 1$, we have the pointwise estimate
\begin{equation}\label{eq:H-bound-E}
    g_{\HH}(x)(v, v) \ge \frac{\|v\|^2}{1 - \|x\|^2} \ge \|v\|^2,
\end{equation}
which implies that the Euclidean length of any rectifiable curve in $\BB^3$ is bounded from above by its hyperbolic length in $\HH_K^3$.

Before turning to the main convergence theorems of this subsection, we first state two auxiliary propositions. Since their proofs are lengthy and somewhat tedious, we postpone the proofs to \cref{sec:appendix-convergence} to avoid distracting from the main line of argument.

\begin{proposition}\label{prop:common-W-space}
Let $K_n \subset \overline{\BB^3}$ be a sequence of full-dimensional compact convex sets converging in the Hausdorff metric to a nonempty compact convex set $K \subset \overline{\BB^3}$. Then $(\mathcal{B}(K), d_{\mathcal{B}(K)})$ and all $(\mathcal{B}(K_n), d_{\mathcal{B}(K_n)})$ can be isometrically realized as subsets of a common compact metric space $(W, d_W)$ equipped with a $1$-Lipschitz map $\Pi : W \to \overline{\BB^3}$, such that
$$
  \Pi|_{\mathcal{B}(K_n)} = \pi_{K_n}, \qquad \Pi|_{\mathcal{B}(K)} = \pi_K,
$$
and
\begin{equation}\label{ccv:common-space}
  d_H^W(\mathcal{B}(K_n), \mathcal{B}(K)) \to 0 \quad \text{as } n \to \infty.
\end{equation}
\end{proposition}

\begin{proposition}\label{prop:ccv-finite-stability}
Suppose that sequences $p_n, q_n \in \mathcal{B}_{\HH}(K_n)$ converge in $(W, d_W)$ to $p, q \in \mathcal{B}_{\HH}(K)$, respectively. Then the hyperbolic intrinsic distances satisfy
\begin{equation}\label{ccv:finite-stability-equation}
  \delta_{\mathcal{B}(K_n)}(p_n, q_n) \to \delta_{\mathcal{B}(K)}(p, q) \quad \text{in } [0, +\infty] \quad \text{as } n \to \infty.
\end{equation}
In particular, $\delta_{\mathcal{B}(K_n)}(p_n, q_n) \to +\infty$ whenever $p$ and $q$ lie in different connected components of $\mathcal{B}_{\HH}(K)$.
\end{proposition}

We now proceed to prove the main convergence theorems of this subsection. A sequence of subsets $S_n\subset S$ is called an \emph{exhaustion} of $S$ if each compact subset of $S$ is contained in $S_n$ for all sufficiently large $n$.

\begin{lemma}\label{ccv:map-compactness}
Suppose that $K_n$, $K$, and $W$ satisfy the conditions in \cref{prop:common-W-space}. Let $(S,d)$ be a
proper metric space, let $o\in S$, and let $S_n$ be an exhaustion of $S$ with
$o\in S_n$ for every $n$. Assume that continuous maps
$F_n:S_n\to\mathcal{B}_{\HH}(K_n)$ satisfy, for each $R>0$,
\begin{equation}\label{eq:ccv-asymptotic-distortion}
    \delta_{\mathcal{B}(K_n)}(F_n(x),F_n(y))
 \leq d(x,y)+\varepsilon_n(R),\qquad \forall x,y\in A_R:=\overline B_d(o,R),
\end{equation}
where $\varepsilon_n(R)\to0$, and all points $\pi_{K_n}(F_n(o))$ are contained in a fixed compact subset of $\BB^3$. Then, after passing to a subsequence, $F_n$ converges to a continuous map $F \colon S \to \mathcal{B}_{\HH}(K)$ uniformly on compact subsets of $S$ with respect to the metric $d_W$ such that
$$
  \delta_{\mathcal{B}(K)}(F(x), F(y)) \leq d(x, y) \qquad \forall x, y \in S.
$$
\end{lemma}

\begin{proof}
Fix $R>0$. Since $\sup_n d_{\HH}(0, \pi_{K_n}(F_n(o))) < \infty$, by the upper distance bound \eqref{eq:ccv-asymptotic-distortion}, there exists $r_R\in(0,1)$ such that 
\begin{equation}\label{ccv:source-localization}
  \pi_{K_n}(F_n(A_R)) \subset \{ z \in \BB^3 : \|z\| \leq r_R \}, \qquad \forall n\in\NN.
\end{equation}
By \eqref{eq:H-bound-E}, the intrinsic Euclidean distance is always bounded by the intrinsic hyperbolic distance, which implies that for all $x, y \in A_R$,
\begin{equation}\label{ccv:asymptotic-lipschitz}
  d_W(F_n(x), F_n(y))= d_{\mathcal{B}(K_n)}(F_n(x), F_n(y)) \leq \delta_{\mathcal{B}(K_n)}(F_n(x), F_n(y)) \leq d(x, y) + \varepsilon_n(R)
\end{equation}

For any $\epsilon > 0$, choose $N \in \mathbb{N}$ such that $\varepsilon_n(R) <\epsilon / 2$ for all $n \ge N$. Then $d(x, y) < \epsilon / 2$ implies $d_W(F_n(x), F_n(y)) < \epsilon$ for all $n \ge N$. Since each of the finitely many remaining maps $F_1, \dots, F_{N-1}$ is uniformly continuous on the compact set $A_R$, the sequence $F_n$ is equicontinuous on $A_R$. Since $(W, d_W)$ is compact and $(S, d)$ is proper, the Arzel\`a--Ascoli theorem combined with a standard diagonal argument over radii $R \in \mathbb{N}$ yields a subsequence of $F_n$ converging locally uniformly to a continuous map $F \colon S \to W$. The Hausdorff convergence $d_H^W(\mathcal{B}(K_n), \mathcal{B}(K)) \to 0$ ensures that $F(S) \subset \mathcal{B}(K)$. Furthermore, the compatibility of coordinate projections $\Pi|_{\mathcal{B}(K_n)} = \pi_{K_n}$ and $\Pi|_{\mathcal{B}(K)} = \pi_K$, together with \eqref{ccv:source-localization}, guarantees that
$$
\pi_K(F(A_R)) \subset \{ z \in \BB^3 : \|z\| \leq r_R \} \subset \BB^3,
$$
which proves that $F(S) \subset \mathcal{B}_{\HH}(K)$. Finally, applying \cref{prop:ccv-finite-stability} to the convergent sequences $F_n(x) \to F(x)$ and $F_n(y) \to F(y)$ yields
$$
\delta_{\mathcal{B}(K)}(F(x), F(y)) = \lim_{n \to \infty} \delta_{\mathcal{B}(K_n)}(F_n(x), F_n(y)) \leq d(x, y).
$$
This completes the proof.
\end{proof}

\begin{theorem}\label{ccv:realization-package}
Let $(S,d)$ be a complete connected $\CBB(-1)$ surface, let $o\in S$,
and let $S_n$ be an exhaustion of $S$ with $o\in S_n$ for every $n$. Equip each $S_n$ with a metric $d_n$ which is compatible with its topology. Let $K_n\subset\overline{\BB^3}$ be full-dimensional compact convex sets converging to a compact convex set $K\subset \overline{\BB^3}$ in the Hausdorff distance and $W$ be as defined in \cref{prop:common-W-space}. Assume that
\begin{enumerate}[label=\textup{(\roman*)}]
 \item The metrics $d_n$ on $S_n$ converge to $d$ uniformly on $A\times A$ for every compact $A\subset S$;
 \item There exist continuous maps $F_n:S_n\to\mathcal{B}_{\HH}(K_n)$ such that for every compact $A\subset S$,
 $$
  \sup_{x,y\in A}\left|\delta_{\mathcal{B}(K_n)}(F_n(x),F_n(y))-d_n(x,y)\right|\to0,
 $$
 and all $\pi_{K_n}\circ F_n(o)$ are contained in a fixed compact subset of $\BB^3$.
\end{enumerate}
Then $\dim K\geq2$ and after passing to a subsequence, $F_n$ converge locally uniformly on $S$ with respect to the metric $d_W$ to an isometric embedding
$$
 F:(S,d)\to(\mathcal{B}_{\HH}(K),\delta_{\mathcal{B}(K)}).
$$
In particular, $F(S)$ is a connected component of
$\mathcal{B}_{\HH}(K)$.
\end{theorem}

\begin{proof}
By the Hopf--Rinow theorem, $S$ is proper. The conditions \textup{(i)} and \textup{(ii)} imply that on each closed ball $A_R=\overline{B}_d(o,R)$, there exists $\varepsilon_n(R)\to0$ such that
$$
 \left|\delta_{\mathcal{B}(K_n)}(F_n(x),F_n(y))-d(x,y)\right|
 \leq\varepsilon_n(R),\qquad \forall x,y\in A_R.
$$
Then, \cref{ccv:map-compactness} supplies a locally uniform limit
$F:S\to\mathcal{B}_{\HH}(K)$. By \cref{prop:ccv-finite-stability}, we have
$$
 \delta_{\mathcal{B}(K)}(F(x),F(y))
 =\lim_n\delta_{\mathcal{B}(K_n)}(F_n(x),F_n(y))
 =d(x,y).
$$
Thus $F:S\to\mathcal{B}_{\HH}(K)$ is an isometric embedding. Since $S$ cannot be embedded into an interval or a point, we must have $\dim K\geq2$. The isometry puts $F(S)$ in the component $\Sigma_0$ containing
$F(o)$, and completeness ensures the closedness of $F(S)$ in $\Sigma_0$. On the other hand, since both $S$ and $\Sigma_0$ are topologically surfaces without boundary, the invariance of domain theorem then implies that $F(S)$ is also open in $\Sigma_0$. Thus, $F(S)=\Sigma_0$ is a connected component of $\mathcal{B}_{\HH}(K)$, which completes the proof.
\end{proof}

We then turn to the case where $S$ is a surface with non-empty boundary and consider realizations that preserve the intrinsic path metrics within their images.

\begin{theorem}\label{ccv:bordered-intrinsic-realization}
Let $(S,d)$ be a compact $\CBB(-1)$ surface with nonempty
boundary. Fix $o\in S$. Let $K_n\subset\overline{\BB^3}$ be full-dimensional compact convex sets converging to a compact convex set $K\subset \overline{\BB^3}$ in the Hausdorff distance and $W$ be as defined in \cref{prop:common-W-space}. Suppose that $d_n$ are metrics 
on $S$ satisfying
$$
 \varepsilon_n:=\sup_{x,y\in S}|d_n(x,y)-d(x,y)|\to0,
$$
and $F_n:S\to\mathcal{B}_{\HH}(K_n)$ are continuous injections such that
\begin{equation}\label{eq:bordered-intrinsic-realization-1}
     \rho_n(F_n(x),F_n(y))=d_n(x,y),\qquad \forall x,y\in S,
\end{equation}
where $\rho_n$ is the intrinsic path metric on $F_n(S)$ induced by $\delta_{\mathcal{B}(K_n)}$, and all $\pi_{K_n}\circ F_n(o)$ are contained in a fixed compact subset of $\BB^3$.
Then we have $\dim K\geq2$ and after passing to a subsequence, $F_n$ converge uniformly on $S$ with respect to the metric $d_W$ to a continuous map 
$$
F:(S,d)\to(\mathcal{B}_{\HH}(K),\delta_{\mathcal{B}(K)})
$$
such that $F|_{\Int S}$ is an intrinsic isometric embedding, where $\Int S$ is equipped with the path metric induced by $d$.

\end{theorem}

\begin{proof}
Since
$$
 \delta_{\mathcal{B}(K_n)}(F_n(x),F_n(y))
 \leq \rho_n(F_n(x),F_n(y))
 = d_n(x,y)
 \leq d(x,y)+\varepsilon_n,
$$
\cref{ccv:map-compactness} yields a subsequence of $F_n$ converging uniformly in $W$ to a continuous map $F:S\to\mathcal{B}_{\HH}(K)$.
Set $\beta(x):=\operatorname{dist}_d(x,\partial S)$. By compactness each $F_n$ is a topological embedding, and by invariance of domain $F_n(\Int S)$ is open in $\mathcal{B}_{\HH}(K_n)$. By \eqref{eq:bordered-intrinsic-realization-1}, any rectifiable curve $\gamma$ in $\mathcal{B}_{\HH}(K_n)$ joining $F_n(x)$ to $F_n(y)$ that is contained in $F_n(S)$ has length at least $d_n(x,y)$. If $\gamma$ leaves $F_n(S)$, then its length in $\mathcal{B}_{\HH}(K_n)$ satisfies
$$
 l_{\delta_{\mathcal{B}(K_n)}}(\gamma)\geq \inf_{z\in\partial S}\rho_n(F_n(x),F_n(z))= \inf_{z\in\partial S}d_n(x,z) \geq \beta(x)-\varepsilon_n.
$$
Taking the infimum over all joining curves from $F_n(x)$ to $F_n(y)$ yields
$$
 \min\{d(x,y),\,\beta(x)\}-\varepsilon_n
 \leq \delta_{\mathcal{B}(K_n)}(F_n(x),F_n(y))
 \leq d(x,y)+\varepsilon_n.
$$
Passing to the limit via \cref{prop:ccv-finite-stability}, we then obtain
\begin{equation}\label{ccv:interior-distance-control}
    \min\{d(x,y),\,\beta(x)\}
 \leq \delta_{\mathcal{B}(K)}(F(x),F(y))
 \leq d(x,y),\qquad \forall x,y\in S.
\end{equation}
Since $\beta(x)>0$ for all $x\in\Int S$, \eqref{ccv:interior-distance-control} implies that $F|_{\Int S}$ is a topological embedding. In particular,  we have $\dim K\geq 2$ since a surface cannot be embedded into a point or an interval. Moreover, $F$ is a local isometry on $\Int S$. To see this, for any $a\in\Int S$, choose $0<r<\beta(a)/3$. For any points $x,y\in B_d(a,r)$, we have 
$$
d(x,y)<2r<\beta(a)-r<\beta(x).
$$
It follows from \eqref{ccv:interior-distance-control} that $F$ preserves distances on $B_d(a,r)$ and thus is a local isometry on $\Int S$. Since $F|_{\Int S}$ is both a topological embedding and a local isometry, it is clear that $F|_{\Int S}$ is an isometry with respect to the induced length metrics, which completes the proof.

\end{proof}

\section{Convex Cap Realizations}\label{sec:Convex-Cap-Realizations}
We assume throughout this section that $(S,d_0)$ is acompact genus zero $\CBB(-1)$ surface with $m$ boundary components. The boundary of $S$ is the disjoint union
$$
   \partial S=C_1\sqcup\cdots\sqcup C_m,
   \qquad m\geq 1,
$$
where each $C_i$ is a topological circle. 
\begin{definition}
\label{def:eq-capped-realization}
Let $(S,d_0)$ be acompact genus zero $\CBB(-1)$ surface with boundary 
$$
   \partial S=C_1\sqcup\cdots\sqcup C_m,
   \qquad m\geq 1.
$$
An intrinsic isometric embedding $F$ from $(S,d_0)$ to the boundary $(\partial K,d_{\partial K})$ of a compact convex body $K\subset\mathbb H^3$ is called a \emph{convex cap realization} of $(S,d_0)$ in $\HH^3$ if there exist totally geodesic planes $P_1,\ldots,P_m$ in $\HH^3$ such that
\begin{enumerate}
    \item[(1)] $E_i:=K\cap P_i$, $i\in\{1,\cdots,m\}$, are disjoint convex disks, and $F(S)\cap E_i=F(C_i)=\partial E_i$;
    \item[(2)] $\partial K=F(S)\cup E_1\cup\cdots\cup E_m$;
    \item[(3)] $\Pi_i(K)\subset E_i$ for each $i\in\{1,\cdots,m\}$, where $\Pi_i$ is the projection from $\HH^3$ onto the geodesic plane $P_i$.
\end{enumerate}
\end{definition}

Our goal is to establish the following convex cap realization theorem forcompact genus zero $\CBB(-1)$ surfaces with boundary.

\begin{theorem}\label{thm:hyperbolic-convex-cap}
    Let $(S,d_0)$ be acompact genus zero $\CBB(-1)$ surface with boundary. Then at least one of the following conclusions holds:
    \begin{enumerate}
        \item[(1)] $\Int S$ admits an intrinsic isometric embedding into the metric double $D(\Omega)$ of some closed convex domain $\Omega\subset\mathbb{H}^2$.
        \item[(2)] $S$ admits a convex cap realization in $\HH^3$.
    \end{enumerate}
\end{theorem}

When $m=1$, we glue two copies of $S$ isometrically along their boundaries. The resulting metric space $D(S)$ is a $\CBB(-1)$ surface which is homeomorphic to $\SS^2$. Now, suppose that $\Int S$ cannot be isometrically embedded in any planar metric double. Then, by Alexandrov's realization theorem, $D(S)$ can be isometrically embedded into $\mathbb{H}^3$ as the boundary of a compact convex body $K$. By Pogorelov's rigidity theorem, the reflection map on $D(S)$ exchanging two copies of $S$ extends to an isometry $F$ of $\mathbb{H}^3$. Since the boundary of $S$ is invariant under $F$, $F$ must be a hyperbolic reflection with respect to some totally geodesic plane $P$ and therefore $\partial S$ lies in the plane $P$. Moreover, reflection invariance and the convexity of $K$ imply $\Pi_P(K)\subset K\cap P$. This proves \cref{thm:hyperbolic-convex-cap} for the case where $m=1$.

When $m=2$, the boundary $\partial S$ consists of exactly two connected components, $C_1$ and $C_2$. For each $n\in\mathbb{Z}$, let $S^{(n)}$ be an isometric copy of $S$, with $C_1^{(n)}$ and $C_2^{(n)}$ denoting the corresponding copies of $C_1$ and $C_2$, respectively. For each $k\in\mathbb{Z}$, we glue $S^{(2k)}$ to $S^{(2k+1)}$ by identifying $C_1^{(2k)}$ with $C_1^{(2k+1)}$, and glue $S^{(2k)}$ to $S^{(2k-1)}$ by identifying $C_2^{(2k)}$ with $C_2^{(2k-1)}$ via the natural isometries. By \cref{cor:perelman-double} and \cref{thm:Globalization-theorem}, the resulting space $\widehat{S}$ is a $\CBB(-1)$ surface. By construction, $\widehat{S}$ is complete, planar, and has exactly two ends. 
Now suppose that $\Int S$ cannot be isometrically embedded into any planar metric double. By Alexandrov's realization theorem, $\widehat{S}$ can be isometrically realized as the boundary $\partial X$ of a noncompact convex body $X\subset\HH^3_P$. Along either end of $\widehat{S}$, the gluing curves have uniformly bounded diameter in the intrinsic metric of $\widehat{S}$. Consequently, their images on $\partial X$ have uniformly bounded hyperbolic diameter while escaping every compact subset of $\HH^3_P$. It follows that the ideal boundary $\overline{X}\cap\partial\HH^3_P$ consists of exactly two points.
Let $f_1$ and $f_2$ be the isometries of $\widehat{S}$ which are the reflections along $C_1^{(0)}$ and $C_2^{(0)}$, respectively. By the rigidity theorem in \cite{luo2026rigidity}, each $f_i$, $i\in\{1,2\}$, extends to an isometry $F_i$ of $\mathbb{H}^3$. By similar arguments as in the case $m=1$, this gives an intrinsic isometric embedding of $S^{(0)}$ into $\HH^3$ such that each boundary component $C_i^{(0)}$ of $S^{(0)}$ lies in some totally geodesic plane $P_i$ and $\Pi_i(X)\subset X\cap P_i$. This proves \cref{thm:hyperbolic-convex-cap} for the case where $m=2$.

By the above arguments we only need to consider the case $m\geq 3$. Let $T_m$ be the $m$-regular tree. For each $v\in T_m$, we can similarly define an isometric copy $S^{(v)}$ of $S$ and then glue these surfaces along their boundaries such that $S^{(v)}$ and $S^{(w)}$ are glued if and only if $v\sim w$ in $T_m$ and $S^{(v)}$ and $S^{(w)}$ are glued in at most one boundary component for all $v,w\in T_m$. In this way, we again obtain a planar $\CBB(-1)$ surface $\widehat{S}$ which can be embedded into $\HH^3$ as a convex surface. However, in this case, $\widehat{S}$ has uncountably many ends, and the current rigidity results we know are not strong enough to cover this situation. Thus, one has to find another way to prove \cref{thm:hyperbolic-convex-cap} for the case where $m\geq 3$. This is the main problem we will solve in this section.

Our strategy is as follows: We first double the surface $S$ to obtain a closed $\CBB(-1)$ surface $(\Sigma,d_{\Sigma})$ with genus $m-1$. Then by using Richard's results on smoothing Alexandrov surfaces via Ricci flow, we obtain a sequence of smooth metrics $\{\widehat{d}_i\}_i$ on $\Sigma$ with curvature strictly larger than $-1$ which is uniformly convergent to $(\Sigma,d_{\Sigma})$. Let $d_i$ be the induced path metric on $S$ in $(\Sigma,\widehat{d}_i)$. By applying a theorem of Schlenker on hyperbolic manifolds with strictly convex boundary, we show that \cref{thm:hyperbolic-convex-cap} holds if we further assume that $(S,d_0)$ is a Riemannian surface with boundary when $m\geq3$. In particular, there exists an intrinsic isometric embedding $f_i$ from $(S,d_i)$ to the boundary $\partial K_i$ of some compact convex body $K_i$ in $\HH^3$ such that each boundary component $C_i$ of $S$ is embedded into some totally geodesic plane in $\HH^3$. Finally, under the assumption that $(\Int S,d_{\Int S})$ cannot be isometrically embedded in any planar metric double surface, we show that $K_i$ converges to some compact convex body $K$ in $\HH^3$ and $f_i$ converges to an intrinsic isometric embedding of $(S,d_0)$ into $\partial K$, which completes the proof.

\subsection{Equivariant construction of convex cap realizations}
\label{sec:eq-construction}

The main result in this section is the following:
\begin{theorem}\label{thm:eq-capped-realization}
Let $(S,d_0)$ be acompact genus zero $\CBB(-1)$ surface with $m$ boundary components. If $m\ge3$, then there exists a sequence of smooth path metrics $\{d_j\}_{j=1}^{\infty}$ on $S$ such that
$$
  \sup_{x,y\in S}|d_j(x,y)-d_0(x,y)|\to 0\quad\text{as}\quad j\to\infty
$$
and there exists a convex cap realization $F_j$ for each $(S,d_j)$.
\end{theorem}

\subsubsection{Equivariant Ricci smoothing}
    We first introduce the concept of Ricci flow initiated from metric spaces.
    \begin{definition}
        A Ricci flow $(M, g(t))_{t \in (0,T)}$ on a compact manifold $M$ is said to have the metric space $(X, d)$ as metric initial condition if the Riemannian distances $d_{g(t)}$ uniformly converge as $t$ goes to $0$ (as functions $M \times M \to \mathbb{R}$) to a distance $d$ on $M$.
    \end{definition} 
    For closed $\CBB(-1)$ surfaces, the following existence and uniqueness result of Ricci flow is known; see \cite[Section~1]{Richard2018}, where the uniqueness result we use is Proposition 0.6.
    \begin{proposition}\label{prop:Ricci-flow}
        Let $(\Sigma,d)$ be a closed $\CBB(-1)$ surface. Then 
        \begin{enumerate}
            \item[(1)] There exists $T>0$ and a Ricci flow $(\Sigma, g(t))_{t \in (0,T)}$ having $(\Sigma, d)$ as the metric initial condition such that $K_{g(t)}(x)\geq-1$ for every $(x,t) \in \Sigma \times (0,T)$, where $K_{g(t)}$ is the Gauss curvature of $(\Sigma,g(t))$ at $x$.
            \item[(2)] Suppose that $\{g_1(t)\}_{t \in (0,T)}$ and $\{g_2(t)\}_{t \in (0,T)}$ are two Ricci flows on $\Sigma$ having $(\Sigma, d)$ as the metric initial condition, and each of them is smooth with respect to some differential structure on $\Sigma$. Assume furthermore that 
            $$ 
            K_{g_i(t)}(x) \geq -1,\quad \forall (x,t) \in \Sigma \times (0,T).
            $$
            Then the two smooth structures on $\Sigma$ coincide and $g_1(t) = g_2(t)$ for all $t \in (0,T)$.
        \end{enumerate}
    \end{proposition} 

    Let $\Sigma=D(S)$ be the metric double of $S$. The Euler characteristic identity of $\Sigma$ is
    $$
     \chi(\Sigma)=2\chi(S)=4-2m.
    $$
    This shows that $\Sigma$ is a closed orientable surface of genus $g=m-1$. By \cref{cor:perelman-double}, $(\Sigma,d_{\Sigma})$ is a $\CBB(-1)$ surface. The map $\iota$ interchanging the two copies of $S$ is an
    isometric involution and
    $$
      \operatorname{Fix}(\iota)=C_1\sqcup\cdots\sqcup C_m.
    $$
    We are now ready to prove the following theorem:
    \begin{theorem}
    \label{thm:eq-richard-smoothing}
    There exists a
    sequence of smooth metrics $\{h_j\}_{j=1}^{\infty}$ on $\Sigma$
    such that for each $j\in\mathbb{N}$,
    \begin{enumerate}
    \item[(1)] $K_{h_j}(x)>-1$ for all $x\in \Sigma$;
    \item[(2)] $\iota$ is a smooth isometric involution of $(\Sigma,h_j)$;
    \item[(3)] the metrics $h_j$ converge to $d_{\Sigma}$ uniformly, i.e.,
    $$
      \sup_{x,y\in \Sigma}|d_{h_j}(x,y)-d_\Sigma(x,y)|\to 0,\quad\text{as} \quad j\to\infty.
    $$
    \end{enumerate}
    \end{theorem}
    \begin{proof}
        By \cref{prop:Ricci-flow} (1), there exists $T>0$ and a Ricci flow $(\Sigma, g(t))_{t \in (0,T)}$ having $(\Sigma, d_{\Sigma})$ as the metric initial condition such that $K_{g(t)}(x)\geq-1$ for every $(x,t) \in \Sigma \times (0,T)$. Pull the flow back by the topological involution
        $\iota$. Clearly, the resulting flow $\{\iota^*g(t)\}_{t\in(0,T)}$ is smooth with respect to the pulled-back differential structure. Since $\iota$ is an isometric involution of $(\Sigma,d_{\Sigma})$, we see that $\{\iota^*g(t)\}_{t\in(0,T)}$ is also a Ricci flow  having $(\Sigma, d_{\Sigma})$ as the metric initial condition and $K_{\iota^*g(t)}(x)\geq-1$ for every $(x,t) \in \Sigma \times (0,T)$. It then follows from \cref{prop:Ricci-flow} (2) that $\iota^*g(t)=g(t)$ for every $t\in(0,T)$. Thus, $\iota$ is a smooth isometric involution of $(\Sigma,g(t))$ for every $t\in(0,T)$. Now choose a sequence $t_j\to0^+$ and a sequence $\lambda_j\to1^+$. Set $h_j=\lambda_j g(t_j)$. Then $K_{h_j}(x)\geq-\frac{1}{\lambda_j}>-1$ for all $x\in \Sigma$, $\iota$ is a smooth isometric involution of $(\Sigma,h_j)$ since $\iota$ is a smooth isometric involution of $(\Sigma,g(t_j))$, and the metrics $h_j$ converge to $d_{\Sigma}$ uniformly since $\lambda_j\to 1$ and the metrics $g(t_j)$ converge to $d_{\Sigma}$ uniformly as $j\to\infty$. This completes the proof.
    \end{proof}
    Now let $\Sigma$ be the metric double of $S$ and let $S^+,S^-$ be the two copies of $S$ in $\Sigma$. Let $d_j$ be the induced path metric of $(\Sigma,h_j)$ on $S^+$. The following corollary is immediate: 
    \begin{corollary}
    \label{cor:eq-half-convergence}
    For all $x,y\in S^+$, $d_j(x,y)=d_{h_j}(x,y)$ and
    $$
     \sup_{x,y\in S}|d_j(x,y)-d_0(x,y)|\to0\quad\text{as} \quad j\to\infty.
    $$
    \end{corollary}

\subsubsection{A smooth equivariant handlebody model}

\begin{definition}
Let $M$ be a smooth manifold. We denote by
$$
\operatorname{Diff}_0(M)
$$
the identity component of the diffeomorphism group
$\operatorname{Diff}(M)$; equivalently,
$\operatorname{Diff}_0(M)$ consists of those diffeomorphisms
$\phi:M\to M$ that are smoothly isotopic to the identity.
Thus, two Riemannian metrics $g_1$ and $g_2$ on $M$ are said to be
equivalent up to isotopy if there exists
$\phi\in\operatorname{Diff}_0(M)$ such that
$$
\phi^*g_2=g_1.
$$
The same convention applies to metrics on $\partial M$.
\end{definition}

\begin{theorem}[(Schlenker {\cite[Theorem~0.1 and
Definition~0.5]{MR2208419}})]
\label{thm:eq-schlenker}
Let $M$ be a compact connected $3$-manifold whose interior admits a
complete convex co-compact hyperbolic structure. Define
$$
\mathcal{G}(M):=\left.\left\{
\begin{array}{c}
\text{hyperbolic metrics $g$ on $M$ such that $\partial M$ is}\\
C^\infty\text{-smooth and strictly convex}
\end{array}
\right\}\right/\operatorname{Diff}_0(M),
$$
where $\operatorname{Diff}_0(M)$ acts by pullback, and define
$$
\mathcal H(\partial M):=\left.
\left\{  C^\infty \text{ Riemannian metrics }
h \text{ on }  \partial M \text{ such that } K_h>-1
\right\}
\right/\operatorname{Diff}_0(\partial M).
$$
Then the map
$$
F:\mathcal G(M)\to\mathcal H(\partial M),
\qquad
[g]\longmapsto[I_g],
$$
where $I_g$ denotes the metric induced by $g$ on $\partial M$, is a
homeomorphism.

In particular, every $C^\infty$ Riemannian metric $h$ on $\partial M$
with $K_h>-1$ is realized, up to isotopy, as the induced boundary metric of a hyperbolic metric on $M$ with smooth strictly convex boundary, and this
hyperbolic metric is unique up to isotopy.
\end{theorem}

Now we consider applying the above theorem of Schlenker to handlebodies. For a genus $g$ handlebody $H_g$, where $g=m-1$, we see that $\partial H_g$ is homeomorphic to $\Sigma$. Hence we may identify $\partial H_g$ with $\Sigma$. We first show that there exists a smooth involution on $H_g$ that is compatible with the involution $\iota$ on $\Sigma$.
\begin{lemma}
\label{lem:eq-smooth-handlebody}
Let $g=m-1$. Then there exists a smooth genus-$g$ handlebody $H_g$, a smooth
involution $\rho:H_g\to H_g$, and a smooth boundary identification
$$
 b:(\partial H_g,\rho|_{\partial H_g})\to(\Sigma,\iota)
$$
such that $b(\rho(x))= \iota(b(x))$ for all $x\in \partial H_g$.
Moreover,
$$
 \rho_*=-\operatorname{id}
 \quad\text{on }H_1(H_g;\mathbb R).
$$
\end{lemma}
\begin{proof}
Let $B^+$ and $B^-$ be two copies of the $3$-ball. Choose pairwise
disjoint closed disks
$$
D_1^\pm,\ldots,D_m^\pm\subset\partial B^\pm,
$$
and identify $D_i^+$ with $D_i^-$ by the identity map for every $i$. Let $H_g$ be the resulting quotient space. After an equivariant smoothing of the corners, $H_g$ is a smooth handlebody of genus $g=m-1$. The map interchanging $B^+$ and $B^-$ is then a smooth involution $\rho:H_g\to H_g$. 
It is clear that the boundary $\partial H_g$ is the double of a planar surface with $m$
boundary components, and $\rho|_{\partial H_g}$ is the involution
exchanging its two copies. By the corresponding description of
$(\Sigma,\iota)$, there is an equivariant identification
$$
 b:(\partial H_g,\rho|_{\partial H_g})\to(\Sigma,\iota).
$$
Thus, $b(\rho(x))= \iota(b(x))$ for all $x\in \partial H_g$. 

It remains to determine the action of $\rho$ on $H_1(H_g;\mathbb R)$.
The handlebody $H_g$ admits a $\rho$-equivariant deformation retraction
onto the graph $\Gamma$ consisting of two vertices joined by $m$ edges
$e_1,\ldots,e_m$. Orient all the edges from one vertex to the other.
Since $\rho$ exchanges the two vertices and reverses each edge, it acts
as $-\operatorname{id}$ on $C_1(\Gamma;\mathbb R)$, and hence also on
$H_1(\Gamma;\mathbb R)$. Therefore, we have
$$
\rho_*=-\operatorname{id}
\qquad\text{on }H_1(H_g;\mathbb R).
$$
\end{proof}

\begin{proposition}
\label{prop:eq-marked-realization}
Assume $g\ge2$.  For every $j$ there is a hyperbolic metric $G_j$ on the fixed smooth
handlebody $H_g$ whose boundary is smooth and strictly convex, and the induced metric of $G_j$ on $\partial H_g$ is exactly $h_j$ under the fixed boundary identification.
\end{proposition}
\begin{proof}
By \cref{thm:eq-schlenker}, there exists a hyperbolic metric $G'_j$ on $H_g$ and a
boundary diffeomorphism $a\in\operatorname{Diff}_0(\partial H_g)$ such
that $a^*I_{G'_j}=h_j$. We claim that $a$ can be extended to a diffeomorphism $A$ of $H_g$. If this is true, then the proof is completed by setting $G_j:=A^*G'_j$. 

We now prove the claim. Choose a smooth isotopy $a_{t}$, $0\leq t\leq1$, from
$\operatorname{id}_{\partial H_g}$ to $a$, and a smooth collar embedding
$$
\Phi:\partial H_g\times[0,2\varepsilon)\to H_g
$$
with $\Phi(x,0)=x$ for every $x\in\partial H_g$.
Let $\chi:[0,2\varepsilon)\to[0,1]$ be a smooth function with $\chi(0)=1$ and
$\chi=0$ on $[\varepsilon,2\varepsilon)$. Define
$$
A\bigl(\Phi(x,s)\bigr)
=
\Phi\bigl(a_{\chi(s)}(x),s\bigr)
$$
on the collar $\Phi(\partial H_g\times[0,2\varepsilon))$ and let $A=\operatorname{id}_{H_g}$ outside $\Phi(\partial H_g\times[0,2\varepsilon))$. By the choice of $\chi$, we see that $A$ is well-defined and is a diffeomorphism of $H_g$ with $A|_{\partial H_g}=a$. This completes the proof.
\end{proof}

\begin{proposition}
\label{prop:eq-isometric-involution}
Assume $g\ge2$.  For every $j$ there is a $G_j$-isometric involution
$$
 R_j:H_g\to H_g
$$
such that
$$
  R_j|_{\partial H_g}=\iota,
  \qquad R_j\simeq\rho.
$$
\end{proposition}
\begin{proof}
    By \cref{lem:eq-smooth-handlebody}, $\rho$ is a smooth involution. Since $\rho|_{\partial H_g}=\iota$ and $\iota^*h_j=h_j$, we see that the induced boundary metric $I_{\rho^*G_j}=I_{G_j}=h_j$. Then, by \cref{thm:eq-schlenker} there exists a diffeomorphism $\psi_j\in\operatorname{Diff}_0( H_g)$ such that $\psi_j^*(\rho^*G_j)=G_j$. Set $R_j=\rho\circ\psi_j$. Then $R_j$ is a $G_j$-isometry and $R_j$ is isotopic to $\rho$. It follows that the induced boundary map $R_j|_{\partial H_g}$ is an $h_j$-isometry isotopic to $\iota$ and consequently $\iota^{-1}\circ(R_j|_{\partial H_g})$ is an $h_j$-isometry isotopic to the identity. 
    By uniformization, the conformal class of $h_j$ contains a unique hyperbolic metric. Every $h_j$-isometry preserves this conformal class and hence, by uniqueness, is also an isometry with respect to the hyperbolic metric. Since $\partial H_g$ is a closed surface of genus at least two, any such isometry isotopic to the identity must be the identity itself (see \cite[Proposition~7.7]{FarbMargalit}). Thus, $\iota^{-1}\circ(R_j|_{\partial H_g})=\operatorname{id}_{\partial H_g}$ and $R_j|_{\partial H_g}=\iota$.
    It remains to check that $R_j$ is an involution. Since $R_j|_{\partial H_g}=\iota$ and $\iota$ is an involution, we see that $R_j^2|_{\partial H_g}=\iota^2=\operatorname{id}_{\partial H_g}$. Then $R_j^2=\operatorname{id}_{H_g}$ follows from the standard fact that any isometry of a connected Riemannian manifold fixing the
    boundary pointwise must be the identity.
\end{proof}

We now determine the fixed-point set of the involution on the handlebody.

\begin{theorem}\label{thm:eq-fixed-disks}
Let $R:H_g\to H_g$ be an isometric involution such that
$$
R|_{\partial H_g}=\iota,
\qquad
R\simeq\rho.
$$
Then
$$
\operatorname{Fix}(R)
=
E_1\sqcup\cdots\sqcup E_m,
\qquad
\partial E_i=C_i,
$$
where $E_1,\ldots,E_m$ are pairwise disjoint properly embedded totally
geodesic disks.
Moreover,
$$
H_g\setminus\operatorname{Fix}(R)
=
B^+\sqcup B^-,
$$
where each component is homeomorphic to a $3$-ball with $m$ pairwise
disjoint closed disks removed from its boundary, and the involution $R$ exchanges $B^+$ and $B^-$.
\end{theorem}

\begin{proof}
First, note that $R$ is orientation reversing since $R\simeq\rho$ and $\rho$ is orientation reversing. Let
$p\in\operatorname{Fix}(R)$. Then the differential $dR_p: T_p H_g\to T_pH_g$ is an involution with $\det (dR_p)=-1$. In normal coordinates centered at $p$, we have
$$
R(\exp_p v)=\exp_p(dR_pv).
$$
Thus, $\operatorname{Fix}(R)$ is locally the exponential image of the
$+1$-eigenspace of $dR$. Since $\det (dR_p)=-1$, the dimension of the $+1$-eigenspace of $dR_p$ is either $0$ or $2$. It follows that each component of $\operatorname{Fix}(R)$ is either an isolated point or a totally geodesic surface (possibly with boundary). In particular, each $C_i\subset \operatorname{Fix}(\iota)\subset\operatorname{Fix}(R)$ lies in the boundary of some totally
geodesic surface component of $\operatorname{Fix}(R)$.

Now, let $E_1,\ldots,E_c$ be the components of $\operatorname{Fix}(R)$ having non-empty intersection with
$\partial H_g$, and let $l_\alpha$ be the number of boundary components
of $E_\alpha$. Since
$\operatorname{Fix}(R)\cap\partial H_g
=C_1\sqcup\cdots\sqcup C_m$, we have
$\sum_{\alpha=1}^c l_\alpha=m$. Consider any connected compact surface $E$ with $l \ge 1$ boundary components. It is clear that the $\FF_2$-homology groups of $E$ satisfy $H_0(E;\mathbb F_2)\cong\mathbb F_2$ and $H_2(E;\mathbb F_2)=0$. Consequently, the Euler characteristic is $\chi(E) = 1 - b_1(E; \mathbb{F}_2)$. Let $b_k(E;\mathbb F_2):=\dim_{\FF_2} H_k(E;\FF_2)$ be the $k$-th $\FF_2$-Betti number of $E$. The total $\FF_2$-Betti number of $E$ is then given by
$$
\beta(E;\FF_2):=\sum_k b_k(E;\mathbb F_2)
=1+b_1(E;\mathbb F_2)=2-\chi(E).
$$
By the classification of surfaces, $\chi(E) = 2 - 2h - l$ if $E$ is orientable of genus $h \ge 0$, and $\chi(E) = 2 - h - l$ if $E$ is non-orientable of genus $h \ge 1$. In either case, $\chi(E) \leq 2 - l$. Therefore, we conclude that $\beta(E;\FF_2)\geq l$ with equality if and only if $E$ is orientable and has genus zero. Since each $E_{\alpha}$ is a connected component of $\operatorname{Fix}(R)$, it follows that
\begin{equation}\label{eq:m-less-than}
    m=\sum_{\alpha=1}^c l_\alpha\leq\sum_{\alpha=1}^c \beta(E_\alpha;\FF_2)\leq\beta(\operatorname{Fix}(R);\FF_2).
\end{equation}
On the other hand, Smith's inequality (see, e.g., \cite[Theorem~4.4]{Floyd1952}) implies that
\begin{equation}\label{eq:m-greater-than}
\beta
\bigl(\operatorname{Fix}(R);\FF_2\bigr)
\leq
\beta(H_g;\FF_2)=1+g=m,
\end{equation}
where $\beta(H_g;\mathbb{F}_2)=1+g$ follows from the fact that a genus-$g$ handlebody deformation retracts onto a wedge of $g$ circles. Consequently, both inequalities in \eqref{eq:m-less-than} and \eqref{eq:m-greater-than} must actually be equalities. We then conclude that $\beta(E_\alpha;\mathbb{F}_2)=l_{\alpha}$ for all $\alpha\in\{1,\ldots,c\}$, which implies that every $E_\alpha$ is orientable and of genus zero. Moreover, the fact that $\sum_{\alpha=1}^c \beta(E_\alpha;\mathbb{F}_2) = \beta(\operatorname{Fix}(R);\mathbb{F}_2)$ forces $\operatorname{Fix}(R)=E_1\sqcup\cdots\sqcup E_c$.

We next show that each $E_\alpha$ is a closed disk. Since $R\simeq\rho$, by \cref{lem:eq-smooth-handlebody},
$$
R_*=\rho_*=-\operatorname{id}
\quad\text{on }H_1(H_g;\mathbb R).
$$
Hence, the Lefschetz number
$$
L(R)=1-\operatorname{tr}\bigl(R_*|_{H_1(H_g;\mathbb R)}\bigr)=1+g=m.
$$
By \cite{Illman1978}, there exists an equivariant triangulation of $H_g$ with respect to $R$. After passing to its
barycentric subdivision if necessary, we may then assume that
$\operatorname{Fix}(R)$ is a subcomplex of $H_g$ and that every simplex preserved by $R$ is fixed by $R$ pointwise. The Hopf trace formula then implies that $L(R)=\chi(\operatorname{Fix}(R))$.
Since each $E_\alpha$ is of genus zero,
$$
\chi(\operatorname{Fix}(R))
=\sum_{\alpha=1}^c(2-l_{\alpha})
=2c-m.
$$
It follows that $2c-m=m$ and hence $c=m$. Thus, each $E_\alpha$ has exactly one boundary component and is
therefore a closed disk. After relabeling, we obtain
$$
\operatorname{Fix}(R)
=
E_1\sqcup\cdots\sqcup E_m,
\qquad
\partial E_i=C_i.
$$

It remains to identify the complementary components. Cut $H_g$ along $E_1\sqcup\cdots\sqcup E_m$ and denote the resulting compact $3$-manifold by $\widehat H$. Since cutting a
handlebody along a properly embedded disk produces one or two handlebodies, we see that every connected component of $\widehat H$ is a handlebody. For each $i$, cutting along $E_i$ produces two copies $E_i^+$ and $E_i^-$.  Since $\partial H_g\setminus\bigcup_{i=1}^m C_i$ is the disjoint union of the two planar surfaces $S^+$ and $S^-$, whose closures $\overline{S^\pm}$ have exactly $m$ boundary components,
we have $\partial\widehat H=\Sigma^+\sqcup\Sigma^-$, where the surfaces
$$
\Sigma^\pm
=
\overline{S^\pm}
\cup
E_1^\pm\cup\cdots\cup E_m^\pm
$$
are $2$-spheres.
Note that every handlebody has exactly one boundary component. Hence $\Sigma^+$ and $\Sigma^-$
belong to distinct components of $\widehat H$, and $\widehat H$ has exactly two connected components. Since every handlebody
with boundary $S^2$ is a $3$-ball, we see that
$H_g\setminus\operatorname{Fix}(R)=B^+\sqcup B^-$, where $\overline{B^\pm}\cong B^3$ and
$\partial\overline{B^\pm}=\Sigma^\pm$.
Finally, since $R|_{\partial H_g}=\iota$ exchanges $S^+$ and $S^-$, the
involution $R$ exchanges $B^+$ and $B^-$.

\end{proof}

\subsubsection{Convex cap realizations for $(S,d_j)$}
\label{sec:proof-high-genus-rewrite}

\begin{proposition}
\label{prop:eq-convex-development}
Let $(H_g,G_j)$ be the hyperbolic handlebody obtained in
\cref{prop:eq-marked-realization}, and equip its universal
cover $\widetilde H_g$ with the lifted metric $\widetilde G_j$.
Then the developing map
$$
 \operatorname{dev}_j:
 (\widetilde H_g,\widetilde G_j)\longrightarrow\HH^3
$$
is an isometry onto a closed convex subset $X_j\subset\HH^3$
with nonempty interior. Moreover, any two developing maps for $G_j$
differ by post-composition with an isometry of $\HH^3$.
\end{proposition}

\begin{proof}
Since $H_g$ is compact with smooth, strictly convex boundary, 
$(H_g,G_j)$ is locally convex and complete with respect to its 
induced length metric. By Proposition~8.3.2 and the subsequent
local convexity criterion in \cite[Section~8.3]{ThurstonGT3M},
$\operatorname{dev}_j$ is a homeomorphism onto a convex subset 
$X_j\subset\HH^3$. Since $\operatorname{dev}_j$ is a local isometry 
onto a convex subset of $\HH^3$, it is a global isometry. 
The completeness of $\widetilde G_j$ implies that $X_j$ is closed, 
and the fact that $\operatorname{dev}_j$ is a local homeomorphism at 
interior points ensures that $X_j$ has nonempty interior. 
Finally, the uniqueness up to isometry follows from the general 
construction of developing maps; see \cite[Section~3.5]{ThurstonGT3M}.
\end{proof}



\begin{lemma}\label{lem:eq-complete-wall} 
    Let $\widetilde E$ be a lift of one of the totally geodesic disks $\{E_{i,j}\}_{i=1}^m$ in $\operatorname{Fix}(R_j)$ obtained in \cref{thm:eq-fixed-disks}. Then there is a unique lift $\widetilde R:\widetilde{H_g}\to \widetilde{H_g}$ of $R_j\circ\pi$ fixing $\widetilde E$ pointwise, where $\pi:\widetilde{H_g}\to H_g$ is the universal covering map. 
    Moreover, let $P\subset\mathbb H^3$ be the totally geodesic plane containing $\operatorname{dev}_j(\widetilde E)$ and let $\sigma_P$ be the hyperbolic reflection with respect to $P$; we have
    $$ 
    \operatorname{dev}_j\circ\widetilde R = \sigma_P\circ\operatorname{dev}_j \quad\text{and}\quad X_j\cap P=\operatorname{dev}_j(\widetilde E).
    $$
    
\end{lemma} \begin{proof} 

Fix $\widetilde x\in\widetilde E$. Since $\widetilde{H_g}$ is simply connected, the theorem implies that there is a unique lift $\widetilde R$ of $R_j\circ\pi$ such that $\widetilde R(\widetilde x)=\widetilde x$.
Since $R_j$ fixes $E_{i,j}$ pointwise, uniqueness of path lifting shows that $\widetilde R$ fixes $\widetilde E$ pointwise. 

Consider the map 
$$
\operatorname{dev}_j\circ\widetilde R\circ\operatorname{dev}_j^{-1}: X_j\to X_j.
$$
It is a hyperbolic local isometry on the connected open set $\operatorname{Int}X_j$. Thus, $\operatorname{dev}_j\circ\widetilde R\circ\operatorname{dev}_j^{-1}$ extends uniquely to an isometry $A$ of $\mathbb H^3$. By \cref{prop:eq-isometric-involution}, $R_j$ is isotopic to $\rho$ and therefore is orientation-reversing. It follows that $\widetilde R$ and $\operatorname{dev}_j\circ\widetilde R\circ\operatorname{dev}_j^{-1}$ are also orientation-reversing, and hence so
is the global isometry $A$. Moreover, since $\widetilde R$ fixes
$\widetilde E$ pointwise, we see that $A$ fixes the two-dimensional set $\operatorname{dev}_j(\widetilde E)\subset P$ pointwise. Thus, we must have $A=\sigma_P$, and hence
$$ 
\operatorname{dev}_j\circ\widetilde R(x) = \sigma_P\circ\operatorname{dev}_j(x),\quad\forall x\in\widetilde{H_g}.
$$
This also implies that
\begin{equation}\label{eq:Xj-cap-P}
    X_j\cap P = \operatorname{Fix}(\sigma_P|_{X_j}) =\operatorname{dev}_j(\operatorname{Fix}\widetilde R).
\end{equation}

It remains to show that $X_j\cap P=\operatorname{dev}_j(\widetilde{E})$, or equivalently, $\operatorname{Fix}\widetilde R=\widetilde{E}$, since the inverse image of $X_j\cap P$ under $\operatorname{dev}_j$ is exactly $\operatorname{Fix}\widetilde R$ by \eqref{eq:Xj-cap-P}. Clearly, the set $X_j\cap P$ is convex and therefore connected. Thus, its inverse image $\operatorname{Fix}\widetilde R$ under $\operatorname{dev}_j$ is a connected subset. Note that
$$ 
\operatorname{Fix}\widetilde R \subset \pi^{-1}(\operatorname{Fix}R_j) = \pi^{-1}\left(\bigsqcup_{i=1}^m E_{i,j}\right).
$$ 
Since the sets $\{E_{i,j}\}_{i=1}^m$ are mutually disjoint and simply connected, the connected components of $\pi^{-1}\left(\bigsqcup_{i=1}^m E_{i,j}\right)$ are precisely the lifts of the disks $E_{k,j}$. It follows that $\operatorname{Fix}\widetilde R$ is contained in $\widetilde E$. The reverse inclusion is immediate, and hence $X_j\cap P=\operatorname{dev}_j(\widetilde E)$. This completes the proof.
\end{proof} 

We now construct a convex cap realization for $(S,d_j)$. Let $B_j^+$ be a component of $H_g\setminus\operatorname{Fix}(R_j)$ and let $\pi:\widetilde H_g\to H_g$ be the universal covering map. Since $\overline{B_j^+}$ is a $3$-ball, every component of $\pi^{-1}(\overline{B_j^+})$ maps homeomorphically onto $\overline{B_j^+}$. Fix one such component and denote it by $\widetilde B_j^+$. Then, $\partial\widetilde B_j^+$ is the union of an isometric copy $\widetilde S_j^+$ of $S$ and $m$ mutually disjoint disks $\widetilde E_{i,j}$, where for each $i\in\{1,\cdots,m\}$, $\widetilde E_{i,j}$ is the unique lift of $E_{i,j}$ contained in $\partial\widetilde B_j^+$. Let $P_{i,j}$ be the totally geodesic plane in $\HH^3$ containing $\operatorname{dev}_j(\widetilde E_{i,j})$. By \cref{lem:eq-complete-wall} we have 
$$
\operatorname{dev}_j(\widetilde B_j^+)\cap P_{i,j}=X_j\cap P_{i,j} =\operatorname{dev}_j(\widetilde E_{i,j}).
$$
It follows that $\operatorname{dev}_j(\Int \widetilde B_j^+)\cap P_{i,j}=\varnothing$ and therefore $\operatorname{dev}_j(\widetilde B_j^+)$ lies in some closed half-space $H_{i,j}^+$ bounded by $P_{i,j}$.

\begin{proposition}\label{prop:eq-halfspace-chamber}
One has
$$
\operatorname{dev}_j(\widetilde B_j^+)
=
X_j\cap\bigcap_{i=1}^m H_{i,j}^+.
$$
In particular, $\operatorname{dev}_j(\widetilde B_j^+)$ is a full-dimensional compact convex body such that for each $i$, the projection of $\operatorname{dev}_j(\widetilde B_j^+)$ to the geodesic plane $P_{i,j}$ is contained in $\operatorname{dev}_j(\widetilde E_{i,j})$.
\end{proposition}

\begin{proof}
Set $K_j=X_j\cap\bigcap_{i=1}^m H_{i,j}^+$. By definition of $H_{i,j}^+$, for each $i\in\{1,\cdots,m\}$ we have $\operatorname{dev}_j(\widetilde B_j^+) \subset H_{i,j}^+$. Therefore, 
$\operatorname{dev}_j(\widetilde B_j^+)\subset K_j$. For the reverse inclusion, let $x\in K_j$ and $y\in \operatorname{dev}_j(\Int\widetilde B_j^+)$.
Since $K_j$ is convex, the hyperbolic geodesic segment $\gamma=[y,x]$ is contained in $K_j$.
Suppose that $x\notin\operatorname{dev}_j(\widetilde B_j^+)$. Let $q$ be the first point at which
$\gamma$ leaves $\operatorname{dev}_j(\widetilde B_j^+)$. Since $y\in\operatorname{Int}X_j$ and
$X_j$ is a convex set with nonempty interior, every point of the segment
$[y,x)$ lies in $\operatorname{Int}X_j$. In particular, $q\in\operatorname{Int}X_j$. Hence $q$ cannot lie on the developed outer boundary
$\operatorname{dev}_j(\widetilde S_j^+)\subset\partial X_j$. It follows that $q$ lies in the interior of some fixed disk $\operatorname{dev}_j(\widetilde E_{i,j})$. Since $y\notin P_{i,j}$, we see that the geodesic segment $[y,x]$ must cross the plane $P_{i,j}$ and the geodesic segment $(q,x]$ lies in the open half-space opposite to $H_{i,j}^+$, which is a contradiction since $(q,x]\subset K_j\subset H_{i,j}^+$. Finally, note that $X_j$ is invariant under the reflection with respect to each geodesic plane $P_{i,j}$. This implies that the projection of $\operatorname{dev}_j(\widetilde B_j^+)$ to the geodesic plane $P_{i,j}$ is contained in $\operatorname{dev}_j(\widetilde E_{i,j})=X_j\cap P_{i,j}$.
This completes the proof.
\end{proof}

\begin{corollary}\label{cor:eq-boundary-decomposition} 
    For each $j\in\mathbb{N}$, there exists a compact convex body $K_j$ in $\HH^3$ and an intrinsic isometric embedding $F_j$ from $(S,d_j)$ into $(\partial K_j,d_{\partial K_j})$ such that
    $$ 
    \partial K_j = F_j(S)\cup B_{1,j}\cup\cdots\cup B_{m,j}, 
    $$ 
    where $\{B_{i,j}\}_{i=1}^m$ are mutually disjoint compact convex disks, each of which is the intersection of some totally geodesic plane in $\HH^3$ with $K_j$, and $F_j(S)\cap B_{i,j}=\partial B_{i,j}$ for each $i\in\{1,\cdots,m\}$.
\end{corollary}
    
\begin{proof} 

    Set $K_j=\operatorname{dev}_j(\widetilde B_j^+)$, $S_j=\operatorname{dev}_j(\widetilde S_j^+)$, and $B_{i,j} = \operatorname{dev}_j(\widetilde E_{i,j})$ for each $i\in\{1,\cdots,m\}$. By the previous argument, we see that $\{B_{i,j}\}_{i=1}^m$ are mutually disjoint compact convex disks, each of which is the intersection of the totally geodesic plane $P_{i,j}\subset\HH^3$ with $K_j$, and $S_j\cap B_{i,j}=\partial B_{i,j}$ for each $i\in\{1,\cdots,m\}$. The decomposition follows from the fact that
    $$ 
    \partial\widetilde B_j^+ = \widetilde S_j^+ \cup \widetilde E_{1,j} \cup\cdots\cup \widetilde E_{m,j}. 
    $$
    Finally, the boundary identification $b^{-1}$ sends $(S,d_j)$ isometrically onto the corresponding half-boundary $S^+$ of $(H_g,G_j)$. Its lift to $\widetilde S_j^+$, followed by $\operatorname{dev}_j$ gives an isometry $F_j:(S,d_j)\to (S_j,d_{S_j})$, where $d_{S_j}$ is the induced path metric on $S_j$ in $\HH^3$. This completes the proof.
    \end{proof}

    \begin{proof}[Proof of \cref{thm:eq-capped-realization}]
        \cref{cor:eq-half-convergence} gives the uniform convergence
        $d_j\to d_0$. \cref{prop:eq-halfspace-chamber} and \cref{cor:eq-boundary-decomposition} ensure all the remaining conditions in \cref{def:eq-capped-realization}.
    \end{proof}

\subsection{Proof of the high-genus strong existence theorem}

In this section, we complete the proof of \cref{thm:hyperbolic-convex-cap} by establishing the remaining case $m \ge 3$, which we restate as follows.

\begin{theorem}
\label{thm:final-high-genus}
Let $(S,d_0)$ be acompact genus zero $\CBB(-1)$ surface with $m\ge 3$ boundary components, and equip $\Int S$ with the path metric $d_{\Int S}$ induced by $d_0$. Assume that $(\Int S, d_{\Int S})$ admits no intrinsic isometric embedding into $D(\Omega)$ for any closed convex domain $\Omega\subset\HH^2$. Then $(S,d_0)$ admits a convex cap realization in $\HH^3$.
\end{theorem}
\begin{proof}[Proof of \cref{thm:hyperbolic-convex-cap} for $m\geq 3$]
    Let $\partial S=\bigsqcup_{i=1}^m \gamma_i$ be the disjoint union of $m\ge 3$ Jordan curves. By \cref{thm:eq-capped-realization}, there exists a sequence of metrics $\{d_j\}_{j=1}^{\infty}$ on $S$ satisfying
    $$
    \varepsilon_j:=\sup_{x,y\in S}|d_j(x,y)-d_0(x,y)|\to 0 \quad \text{as } j\to\infty,
    $$
    such that each $(S,d_j)$ admits a convex cap realization $F_j:(S,d_j)\to (\partial K_j, d_{\partial K_j})$ in the sense of \cref{def:eq-capped-realization}. Specifically, $K_j\subset\HH^3$ is a compact convex body, $F_j$ is an isometry onto $F_j(S)$ equipped with its induced path metric, and there exist supporting totally geodesic planes $P_{1,j},\ldots,P_{m,j}\subset\HH^3$ with orthogonal projections $\Pi_{i,j}:\HH^3\to P_{i,j}$ satisfying:
    \begin{enumerate}
        \item[(1)] $E_{i,j}:=K_j\cap P_{i,j}$, $i\in\{1,\cdots,m\}$, are disjoint convex disks, and $F_j(S)\cap E_{i,j}=F_j(\gamma_i)=\partial E_{i,j}$;
        \item[(2)] $\partial K_j=F_j(S)\cup E_{1,j}\cup\cdots\cup E_{m,j}$;
        \item[(3)] $\Pi_{i,j}(K_j)\subset E_{i,j}$ for each $i\in\{1,\cdots,m\}$.
    \end{enumerate}

    Now, fix a basepoint $x_*\in\Int S$ and a basepoint $o\in\HH^3$. After post-composing each $F_j$ with a hyperbolic isometry, we may normalize the sequence so that $F_j(x_*)=o$. Fix $R>\operatorname{diam}_{d_0}(S)+1$. Then for all sufficiently large $j$, we have
    $$
    d_{\HH}(o,F_j(x))\leq d_j(x_*,x)<R \qquad \forall  x\in S.
    $$
    Since $K_j$ is the convex hull of $F_j(S)$, we also have $K_j\subset\overline{B}_{d_\HH}(o,R)$. By \cref{thm:blashke-rewrite}, after passing to a subsequence, $K_j$ converges in Hausdorff distance to a nonempty compact convex subset $K\subset\HH^3$. Since every plane $P_{i,j}$ intersects the fixed ball $\overline{B}_{d_\HH}(o,R)$, after passing to a further subsequence, we may assume that $P_{i,j}$ converges to a geodesic plane $P_i$ for each $i\in\{1,\ldots,m\}$. Denote by $\Pi_i$ the orthogonal projection onto $P_i$.

    By \cref{ccv:bordered-intrinsic-realization}, we obtain that $\dim K\geq 2$ and that $F_j$ converges uniformly on $S$ with respect to $d_W$ to a continuous map $F:S\to\mathcal{B}_{\HH}(K)$ whose restriction to $\Int S$ is an intrinsic isometric embedding. If $\dim K=2$, then $\mathcal{B}_{\HH}(K)$ is the metric double $D(K)$ and $\Int S$ admits an intrinsic isometric embedding into $D(K)$, which is a contradiction. Thus, we have $\dim K=3$ and $\mathcal{B}_{\HH}(K)=\partial K$, which implies that $F$ is an actual map from $S$ into $\HH^3_K$. Since $F$ is the uniform limit of $F_j$, we see that $F(\gamma_i)\subset P_i$ for each $i$.

    To analyze the boundary behavior, fix $i\in\{1,\ldots,m\}$ and let $S^{(i)}$ denote the partial metric double of $S$ obtained by gluing two copies $S^+$ and $S^-$ of $S$ along $\gamma_i$ only, endowed with the glued metrics $d_j^{(i)}$ and $d_0^{(i)}$. Locally, $(S^{(i)},d_0^{(i)})$ coincides either with $(S,d_0)$ or with its full metric double. It then follows from \cref{cor:perelman-double} and \cref{thm:Globalization-theorem} that $(S^{(i)},d_0^{(i)})$ is a compact $\CBB(-1)$ surface with $2m-2>0$ boundary components. The same conclusion holds for $(S^{(i)},d_j^{(i)})$. It is easy to see that each copy $S^\pm$ is isometrically embedded into $S^{(i)}$ and for every $x^+\in S^+$ and $y^-\in S^-$,
    $$
    d_0^{(i)}(x^+,y^-)=\inf_{z\in \gamma_i}\bigl(d_0(x,z)+d_0(z,y)\bigr),\qquad d_j^{(i)}(x^+,y^-)=\inf_{z\in \gamma_i}\bigl(d_j(x,z)+d_j(z,y)\bigr).
    $$
    Consequently, we have
    $$
    \sup_{u,v\in S^{(i)}}\bigl|d_j^{(i)}(u,v)-d_0^{(i)}(u,v)\bigr|\leq 2\varepsilon_j.
    $$
    Let $\sigma_{i,j}$ and $\sigma_i$ be the reflections across the planes $P_{i,j}$ and $P_i$, respectively. The projection condition $\Pi_{i,j}(K_j) \subset E_{i,j}$ in \cref{def:eq-capped-realization} guarantees that for each $j$, the reflected union
    $$
    K_j^{(i)}:=K_j\cup\sigma_{i,j}(K_j)
    $$
    is convex. Let $F_j^{(i)}:S^{(i)}\to \partial K_j^{(i)}$ be the map defined by
    $$
    F_j^{(i)}(x^+)=F_j(x),\qquad F_j^{(i)}(x^-)=\sigma_{i,j}(F_j(x)).
    $$
    Then $F_j^{(i)}$ is an intrinsic isometric embedding. Moreover, as $j\to\infty$, $K_j^{(i)}$ converges to the set $K^{(i)}:=K\cup\sigma_i(K)$, which is then a compact convex set.
    Since $K^{(i)}$ contains the three-dimensional set $K$, we have $\dim K^{(i)}=3$. Applying \cref{ccv:bordered-intrinsic-realization} to maps $F_j^{(i)}$, we obtain a limit map $G_i:S^{(i)}\to\partial K^{(i)}$ such that $G_i|_{S^+}=F$ and the restriction of $G_i$ to $\Int S^{(i)}$ is an intrinsic isometric embedding. 
    
    We first show that $F$ is an injection on $S$. Since $F|_{\Int S}$ is a topological embedding, $F(\Int S)\cap F(\partial S)=\varnothing$. Otherwise, suppose $F(p)=F(q)$ for some $p\in \Int S$ and $q\in\partial S$ and let $q_n$ be a sequence in $\Int S$ converging to $q$. Then we have $F(q_n)\to F(q)=F(p)$, which implies that $q_n\to p$ since $F|_{\Int S}$ is a topological embedding, a contradiction. Now consider a boundary component $\gamma_i$. Since $\gamma_i\subset\Int S^{(i)}$, $G_i|_{\Int S^{(i)}}$ is a topological embedding, and $G_i|_{S^+}=F$, by a similar argument as above we see that $F$ is injective on $\gamma_i$ and $F(\gamma_i)\cap F(S\setminus \gamma_i)=\varnothing$. This implies that $F$ is an injection. It follows that $F$ is also an intrinsic isometric embedding on the whole space $S$.

    Since for each $j$, $E_{i,j}=K_j\cap P_{i,j}=\Pi_{i,j}(K_j)\subset K_j$, after passing to the Hausdorff limit, we obtain
    $$
    E_i:=K\cap P_i=\Pi_i(K)\subset K.
    $$
    It follows that each $E_i=\Pi_i(K)$ is a nondegenerate convex disk since $\dim K=3$.    Furthermore, by construction of the reflected bodies $K^{(i)}$, we have
    $$
    F(S)\subset\partial K^{(i)} \qquad\text{and}\qquad \partial K^{(i)}\cap P_i=\partial E_i.
    $$
    Hence $F(S)\cap P_i\subset\partial E_i$. It follows that $F(\gamma_i)$ is a Jordan curve embedded in the topological circle $\partial E_i$, which forces $F(\gamma_i)=\partial E_i$ and $F(S)\cap E_i=F(S)\cap P_i=F(\gamma_i)$. Note that the planes $P_i$ must be pairwise distinct since each set $F(S)\cap P_i$ is a single Jordan curve. Therefore, the sets $E_i=K\cap P_i$ can intersect only along their relative boundaries. Since $F$ is injective on $S$ and the components $\gamma_1,\ldots,\gamma_m$ are pairwise disjoint, we conclude that the convex disks $E_1,\ldots,E_m$ are pairwise disjoint. 

    For each $j$, we have $\partial K_j=F_j(S)\cup\bigcup_{i=1}^m E_{i,j}$. Passing to the Hausdorff limit yields
    $$
    \partial K=F(S)\cup\bigcup_{i=1}^m E_i.
    $$
    Thus, the map $F:S\to\partial K$ satisfies all the conditions in \cref{def:eq-capped-realization}, which completes the proof.
\end{proof}

\section{Circle-Type Convex Realizations of Surfaces of Finite Funnel-Type}\label{sec:finite-type-Realizations}


The main theorem of this section is the following:
    \begin{theorem}\label{thm:finite-funnel-end}
		Suppose that $(S,d_0)$ is a planar $\CBB(k)$ surface of finite funnel-type, where $k>-1$. Then $(S,d_0)$ is isometric to a complete convex surface $Y$ that bounds a convex set $X$ in $\mathbb{H}^3_P$ such that $\overline{X} \cap \partial \mathbb{H}^3_P$ is the disjoint union of finitely many disks.
	\end{theorem}

\subsection{Exhaustion by compact $\CBB(-1)$ surfaces with boundary}

In this section, we prove that there exists an exhaustion of $(S,d_0)$ by compact connected $\CBB(-1)$ subsurfaces.
\begin{theorem}\label{thm:CBB-subsurfaces}
Let $(S,d_0)$ be a planar $\CBB(-1)$ surface of finite funnel-type. Then there exists an exhaustion
by compact connected subsurfaces
$$
S_1\subset\operatorname{Int}S_2
\subset\operatorname{Int}S_3\subset\cdots\subset S,
\qquad
\bigcup_{j=1}^{\infty}S_j=S,
$$
such that for each $j$, $(S_j,d_j)$ is a $\CBB(-1)$ surface with boundary, where $d_j$ denotes the intrinsic path metric on $S_j$ induced by $d_0$.
\end{theorem}

We begin with the topological setup. Let $\{e_1,\ldots,e_m\}$ be the set of ends of $S$. Then, the end compactification $\mathcal{E}(S)$ is homeomorphic
to $\SS^2$. We may therefore choose a homeomorphism $\widehat f:\mathcal{E}(S)\to\SS^2$
such that $\widehat f(e_i)=p_i$
for some distinct points $p_1,\ldots,p_m\in\SS^2$. Its restriction
gives a homeomorphism
$$
f:S\to
U:=\SS^2\setminus\{p_1,\ldots,p_m\}.
$$
For each $i$, choose a closed disk neighborhood $B_i$ of $p_i$ in
$\SS^2$ such that the disks $B_1,\ldots,B_m$ are pairwise disjoint,
and set $A_i:=f^{-1}(B_i\setminus\{p_i\})$. Then, $A_i$ is an annular neighborhood of the end $e_i$.

A rectifiable closed curve $\gamma\subset A_i$ is said to
\emph{surround the end} $e_i$ if
$$
[\gamma]\neq0
\qquad\text{in}\quad H_1(A_i;\ZZ)\cong\ZZ.
$$
If $\gamma$ is a Jordan curve, then it is clear that $\gamma$ surrounds $e_i$ if and only if $f(\gamma)$ separates $p_i$ from
$\{p_1,\ldots,p_m\}\setminus\{p_i\}$ in $\SS^2$. For each $i$, choose a sequence of shrinking closed spherical disks
$B_{i,j}\subset B_i$ centered at $p_i$ such that $\bigcap_{j=1}^{\infty}B_{i,j}=\{p_i\}$.
Set
$c_{i,j}:=f^{-1}(\partial B_{i,j})$
and $A_{i,j}:=
f^{-1}\bigl(B_{i,j}\setminus\{p_i\}\bigr)$.
Then $A_{i,j}$ is a closed annular tail of the end $e_i$ with boundary
$c_{i,j}$.

\subsubsection{Shortest curves surrounding the ends}
In this subsection, we prove the following result about the existence of shortest curves surrounding the ends:
\begin{proposition}\label{prop:minimizing-end-curves}
For every $i\in\{1,\ldots,m\}$ and every $j\geq1$, there exists a
rectifiable Jordan curve $\gamma_{i,j}\subset A_{i,j}$ that surrounds
$e_i$ and satisfies
$$
l_{d_0}(\gamma_{i,j})=\inf\left\{
l_{d_0}(\gamma):
\gamma\subset A_{i,j}\text{ is a rectifiable closed curve surrounding }e_i\right\}.
$$
\end{proposition}
Since $(S,d_0)$ is a planar $\CBB(-1)$ surface of finite funnel-type, there exist a closed convex
set $X\subset\HH^3_P$ with nonempty interior and an isometry $F:(S,d_0)\to(\partial X,d_{\partial X})$.
Through this isometry, we identify $S$ with $\partial X$. Set
$$
P_X:=\partial X\cup\partial_\infty X
\subset\overline{\HH^3_P}.
$$
Then $P_X$ is a topological $2$-sphere. For each $i\in\{1,\cdots,m\}$, let $\Lambda_i$ denote the component of $\partial_\infty X$
corresponding to the end $e_i$. Then, by the funnel-type assumption, $\Lambda_i$ contains at least two points.
We first show that a sequence of curves of uniformly bounded length surrounding an end cannot escape to infinity.
\begin{lemma}
\label{lem:end-curve-coercivity}
For every $M<\infty$, there exists a compact subset
$Q_{i,j,M}\subset A_{i,j}$ such that for every rectifiable closed curve
$c\subset A_{i,j}$ satisfying
$$
[c]\neq 0\quad\text{in }H_1(A_{i,j};\mathbb Z)
\qquad\text{and}\qquad
l_{d_0}(c)\leq M,
$$
we have $Q_{i,j,M}\cap c\neq\varnothing$.
\end{lemma}

\begin{proof}
Suppose, to the contrary, that no such compact set exists. Let $Q_1\subset Q_2\subset\cdots\subset A_{i,j}$
be a compact exhaustion of $A_{i,j}$. Then, for each $n$, there is a
rectifiable closed curve $c_n\subset A_{i,j}\setminus Q_n$ such that $$
[c_n]\neq 0\quad\text{in }H_1(A_{i,j};\mathbb Z)
\qquad\text{and}\qquad
l_{d_0}(c_n)\leq M.
$$
Since $A_{i,j}$ is closed in $S$, the curves $c_n$ eventually escape every compact subset of $S$. For every $n$, choose a point $x_n\in c_n$. After
passing to a subsequence, we may assume that $x_n$ converges to a point $q\in\partial\HH^3_P$ in $\overline{\HH^3_P}$. Since the lengths $ l_{d_0}(c_n)$ are uniformly bounded, for every neighborhood $U$ of $q$ in $\overline{\HH^3_P}$, we have $c_n \subset U$ for all sufficiently large $n$.

By definition, $\Lambda_i$ contains at least two points. Choose a point $p\in \Lambda_i\setminus\{q\}$ and a point $b_{i,j}\in S\setminus A_{i,j}$. Let $W$ be a disk neighborhood of $q$ in the topological sphere $P_X$ such that $\partial A_{i,j}\cap W=\varnothing$ and $p\notin W$. For all sufficiently
large $n$, we have $c_n\subset W$, and hence
$$
[c_n]=0
\quad\text{in}\quad
H_1\bigl(P_X\setminus\{b_{i,j},p\};\mathbb Z\bigr).
$$
Since $P_X \setminus \{b_{i,j}, p\}$ is an annulus and the curve $\partial A_{i,j} \subset P_X$ separates $b_{i,j}$ from $p$, the inclusion map $A_{i,j} \hookrightarrow P_X \setminus \{b_{i,j}, p\}$ induces an isomorphism
$$
H_1(A_{i,j}; \mathbb{Z}) \to H_1\bigl(P_X \setminus \{b_{i,j}, p\}; \mathbb{Z}\bigr).
$$
This contradicts $[c_n] \neq 0$ in $H_1(A_{i,j}; \mathbb{Z})$, which completes the proof.
\end{proof}

We can now prove the existence of a shortest curve surrounding each end $e_i$ in $A_{i,j}$.

\begin{proof}[Proof of Proposition~\ref{prop:minimizing-end-curves}]
Fix $i\in\{1,\ldots,m\}$ and $j\geq1$. Define
$$
\mathcal A_{i,j}
:=
\left\{
c:[0,1]\to A_{i,j}:
c\text{ is a rectifiable closed curve and }
[c]\neq0\text{ in }H_1(A_{i,j};\mathbb Z)
\right\}
$$
and set $\mu_{i,j}:=\inf_{c\in\mathcal A_{i,j}}l_{d_0}(c)$.
The class $\mathcal A_{i,j}$ is nonempty since $A_{i,j}$ is an annulus in the length space $S$.

Now, choose a minimizing sequence $\{c_n\}_{n=1}^{\infty}\subset
\mathcal A_{i,j}$ such that
$l_{d_0}(c_n)\to\mu_{i,j}$ as $n\to\infty$.
Let $M=\sup_n l_{d_0}(c_n)$. By \cref{lem:end-curve-coercivity}, there exists a compact set $Q_{i,j,M}\subset A_{i,j}$ such that $c_n\cap Q_{i,j,M}\neq\varnothing$ for every $n$. After reparametrization, we may assume that each $c_n$ is parametrized with constant speed and satisfies $c_n(0) \in Q_{i,j,M}$.
Then, we have $c_n
\subset
\overline{B}_{d_0}\left(Q_{i,j,M},M/2 \right)$ and $\Lip(c_n)=l_{d_0}(c_n)\leq M$.
Since $(S,d_0)$ is a complete locally compact length space, $\overline{B}_{d_0}\left(Q_{i,j,M},M/2 \right)$ is compact. By the Arzelà--Ascoli theorem, after passing to a subsequence, the maps
$c_n$ converge uniformly to a Lipschitz map $\gamma_{i,j}:\SS^1\to A_{i,j}$.
Here $\gamma_{i,j}\subset A_{i,j}$ since $A_{i,j}$ is closed in $S$. The lower semicontinuity of length yields
\begin{equation}\label{eq:minimal}
    l_{d_0}(\gamma_{i,j})\leq\liminf_{n\to\infty}l_{d_0}(c_n)=\mu_{i,j}.
\end{equation}

Next, we verify that $\gamma_{i,j}$ surrounds the end $e_i$. Let $\pi : \mathbb{S}^1 \times [0,\infty) \to \mathbb{S}^1$ denote the projection onto the first factor. Choose a homeomorphism
$h_{i,j}:A_{i,j}\to\SS^1\times[0,\infty)$ and define $\rho_{i,j} := \pi \circ h_{i,j}$. Let $N=\overline B_{d_0}(Q_{i,j,M},M/2)\cap A_{i,j}$. This is a compact subset of $A_{i,j}$ containing all $c_n$ and their uniform limit $\gamma_{i,j}$. Hence $\rho_{i,j}\circ c_n\to\rho_{i,j}\circ\gamma_{i,j}$ uniformly. Since two sufficiently close maps from $\SS^1$ to $\SS^1$ are homotopic, for all sufficiently large $n$, we have
$$
\deg(\rho_{i,j}\circ\gamma_{i,j})=\deg(\rho_{i,j}\circ c_n)\neq0.
$$
Thus $[\gamma_{i,j}]\neq0$ in $H_1(A_{i,j};\mathbb Z)$ and $\gamma_{i,j}$ surrounds $e_i$. It then follows from \eqref{eq:minimal} and the definition of $\mu_{i,j}$ that $l_{d_0}(\gamma_{i,j})=\mu_{i,j}>0$.

It remains to show that $\gamma_{i,j}$ is a Jordan curve. Suppose, for contradiction, that there are two distinct
parameters $s$ and $t$ such that $\gamma_{i,j}(s)=\gamma_{i,j}(t)$.
Cutting $\gamma_{i,j}$ at this point produces two rectifiable closed
curves $\alpha$ and $\beta$ such that
$$
[\gamma_{i,j}]=[\alpha]+[\beta]
\quad\text{in }H_1(A_{i,j};\mathbb Z).
$$
Since $[\gamma_{i,j}]\neq0$, at least one of $[\alpha]$ and $[\beta]$
is nonzero. However, both $\alpha$ and $\beta$ have length strictly less than
$l_{d_0}(\gamma_{i,j})$, contradicting the minimality of
$\gamma_{i,j}$. Therefore, $\gamma_{i,j}$ is a minimizing
rectifiable Jordan curve surrounding the end $e_i$. This completes the proof.
\end{proof}

\subsubsection{Construction of $\CBB(-1)$ subsurfaces}

Fix $i\in\{1,\ldots,m\}$ and $j\geq1$. Let $\gamma_{i,j}$ be the Jordan curve defined in \cref{prop:minimizing-end-curves} and let $E_{i,j}$ be the component
of $S\setminus\gamma_{i,j}$ containing the end $e_i$. Set $S_{i,j}:=S\setminus E_{i,j}$. Then $S_{i,j}$ is the closed side of $\gamma_{i,j}$ away from $e_i$ and $\partial S_{i,j}=\gamma_{i,j}$.
Since $\gamma_{i,j}$ is a Jordan curve in the closed annular
tail $A_{i,j}$ surrounding $e_i$, the component $E_{i,j}$ is contained in $A_{i,j}$.

\begin{lemma}
\label{lem:no-end-side-shortcut}
Let $p\in\gamma_{i,j}$ and $B_p$ be a topological disk neighborhood of $p$ in $S$ such that $B_p\cap\gamma_{i,j}$ is a connected arc and
separates $B_p$ into two connected components. Suppose that $\sigma\subset B_p\cap\overline{E_{i,j}}$ is a rectifiable arc with distinct endpoints $a,b\in B_p\cap\gamma_{i,j}$, and $\sigma_0$ is the subarc of
$B_p\cap\gamma_{i,j}$ joining $a$ to $b$. Then
$$
l_{d_0}(\sigma_0)\leq l_{d_0}(\sigma).
$$
\end{lemma}

\begin{proof}
By our choice of $B_p$, the arcs $\sigma_0$ and
$\sigma$ are homotopic relative to their endpoints inside
$B_p\cap\overline{E_{i,j}}\subset A_{i,j}$.
Replace the subarc $\sigma_0$ of $\gamma_{i,j}$ by $\sigma$. Then the resulting
rectifiable closed curve $c\subset A_{i,j}$ satisfies
$$
[c]=[\gamma_{i,j}]
\neq0
\quad\text{in }H_1(A_{i,j};\mathbb Z).
$$
By the minimality of $\gamma_{i,j}$,
$$
l_{d_0}(\gamma_{i,j})
\leq l_{d_0}(c)=
l_{d_0}(\gamma_{i,j})
-l_{d_0}(\sigma_0)+l_{d_0}(\sigma).
$$
It follows that
$l_{d_0}(\sigma_0)\leq l_{d_0}(\sigma)$.
\end{proof}

Set $S_j:=\bigcap_{i=1}^m S_{i,j}$. It is clear that $S_j$ is a compact subset of $S$ and $\partial S_j=\bigsqcup_{i=1}^m\gamma_{i,j}$ is the disjoint union of rectifiable curves. In particular, the intrinsic path metric $d_j$ on $S_j$ induced by $d_0$ is well-defined and $(S_j, d_j)$ is a compact length space.
\begin{proposition}\label{prop:Sj-is-CBB}
    $(S_j,d_j)$ is a $\CBB(-1)$ surface with boundary.
\end{proposition}
\begin{proof}
    We first show that for every $p\in S_j$, there exists a neighborhood $V_p$ of $p$ in $S$ such that 
    \begin{equation}\label{eq:local-geodesic-convex}
        d_{j}(y,z)=d_0(y,z)\qquad\forall y,z\in S_{j}\cap V_p.
    \end{equation}
   Since this holds trivially for $p\in\Int S_j$, we only need to check the case $p\in\partial S_j$. Let $\gamma_{i,j}$ be the component of $\partial S_j$ containing $p$. Choose a topological disk neighborhood $B_p$ of $p$ as in \cref{lem:no-end-side-shortcut}. After shrinking $B_p$ if necessary, we may assume that $B_p \cap S_j = B_p \cap S_{i,j}$. Let $V_p\Subset B_p$ be a smaller neighborhood of $p$ in $S$ such that every $d_0$-minimizing geodesic joining two points of $V_p$ is contained in $B_p$.

    Let $y,z\in S_{j}\cap V_p$, and let $\eta:[0,1]\to S$ be a $d_0$-minimizing geodesic from $y$ to $z$ in $S$. By the choice of $V_p$, we have $\eta\subset B_p$. For each connected component of $\eta^{-1}(E_{i,j})$, its closure is a closed interval whose image under $\eta$ is a $d_0$-minimizing subarc $\sigma$ with endpoints $a,b\in\gamma_{i,j}$. Let $\sigma_0$ be the subarc of $B_p\cap\gamma_{i,j}$ joining $a$ to $b$. Since $\sigma$ is a subarc of a minimizing geodesic, we have $l_{d_0}(\sigma)=d_0(a,b)\leq l_{d_0}(\sigma_0)$. On the other hand, \cref{lem:no-end-side-shortcut} gives $l_{d_0}(\sigma_0)\leq l_{d_0}(\sigma)$. It follows that $l_{d_0}(\sigma_0)=l_{d_0}(\sigma)$. The open set $\eta^{-1}(E_{i,j})\subset[0,1]$ is a countable union of pairwise disjoint open intervals. On the closure of each such interval, we replace the corresponding subarc $\sigma$ by the arc $\sigma_0\subset\gamma_{i,j}$. The resulting curve $\tilde{\eta} :[0,1]\to S_{i,j}$ from $y$ to $z$ satisfies $l_{d_0}(\widetilde{\eta})= l_{d_0}(\eta)=d_0(y,z)$.
    It follows that $d_{j}(y,z)\leq l_{d_0}(\widetilde{\eta})\leq d_0(y,z)$. The reverse inequality $d_0(y,z)\leq d_{j}(y,z)$ is immediate since every path in $S_{j}$ is also a path in $S$. Therefore, we have
    $$
    d_{j}(y,z)=d_0(y,z),\qquad\forall y,z\in S_{j}\cap V_p.
    $$
    
    Now, we show that $(S_j,d_j)$ is $\CBB(-1)$. Since $(S,d_0)$ is locally $\CBB(-1)$, by \eqref{eq:local-geodesic-convex}, $(S_j,d_j)$ is also locally $\CBB(-1)$. Moreover, by \eqref{eq:local-geodesic-convex}, the topology on $S_j$ induced by $d_j$ agrees with the subspace topology induced by $d_0$. It follows that $(S_j,d_j)$ is a compact, and hence complete, length space. \cref{thm:Globalization-theorem} then implies that $(S_j,d_j)$ is a $\CBB(-1)$ surface with boundary, which completes the proof.   
\end{proof}

\begin{proof}[Proof of \cref{thm:CBB-subsurfaces}]
    Let $S_j=\bigcap_{i=1}^m S_{i,j}$. Clearly, each $S_j$ is a compact subsurface of $S$ with boundary. By our choice of $A_{i,j}$, after passing to a subsequence if necessary, we may assume that $\{S_j\}_{j=1}^{\infty}$ is an exhaustion of $S$. Finally, \cref{prop:Sj-is-CBB} implies that each $(S_j, d_j)$ is a $\CBB(-1)$ space, which completes the proof.
\end{proof}

\subsection{Proof of \cref{thm:finite-funnel-end}}
\subsubsection{Circle-domain coordinates}\label{subsubsec:fin-circle}
We first summarize our results regarding the exhaustion of $(S,d_0)$ by $\CBB(-1)$ subsurfaces. Fix $k>-1$.
Let $(S,d_0)$ be a planar $\CBB(k)$ surface of finite funnel-type. By \cref{thm:CBB-subsurfaces}, there exists an exhaustion of $S$ by connected subsurfaces
$$
S_1\subset\operatorname{Int}S_2
\subset\operatorname{Int}S_3\subset\cdots\subset S,
$$
such that each $(S_j,d_j)$ is acompact genus zero $\CBB(-1)$ surface,
where $d_j$ is the intrinsic path metric induced by $d_0$, and
$\partial S_j=\gamma_{1,j}\sqcup\cdots\sqcup\gamma_{m,j}$
is the union of $m\geq1$ pairwise disjoint rectifiable Jordan curves. It follows from \cref{thm:hyperbolic-convex-cap} that for each $j$, there exists a convex cap realization $F_j$ which is an intrinsic isometric embedding from $(S_j,d_j)$ to the boundary $(\partial K_j,d_{\partial K_j})$ of a compact convex body $K_j\subset\mathbb H^3$ such that there exist totally geodesic planes $P_{1,j},\ldots,P_{m,j}$ in $\HH^3$ satisfying
\begin{enumerate}
    \item[(1)] $E_{i,j}:=K_j\cap P_{i,j}$, $i\in\{1,\cdots,m\}$, are disjoint convex disks,
    \item[(2)] $\partial K_j=F_j(S_j)\cup E_{1,j}\cup\cdots\cup E_{m,j}$, and
    \item[(3)] $F_j(S_j)\cap E_{i,j}=F_j(\gamma_{i,j})=\partial E_{i,j}$ for each $i\in\{1,\cdots,m\}$.
\end{enumerate}
Here $F_j(S_j)$ is equipped with its intrinsic path metric, with respect to which $F_j$ is an isometry. The index $i$ is chosen so that
$\gamma_{i,j}$, and hence the cap $E_{i,j}$, corresponds to the labelled end $e_i$.

First, by the convergence theorem in \cref{subsec:convergence}, we have the following proposition.
\begin{proposition}\label{prop:Kj-convergence}
After composing with hyperbolic isometries and passing to a subsequence, the compact convex sets $K_j$ converge in the Euclidean Hausdorff distance to a compact convex set $K\subset\overline{\BB^3}$. Put $X=K\cap\BB^3$. Then $X$ has nonempty interior, and the maps $F_j$ converge locally uniformly to an isometry from $(S,d_0)$ onto a connected component of $(\partial X,d_{\partial X})$.
\end{proposition}

\begin{proof}
Fix $o\in\Int S_1$ and normalize by hyperbolic isometries so that $F_j(o)=0$ for all $j$. By Blaschke's selection theorem, after passing to a subsequence, $K_j$ converges in the Euclidean Hausdorff metric to a compact convex set $K\subset\overline{\BB^3}$.

Let $\delta_j$ denote the hyperbolic intrinsic metric on $\partial K_j$, and set $\beta_j(x)=d_0(x,\partial S_j)$. For any $x,y\in S_j$, a curve on $\partial K_j$ connecting $F_j(x)$ and $F_j(y)$ that remains within $F_j(S_j)$ has length at least $d_j(x,y)$. On the other hand, any curve leaving $F_j(S_j)$ must first meet $F_j(\partial S_j)$ and therefore has length at least $d_j(x,\partial S_j)\geq\beta_j(x)$. Consequently,
\begin{equation}\label{eq:prop-Kj-convergence}
    \min\{d_j(x,y),\beta_j(x)\}
\leq\delta_j(F_j(x),F_j(y))\leq d_j(x,y).
\end{equation}
On every compact subset $A\subset S$, we have $\inf_{x\in A}\beta_j(x)\to\infty$ as $j\to\infty$. Moreover, $d_j=d_0$ on $A\times A$ for all sufficiently large $j$, since all $d_0$-minimizing geodesics in $S$ with endpoints in $A$ are contained in a fixed compact metric ball, which eventually lies in $S_j$. The inequality \eqref{eq:prop-Kj-convergence} then implies that the pullback distances $F_j^*\delta_j$ converge locally uniformly on $S\times S$ to $d_0$.

Now choose $R_j\to\infty$ such that $\overline B_{d_0}(o,R_j)\subset\Int S_j$ and $d_j=d_0$ on $\overline B_{d_0}(o,R_j)$. Applying \cref{ccv:realization-package} to the restrictions of $F_j$ to $\overline B_{d_0}(o,R_j)$ equipped with $d_0$ yields locally uniform convergence of $F_j$ to an isometry from $(S,d_0)$ onto a connected component of $\mathcal B_{\HH}(K)=\partial X$, with $\dim K\geq2$. If $\dim K=2$, this component would be an open subset of the double of a hyperbolic planar region and would thus contain an open hyperbolic domain, contradicting the fact that $(S,d_0)$ is a $\CBB(k)$ surface for some $k>-1$. Thus $\dim K=3$, which completes the proof.
\end{proof}

Recall from \cref{sec:prelim-finite-type} that the length metric on an oriented $\CBB(-1)$ surface endows its interior with a canonical conformal structure, with respect to which any orientation-preserving intrinsic isometry is conformal. Since $d_j$ and $d_0$ agree locally on $\operatorname{Int}S_j$, the conformal structure of $(\operatorname{Int}S_j, d_j)$ coincides with the restriction of that of $(S,d_0)$. Consequently, as $F_j$ is an intrinsic isometric embedding, we may choose orientations so that the restriction
$$
F_j^\circ
:=
F_j|_{\operatorname{Int}S_j}
:
\operatorname{Int}S_j
\to
F_j(\operatorname{Int}S_j)
$$
is a conformal homeomorphism onto its image.

Let $S$ be a planar surface equipped with a conformal structure, and $S_1$ be a compact subsurface in $S$ with boundary. Then, $\partial S_1$ is the disjoint union of finitely many Jordan curves in $S$.
By \cref{thm:GUT-Koebe}, there exists a conformal homeomorphism $F:S\to\Omega$, where $\Omega$ is a domain in $\widehat{\mathbb{C}}$.
On $\Int S_1$, the conformal structure of $S$ restricts to a conformal structure on $S_1$. It follows that $\Int S_1$ is conformally equivalent to $F(\Int S_1)$, which is a finitely connected Jordan domain in $\widehat{\mathbb{C}}$. By Koebe's uniformization theorem for finitely connected domains, $F(\Int S_1)$ is conformally equivalent to a circle domain $U_1$. Since both $F(\Int S_1)$ and $U_1$ are finitely connected Jordan domains, the finitely connected Carathéodory theorem (see, e.g., \cite[Chapter 15, Theorem 3.4]{Conway-II}) implies that the conformal equivalence extends to a homeomorphism from $F(S_1)$ to $\overline{U_1}$. Thus, there exists a conformal equivalence from $\Int S_1$ to $U_1$ which extends to a homeomorphism from $S_1$ to $\overline{U_1}$. Consequently, each component of $\widehat{\CC}\setminus U_1$ is a closed round disk.

Fix three distinct points
$$
a_1,a_2,a_3\in\operatorname{Int}S_1
$$
and three distinct points
$$
q_1,q_2,q_3\in\widehat{\mathbb C},
\qquad q_3=\infty.
$$
By a similar argument as above, we see that for
each $j$ there exists a normalized conformal homeomorphism
$$
u_j:\operatorname{Int}S_j\to U_j:=\widehat{\mathbb C}\setminus\bigcup_{i=1}^m D_{i,j}\qquad\text{with} \qquad u_j(a_k)=q_k,\ k\in\{1,2,3\},
$$
where $D_{1,j},\ldots,D_{m,j}$ are pairwise disjoint closed round disks. Since $\infty\in U_j$, all the disks $D_{i,j}$ are
bounded in the complex plane
$\mathbb C=\widehat{\mathbb C}\setminus\{\infty\}$.
As in the above argument, we see that each $u_j$ extends uniquely to a homeomorphism
$$
\overline{u}_j:S_j\to \overline{U_j},
\qquad
\overline{U_j}=\widehat{\mathbb C}
\setminus\bigcup_{i=1}^m\operatorname{Int}D_{i,j}.
$$
After labeling the complementary disks suitably, we may assume that $\overline{u}_j(\gamma_{i,j})=\partial D_{i,j}$ for each $i\in\{1,\cdots,m\}$. For simplicity, we continue to write $u_j$ for this extension.

Define
$$
\varphi_j
:=
F_j\circ u_j^{-1}
:
\overline{U_j}
\to
F_j(S_j)\subset\partial K_j.
$$
Then $\varphi_j$ is a homeomorphism onto $F_j(S_j)$ and
$$
\varphi_j(\partial D_{i,j})
=
F_j(\gamma_{i,j})
=
\partial E_{i,j}.
$$
Since each $F_j$ is an intrinsic isometric embedding, we see that $\varphi_j|_{U_j}$ is a conformal homeomorphism from $U_j$ to $F_j(\operatorname{Int}S_j)$.

\begin{theorem}\label{thm:fin-koebe-kernel}
After passing to a subsequence, there exist pairwise disjoint closed
round disks $D_1,\ldots,D_m\subset\widehat{\mathbb C}$ and a normalized conformal homeomorphism $u:S\to U:=
\widehat{\CC}\setminus\bigcup_{i=1}^m D_i$ such that:
\begin{enumerate}[label=\textup{(\roman*)}]
    \item
    For every $i\in\{1,\ldots,m\}$, $D_{i,j}\to D_i$ with respect to the spherical Hausdorff distance.

    \item
    The maps $u_j$ converge locally uniformly to $u$ on $S$ and the inverse maps $u_j^{-1}$ converge locally uniformly to $u^{-1}$
    on $U$.


\end{enumerate}
\end{theorem}
\begin{proof}
Since $(S,d_0)$ is a planar surface with $m$ ends, by \cref{thm:GUT-Koebe} and \cref{cor:finite-funnel-type}, there exists a conformal homeomorphism
$$
u:S\to U:=\widehat{\mathbb{C}}\setminus\bigcup_{i=1}^m D_i,
\qquad
u(a_k)=q_k,\quad k\in\{1,2,3\},
$$
where $D_1,\ldots,D_m$ are pairwise disjoint closed round disks in $\widehat{\mathbb{C}}$, labelled so that $D_i$ corresponds to the end $e_i$. For each $j\geq 1$, set $\Omega_j:=u(\operatorname{Int}S_j)$ and $g_j=u_j\circ u^{-1}$ on $\Omega_j$. By the choice of $\{S_j\}_{j=1}^\infty$, we have $\Omega_j\uparrow U$ and $g_j(q_k)=q_k$ for each $k\in\{1,2,3\}$.
By the normality of three-point-normalized univalent maps and a diagonal argument along an exhaustion of $U$, we may assume, after passing to a subsequence, that $g_j$ converges locally uniformly on $U$ in the spherical metric to a meromorphic map $g$ fixing the three points $q_k$. The limit $g$ is nonconstant and hence univalent by Hurwitz's theorem.

By the generalized Carath\'eodory kernel convergence theorem (see \cite[Chapter~15, Theorem~4.7]{Conway-II}), we see that $U_j\to U':=g(U)$ in the kernel sense and $g_j^{-1}\to g^{-1}$ locally uniformly on $U'$. Since $\infty\in U'$, the disks $D_{i,j}$ eventually lie in a
fixed bounded subset of $\CC$. Passing to a further subsequence, we may assume that $D_{i,j}\to D_i'$ in the spherical Hausdorff distance, where each $D_i'$ is a closed round disk or a point. Let $E=\bigcup_i D_i'$. It follows that $U'$ is the component of
$\widehat{\CC}\setminus E$ containing $\infty$. Since each component of $\widehat{\CC}\setminus U'$ must contain some $D_i'$ and $U'$ is $m$-connected, the $m$ sets $D_i'$ belong to distinct components of $\widehat{\CC}\setminus U'$ and hence are pairwise disjoint.
Consequently, $U'=\widehat{\CC}\setminus\bigcup_{i=1}^m D_i'$ is a circle domain. Koebe's uniqueness theorem for finitely connected circle domains implies that $g$ is a M\"obius transformation. Since $g$ fixes three distinct points $q_i$, we have $g=\id_U$, and the labeling of the ends identifies $D_i'=D_i$. Finally, composing $g_j$ with $u$ and $g_j^{-1}$ with $u^{-1}$ yields the desired locally uniform convergence of $u_j$ and $u_j^{-1}$, respectively. This completes the proof.
\end{proof}

\subsubsection{Transboundary extremal length and equicontinuity}\label{sec:fin-source-el}

We work in the Klein model $\HH^3_K=\BB^3$. Let $K_j,K,X$ be defined as in \cref{prop:Kj-convergence}. After applying a hyperbolic isometry, we may assume that
$0\in\Int X$. Since $K_j$ converges in Hausdorff distance to $K$ in $\RR^3$, it follows that there exists $r\in(0,1)$ and $j_0$ such that $r\BB^3\subset K_j\subset\BB^3$ for all $j\ge j_0$.

\begin{lemma}\label{lem:fin-sep-length}
Suppose $K$ is a compact convex body such that $r\BB^3\subset K\subset\BB^3$ for some $r>0$, and the Euclidean diameters of the compact subsets $A,B\subset\partial K$
satisfy
$$
 \diam_\EE A\ge a>0,
 \qquad \diam_\EE B\ge b>0.
$$
Then there exists a constant $c=c(r,a,b)>0$ such that for each rectifiable closed curve $\gamma:\SS^1\to\partial K$ separating $A$ and $B$
in $\partial K$, its Euclidean length satisfies $l_{d_\EE}(\gamma)\geq c$.
\end{lemma}

\begin{proof}
Suppose, for contradiction, that there exist $K_n,A_n,B_n$ satisfying the hypotheses with the same constants $r,a,b$, along with rectifiable closed curves $\gamma_n$ separating $A_n$ and $B_n$ in $\partial K_n$, such that $l_{d_\EE}(\gamma_n)\to 0$.

For each $n$, the radial projection $R_n:\partial K_n\to\SS^2$, defined by $R_n(x)=x/\|x\|$, and its inverse are $L$-Lipschitz with respect to the Euclidean metrics, where $L=L(r)$ depends only on $r$. Consequently, we have $l_{d_\EE}(R_n\circ\gamma_n)\to 0$ and
\begin{equation}\label{eq:lower-bound-diameter}
    \min\{\diam_\EE R_n(A_n),\diam_\EE R_n(B_n)\}
\geq \frac{\min\{a,b\}}{L}>0.
\end{equation}
Since $l_{d_\EE}(R_n\circ\gamma_n)\to 0$, we may choose open spherical disks $V_n\subset\SS^2$ containing the images of $R_n\circ\gamma_n$ such that $\diam_\EE V_n\to 0$. As $R_n\circ\gamma_n$ separates $R_n(A_n)$ and $R_n(B_n)$ on $\SS^2$, at least one of the sets $R_n(A_n)$ or $R_n(B_n)$ must be contained in $V_n$. This contradicts the uniform lower bound in \eqref{eq:lower-bound-diameter}.
\end{proof}

Let $\Phi: \mathbb{H}^3_P \to \mathbb{H}^3_K$, defined by $\Phi(x) = \frac{2x}{1+|x|^2}$, be an isometry between the two models. On the Klein ball model, we let $d_{\EE}'$ be the pullback by
$\Phi^{-1}$ of the ambient Euclidean metric in the Poincar\'e ball, and use on convex surfaces
the induced path metric from $d_{\EE}'$. This choice is because the ambient Euclidean metric on $\mathbb{H}^3_P$ is conformal to the hyperbolic metric on $\mathbb{H}^3_P$, and thus is more suitable for the extremal length argument.
The following proposition collects some results in \cite[Section~4 and Appendix A]{LuoWu}.

\begin{proposition}\label{prop:luowu-area}
(1). For all $x,y\in\RR^3$, we have $\|\Phi(x)-\Phi(y)\|\leq 2\|x-y\|$.

(2) For every hyperbolic convex surface $\partial K\subset\HH_K^3$, the area of $\partial K$ in the metric $d_{\EE}'$ satisfies $\Area_{d_{\EE}'}(\partial K)\le16\pi$.

(3). The induced path metric from $d_{\EE}'$ and the intrinsic hyperbolic metric on $\partial K$ determine the same conformal structure on $\partial K$.
\end{proposition}

\begin{theorem}\label{thm:fin-equicontinuous}
The maps
$$
 \varphi_n:(\overline U_n,d_{\SS})
 \to (\ol{\BB^3},d_\EE)
$$
are equicontinuous. That is, for every $\varepsilon>0$, there exists $\delta>0$ such that for all $n\in\mathbb N$ and all $x,y\in \overline{U_n}$, we have $\|\varphi_n(x)-\varphi_n(y)\|<\varepsilon$ whenever $d_{\SS}(x,y)<\delta$.
\end{theorem}

\begin{proof}
    Suppose, for contradiction, that there exists $\varepsilon_0>0$ and a sequence $x_k,y_k\in \overline{U_{n_k}}$ such that $d_{\mathbb S^2}(x_k,y_k)\to 0$ and $\|\varphi_{n_k}(x_k)-\varphi_{n_k}(y_k)\|\geq\varepsilon_0$ for every $k\in\NN$. Since each $\varphi_n$ is uniformly continuous, we only need to consider the case where $n_k\to\infty$. Without loss of generality, we may assume that $n_k=k$, $x_k,y_k\in U_{n_k}$, and both $\{x_k\}$ and $\{y_k\}$ converge to a point $z\in \widehat{\CC}$. If $z\in U$, then the kernel convergence of the inverse maps $u_j^{-1}$ and compact convergence of $F_j\to F$ would imply that $\|\varphi_{k}(x_k)-\varphi_{k}(y_k)\|\to 0$, which is a contradiction. Hence, $z$ lies in some component $D_i$ of $\widehat{\CC}\setminus U$. By \cref{thm:fin-koebe-kernel}, $D_{i,k}$ converges to $D_i$ in the spherical Hausdorff distance.
    
    Fix a nondegenerate closed round disk $B_*\subset U$. Then $B_*$ is contained in $U_n$ for all sufficiently large $n$. Let $\Gamma_k$ be the family of simple closed curves in $\mathcal{E}(U_{k})$ separating $\{x_k,y_k\}$ from $B_*$. By \cref{prop:sourse_vanishing}, we see that $\EL_{U_{k};D_{i,k}}(\Gamma_k)\to 0$ as $k\to\infty$.

    Let $r$ be a positive number such that $r\BB^3\subset K_n\subset\BB^3$ for all sufficiently large $n$. Since $F$ is injective, we see that there exists a uniform lower bound $a>0$ for $\diam_\EE \varphi_n(B_*)$.
    We now proceed to show that there exists a constant $c_*>0$ such that
    \begin{equation}\label{eq:fin-positive-lower-bound}
        \EL_{U_{k};D_{i,k}}(\Gamma_k)\ge c_*>0,
    \end{equation}
    which would yield the desired contradiction. 
    Let $d'_{k}$ be the induced path metric on $F_k(\Int S_{k})$ from $d'_{\EE}$ and let $d_k''$ be the pull-back metric of $d_k'$ on $U_{k}$ under $\varphi_{k}$. Let $\mm_k=(d_k'',\mu_k)$, where $\mu_k$ assigns the value $0$ to the end $e_{D_{i,k}}$ of $U_{k}$ corresponding to $D_{i,k}$ and assigns the value $1$ to all the other ends of $U_{k}$. Then, $\mm_k$ is an extended metric which is admissible for $\EL_{U_k;D_{i,k}}$. Since $A(\mathfrak{m}_k)\leq 16\pi + m$ by \cref{prop:luowu-area}, we have
    $$
    \EL_{U_k;D_{i,k}}(\Gamma_k)\geq \frac{l_{\mm_k}^2(\Gamma_k)}{A(\mm_k)}\geq \frac{l_{\mm_k}^2(\Gamma_k)}{16\pi + m}.
    $$
    Then, it suffices to show that there exists $c_*>0$ such that for every $k$ and every simple closed curve $\gamma$ in $\mathcal{E}(U_{k})$ separating $\{x_k,y_k\}$ from $B_*$, we have $l_{\mm_k}(\gamma)\geq c_*$. If $\gamma$ passes through no end of $U_k$, then $\gamma$ is a simple closed curve lying entirely in $U_k$. By \cref{prop:luowu-area} and \cref{lem:fin-sep-length}, we have 
    $$
    l_{\mm_k}(\gamma)\geq l_{d_{k}''}(\gamma)=l_{d_{\EE}'}(\varphi_k(\gamma))\geq\frac{1}{2} l_{d_{\EE}}(\varphi_k(\gamma))\geq \frac{c}{2}
    $$
    for some constant $c>0$. If $\gamma$ passes through an end of $U_k$ which is different from $e_{D_{i,k}}$, then we have $l_{\mm_k}(\gamma)\geq 1$ by the choice of $\mm_k$. It remains to consider the case where $\gamma$ passes through and only through the end $e_{D_{i,k}}$. Without loss of generality, we may assume that $l_{\mm_k}(\gamma)<+\infty$. Then, one can lift $\gamma$ back to a simple arc $\widetilde{\gamma}$ in $\partial K_k$ such that the two endpoints $p_1,p_2$ of $\widetilde{\gamma}$ lie in the boundary $\partial E_{i,k}$. Let $[p_1p_2]$ be the line segment connecting $p_1$ and $p_2$. Since $E_{i,k}\subset K_k$ is a $2$-dimensional convex set, $\eta:=\widetilde{\gamma}\cup [p_1p_2]$ is a simple closed curve in $\partial K_k$ separating $\varphi_k(B_*)$ and $\{\varphi_k(x_k),\varphi_k(y_k)\}$. By \cref{lem:fin-sep-length}, we have
    $$
    l_{\mm_k}(\gamma)\geq l_{d_{k}''}(\gamma)=l_{d_{\EE}'}(\widetilde \gamma)\geq\frac{1}{2} l_{d_{\EE}}(\widetilde \gamma)\geq \frac{1}{4} l_{d_{\EE}}(\eta)\geq\frac{c}{4}.
    $$
    Then, \eqref{eq:fin-positive-lower-bound} holds for $c_* = \frac{\min\{c/4,1\}^{2}}{16\pi+m}>0$, which completes the proof. 
\end{proof}
The following lemma is a direct consequence of \cref{thm:fin-equicontinuous} and a generalized form of the Arzela-Ascoli theorem \cite[Theorem 5.8]{LuoWu}:

\begin{lemma}\label{lem:fin-varphi}
After passing to a subsequence, the maps $\varphi_n:\overline{U_n}\to\overline{\BB^3}$ converge to a continuous map $\varphi:\overline U\to\overline{\BB^3}$ such that
\begin{equation}
 \varphi_j(z_j)\to\varphi(z)\qquad \text{whenever}\qquad z_j\in\overline{U_j},\quad z_j\to z\in\overline U.
 \label{eq:fin-short-moving-points}
\end{equation}
In particular, we have $\varphi|_{U}=F\circ u^{-1}$.
\end{lemma}

Recall that for each $j$, $F_j(S_j)$ lies in a convex surface $\partial K_j$ in $\HH^3_K$ such that each boundary component $F_j(\gamma_{i,j})$ of $F_j(S_j)$ lies in a geodesic plane $P_{i,j}$. Moreover, we have
$$
\partial K_j=F_j(S_j)\cup\bigsqcup_{i=1}^m E_{i,j},
$$
where $E_{i,j}=K_j\cap P_{i,j}$. In the Klein ball model, each geodesic plane $P_{i,j}$ is the intersection of $\BB^3$ with some Euclidean plane $\Pi_{i,j}$ in $\mathbb{R}^3$. For each $i\in\{1,\cdots,m\}$, let $\Pi_i$ be the Euclidean plane which is the Hausdorff limit of the planes $\Pi_{i,j}$ in $\mathbb{R}^3$ and set $\Lambda_i:=\varphi(\partial D_i)$.

\begin{lemma}\label{lem:fin-lambda-circle}
 $\Lambda_i$ is either a point or a round circle in the sphere $\SS^2$.
\end{lemma}
\begin{proof}
We first show that $\Lambda_i\subset\mathbb{S}^2$. Since $F(S)$ is a connected component of the boundary of some closed convex set $X\subset\HH^3$ by \cref{prop:Kj-convergence}, $F(S)$ is closed in $\mathbb H^3$. Thus, the homeomorphism $F:S\to F(S)$ is proper as a map from $S$ into $\HH^3$. Suppose, for contradiction, that there exists $z\in \partial D_i$ such that $\varphi(z)\in\BB^3$. Let $z_k$ be a sequence in $U$ converging to $z$. Then
$u^{-1}(z_k)$ tends to the end $e_i$ of $S$ and is not contained in any compact subset of $S$. By the continuity of $\varphi$, $F(u^{-1}(z_k))=\varphi(z_k)$ converges to $\varphi(z)\in\BB^3$. Thus, $\{\varphi(z)\}\cup\{F(u^{-1}(z_k)):k\geq1\}$ is a compact set in $\BB^3$ whose preimage $\{u^{-1}(z_k):k\geq1\}$ under $F$ is noncompact, which is a contradiction.

Since $\partial D_{i,j}\to \partial D_i$ in spherical Hausdorff distance and
$\varphi_j(\partial D_{i,j})=\partial E_{i,j}$, the equicontinuity of $\{\varphi_j\}_j$ implies that $\partial E_{i,j}\to \Lambda_i$ in Euclidean Hausdorff distance as $j\to\infty$. Set $E_i:=C_{\EE}(\Lambda_i)$, where $C_{\EE}(\Lambda_i)$ is the convex hull of $\Lambda_i$ in the Euclidean space. Since $E_{i,j}=C_{\EE}(\partial E_{i,j})$ for every $j$ and the Euclidean Hausdorff distance between the two sets $A,B\subset\RR^3$ is always not less than the Euclidean Hausdorff distance between $C_\EE(A)$ and $C_\EE(B)$, we also have $E_{i,j}\to E_i$ in Euclidean Hausdorff distance as $j\to\infty$.

We now show that $\Lambda_i$ is either a point or a round circle. If $\Lambda_i$ is a singleton, then there is nothing to prove. Otherwise, there exist two distinct points $\xi,\eta\in\Lambda_i$. Since $E_{i,j}\to E_i$ in Euclidean Hausdorff distance, we can choose $\xi_j,\eta_j\in\partial E_{i,j}$ converging to $\xi$ and $\eta$, respectively. Then, the points $c_j:=\tfrac12(\xi_j+\eta_j)\in E_{i,j}$ converge to the point $\frac{1}{2}(\xi+\eta)\in\BB^3$. We claim that for  an arbitrary sequence $x_j$ in $\BB^3$,
\begin{equation}\label{eq:fin-1}
     (\forall j,\ x_j\in P_{i,j})\land (x_j\to x\in\Pi_i\cap\BB^3)\implies x\in E_{i}.
\end{equation}
To see this, note that the line segment $[c_j,x_j]\cap\partial E_{i,j}=\varnothing$ for all sufficiently large $j$ since $\partial E_{i,j}$ converges to the ideal set $\Lambda_i$ in Euclidean Hausdorff distance as $j\to\infty$. It then follows from the fact $c_j\in E_{i,j}$ that $[c_j,x_j]\subset E_{i,j}$ for all sufficiently large $j$. Consequently, we have $x\in E_i$ since $E_{i,j}\to E_i$ in Euclidean Hausdorff distance as $j\to\infty$. Since every $x\in\Pi_i\cap\BB^3$ can be approximated by some sequence $x_j\in P_{i,j}$, we have $\Pi_i\cap\BB^3\subset E_i$, which implies that $E_i=\Pi_i\cap\overline{\BB^3}$.
It follows that 
$$
\Lambda_i=C_{\EE}(\Lambda_i)\cap \SS^2=E_i\cap\SS^2=\Pi_i\cap\SS^2,
$$
which implies that $\Lambda_i$ is a round circle. This completes the proof.
\end{proof}

We are now ready to prove prove the main theorem of this section.
\newtheorem*{theorem1}{\cref{thm:finite-funnel-end}}
\begin{theorem1}
	Suppose that $(S,d_0)$ is a planar $\CBB(k)$ surface of finite funnel-type, where $k>-1$. Then $(S,d_0)$ is isometric to a complete convex surface $Y$ that bounds a convex set $X$ in $\mathbb{H}^3_P$ such that $\overline{X} \cap \partial \mathbb{H}^3_P$ is the disjoint union of finitely many disks.
\end{theorem1}

\begin{proof}
Let $F$ and $X$ be defined as in \cref{prop:Kj-convergence}.
Let $R(x)=x/|x|$ be the radial projection to the unit sphere and put
$$
 A_{i,j}=R(E_{i,j}),\qquad A_i=R(E_i),\qquad V=R(F(S)).
$$
As there exists $r\in(0,1)$ such that $r\BB^3\subset K_j$ for all $j\in\NN$, all these sets are well defined and the sets $A_{i,j}$ and $A_i$ are convex in $\SS^2$.
Since we have shown in \cref{lem:fin-lambda-circle} that $E_{i,j}\to E_i$ in Euclidean Hausdorff distance, it follows that $A_{i,j}\to A_i$ in Euclidean Hausdorff distance and $\partial A_i=\Lambda_i$. Then $A_i=\Lambda_i$ if $\Lambda_i$ is a point and $A_i$ is the closed disk bounded by $\Lambda_i$ if $\Lambda_i$ is a circle.

Clearly, the radial projection restricts to a homeomorphism from $F(S)$ to $V$. We claim that $V=\SS^2\setminus\bigcup_iA_i$. First, for each $j$ we have
$$
\SS^2=R\bigl(\varphi_j(\overline{U_j})\bigr)\cup\bigcup_iA_{i,j}.
$$
Since $\varphi_j(\overline{U_j})\to\varphi(\overline U)=F(S)\cup\bigcup_i\Lambda_i$ and $A_{i,j}\to A_i$ in Euclidean Hausdorff distance, we also have $\SS^2=V\cup\bigcup_iA_i$. It remains to show that $V\cap\bigcup_iA_i=\varnothing$. Note that  $V\cap\Lambda_i=\varnothing$ for every $i\in\{1,\cdots,m\}$ since $\Lambda_i\subset\partial_{\infty}X$. It follows that either $V\cap \Int A_i\neq\varnothing$ for some $i$ or $V\cap\bigcup_iA_i=\varnothing$. Suppose, for contradiction, that $V\cap \Int A_i\neq\varnothing$ for some $i$. Choose $p\in S$ such that $v:=R(F(p))\in V\cap \Int A_i$. Let $v_j=R(F_j(p))$; then $v_j\notin A_{i,j}$ since $F_j(p)\notin E_{i,j}$ and $v_j\to v\in\Int A_i$. On the other hand, since $A_{i,j}\to A_i$ in Hausdorff distance and both $A_{i,j}$ and $A_i$ are convex sets in $\SS^2$, we must have $v_j\in A_{i,j}$ for sufficiently large $j$, which is a contradiction.

Since $V$ is homeomorphic to $S$ and $S$ has $m$ ends, the complement $\SS^2\setminus V$ has exactly $m$ connected components. Consequently, the compact connected sets $A_i$ ($i\in\{1,\dots,m\}$) are pairwise disjoint, and we obtain the disjoint union decomposition $\SS^2=V\sqcup\bigsqcup_{i=1}^m A_i$. Now, for each $i$ such that $A_i$ is a disk, let $H_i\subset\HH^3_K$ be the closed hyperbolic half-space bounded by the geodesic plane $\Pi_i\cap\BB^3$ with $\partial_\infty H_i = A_i$, which satisfies $\Int H_i\cap X=\varnothing$. Define
$$
Y:=X\cup\bigcup_{i:\,A_i\text{ is a disk}}H_i.
$$
Then $Y$ is a full-dimensional closed convex set in $\HH^3_K$ satisfying $\partial_{\infty} Y=\bigsqcup_{i=1}^m A_i$. Moreover, by construction we have $R(\partial Y)=V$, which implies that $\partial Y=F(S)$ and that $F:(S,d_0)\to(\partial Y,d_{\partial Y})$ is a convex realization. Finally, since $S$ is of finite funnel-type, \cref{cor:finite-funnel-type} rules out singleton components of $\partial_\infty Y$. Therefore, every $A_i$ is a nondegenerate closed round disk. This completes the proof.
\end{proof}

\section{Circle-Type Convex Realizations of General Planar $\CBB(-1)$ Surfaces}\label{sec:general-Realizations}

\subsection{A complete conformal blow-up on finite-type domains}
\label{subsec:complete-blowup}

Throughout this section, let $(S,d_0)$ be a planar $\CBB(-1)$ surface, and let $\mu_0$ denote its intrinsic area measure. By Reshetnyak's representation theorem \cite[Sections~7.1--7.3]{Reshetnyak}, $(S,d_0)$ admits a canonical conformal structure, i.e., a Riemann surface atlas with holomorphic transition functions and local isothermal coordinates $z=x+iy$, in which the metric and measure take the form
$$
\mathrm{d}s_{d_0} = e^\sigma |\mathrm{d}z|, \qquad \mathrm{d}\mu_0 =e^{2\sigma}\,\mathrm{d}x\mathrm{d}y.
$$
Here, the conformal factor $\sigma \in L^1_{\mathrm{loc}}$ is the difference of two subharmonic functions, and its canonical logarithmic-potential representative uniquely specifies the length element. We refer to \cite[Sections~7.1--7.3]{Reshetnyak} and \cref{sec:appendix} for the detailed representation theorem. For any function $f$ locally integrable with respect to these isothermal coordinates, its conformally invariant distributional Laplacian $\boldsymbol{\Delta} f$ is the distribution defined by setting on any coordinate chart $\mathcal{O}$
$$
\langle \boldsymbol{\Delta} f, \varphi \rangle = \int_{z(\mathcal{O})} (f \circ z^{-1}) \Delta_0 (\varphi \circ z^{-1})\,\mathrm{d}x\mathrm{d}y, \qquad \forall \varphi \in C_c^\infty(\mathcal{O}),
$$
where $\Delta_0 = \partial_x^2 + \partial_y^2$ denotes the standard Euclidean Laplacian in the conformal coordinate. This pairing is invariant under conformal changes of coordinates. Note that on a smooth Riemannian surface, the distribution $\boldsymbol{\Delta} f$ of a smooth function $f$ coincides with the two-form $(\Delta_g f)\,\mathrm{d}A_g$, i.e., 
$$
\langle\boldsymbol\Delta f,\varphi\rangle
=\int_{\mathcal{O}} f\,\Delta_g\varphi\,\mathrm dA_g,\qquad \forall \varphi \in C_c^\infty(\mathcal{O}).
$$

Since $\sigma$ is locally a difference of subharmonic functions, its distributional Laplacian $\boldsymbol{\Delta}\sigma$ is a signed Radon measure. By Reshetnyak’s theorem, the intrinsic curvature measure $\omega_0$ of $(S,d_0)$ coincides with $-\boldsymbol{\Delta}\sigma$. If $S$ is a smooth Riemannian surface with metric $g$, then the Gauss curvature is given by $K_g=-\Delta_g\sigma$. This shows that the intrinsic curvature measure $\omega_0$ is a generalization of the curvature form $K_g\,\mathrm d A_g$ for Riemannian surfaces. The notion of curvature bounded below by $-1$ can then be translated into the distributional inequality $\omega_0 \ge -\mu_0$. In this paper, distributional inequalities are understood in the sense of testing against arbitrary nonnegative test functions. Now, given a continuous function $u$ on a connected open subset $V \subset S$, consider the conformally deformed path metric $d_{u}$ on $V$ induced by the length element $\mathrm{d}s_u = e^u \mathrm{d}s_{d_0}$. If $\boldsymbol{\Delta} u$ is a signed Radon measure, the corresponding area and curvature measures are then given by
$$
\mu_u = e^{2u}\mu_0, \qquad \omega_u = \omega_0 - \boldsymbol{\Delta} u,
$$
respectively. In particular, the perturbed metric satisfies the inequality $\omega_u \ge -\mu_u$ if $\boldsymbol{\Delta} u \leq (e^{2u} - 1)\mu_0$. If $d_u$ is complete, \cref{app:deformation-criterion} implies that $(V,d_u)$ is a $\CBB(-1)$ surface. See \cref{sec:appendix} for more details. 

In this section, our goal is to find a solution $u_V$ of $\boldsymbol{\Delta} u \leq (e^{2u} - 1)\mu_0$ on each domain $V\subset S$ such that $(V,d_{u_V})$ becomes a $\CBB(-1)$ surface and for each exhaustion $V_n$ of $S$, the metrics $d_{u_{V_n}}$ converge locally uniformly to $d_0$.

\begin{lemma}\label{lem:ball-barriers}
Given $p\in S$ and $R>0$, the function
$$
 B_{p,R}(x):=\log\frac{\sinh R}{\cosh R-\cosh d_0(p,x)},
 \qquad x\in B_{d_0}(p,R),
$$
is nonnegative, locally Lipschitz with respect to $d_0$, and satisfies on $B_{d_0}(p,R)$ the distributional inequality
\begin{equation}\label{eq:ball-barrier-inequality}
 \boldsymbol\Delta B_{p,R}
 \leq\bigl(e^{2B_{p,R}}-1\bigr)\mu_0.
\end{equation}
Moreover, the function
\begin{equation}\label{eq:reciprocal-ball-barrier}
 \beta_{p,R}(x):=
 \left(\frac{\cosh R-\cosh d_0(p,x)}{\sinh R}\right)_+,
 \qquad x\in S,
\end{equation}
is globally $1$-Lipschitz on $(S,d_0)$, supported in $\overline{B_{d_0}(p,R)}$, and takes values in $[0,\tanh(R/2)]$.
\end{lemma}

\begin{proof}
Set $h(x)=\cosh d_0(p,x)$. By radial comparison (\cref{app:radial-comparison}), we have
$$
 \boldsymbol\Delta h\leq 2h\mu_0,
 \qquad |\nabla h|^2=h^2-1\quad\text{$\mu_0$-almost everywhere}.
$$
Applying the weak chain rule to the composition $B_{p,R}=-\log(\cosh R-h)+\log\sinh R$ on $B_{d_0}(p,R)$, we compute
$$
\begin{aligned}
\boldsymbol\Delta B_{p,R}
 &\leq\frac{1}{\cosh R-h}\boldsymbol\Delta h + \frac{|\nabla h|^2}{(\cosh R-h)^2}\mu_0\\
 &\leq\left(\frac{2h}{\cosh R-h}+\frac{h^2-1}{(\cosh R-h)^2}\right)\mu_0=(e^{2B_{p,R}}-1)\mu_0.
\end{aligned}
$$

Now, let $\phi(t):=\max\{0,(\cosh R-\cosh t)/\sinh R\}$ for $t\geq 0$. We have $\beta_{p,R}=\phi\circ d_0(p,\cdot)$. The function $\phi$ is continuously differentiable on $[0,R)\cup(R,\infty)$ with derivative
$$
 \phi'(t)=
 \begin{cases}
   -\dfrac{\sinh t}{\sinh R}, & 0<t<R, \\[6pt]
   0, & t>R.
 \end{cases}
$$
Since $|\phi'(t)|\leq 1$ for all $t\neq R$ and $\phi$ is continuous at $t=R$, $\phi$ is $1$-Lipschitz on $[0,\infty)$. It follows that $\beta_{p,R}$ is $1$-Lipschitz on $(S,d_0)$ since the distance function $x\mapsto d_0(p,x)$ is $1$-Lipschitz with respect to $d_0$. Finally, note that $\phi$ is monotonically decreasing on $[0,R]$, attaining its minimum value $0$ on $[R,\infty)$ and its maximum value $\tanh (R/2)$ at $t=0$. Hence $\beta_{p,R}(S)=[0,\tanh(R/2)]\subset[0,1)$, which in turn ensures that $B_{p,R}\geq 0$ on $B_{d_0}(p,R)$.
\end{proof}

\begin{proposition}\label{prop:finite-domain-completion}
Let $V \subsetneq S$ be a non-empty connected open set. Then there exists a non-negative, locally Lipschitz function $u_V$ on $V$ with respect to $d_0$ such that the path metric $\widehat d_V$, defined by
\begin{equation}\label{eq:completed-domain-distance}
\widehat d_V(x,y) = \inf_{\gamma\subset V:x\to y} \int_\gamma e^{u_V} \mathrm{d}s_{d_0},
\end{equation}
is a complete $\CBB (-1)$ metric, and induces the original conformal structure on $V$ induced from $d_0$. Moreover, we have $\widehat d_V \ge d_0$ on $V \times V$, and
\begin{equation}\label{eq:completion-interior-estimate}
0 \leq u_V(x) \leq \log\coth \frac{\delta_V(x)}{2}, \qquad \text{where } \delta_V(x) := \operatorname{dist}_{d_0}(x, S \setminus V).
\end{equation}
\end{proposition}

\begin{proof}
The boundary distance function $\delta_V$ is $1$-Lipschitz and strictly positive on $V$. We define on $S$
\begin{equation}\label{eq:reciprocal-envelope}
f_V(x) = \sup_{p\in V} \beta_{p,\delta_V(p)}(x) 
\end{equation}
and set $u_V = -\log f_V$ on $V$. For each $p \in V$, the metric ball $B_{d_0}(p,\delta_V(p))$ is contained in $V$. By \cref{lem:ball-barriers}, $f_V$ is $1$-Lipschitz, satisfies $0 \leq f_V \leq 1$ on $S$, and vanishes identically on $S \setminus V$. Evaluating the supremum at $p = x$ yields
$$
f_V(x) \ge \beta_{x,\delta_V(x)}(x) = \tanh \frac{\delta_V(x)}{2} > 0.
$$
Consequently, $u_V$ is non-negative, locally Lipschitz on $V$, and satisfies the pointwise bounds in \eqref{eq:completion-interior-estimate}.

We next establish the distributional inequality
\begin{equation}\label{eq:envelope-curvature-inequality}
\boldsymbol{\Delta} u_V \leq (e^{2u_V} - 1)\mu_0.
\end{equation}
Fix a countable dense subset $\{p_j\}_{j=1}^\infty \subset V$. For each $m \ge 1$, set
$$
f_m = \max_{1 \leq j \leq m} \beta_{p_j,\delta_V(p_j)}\ \text{on } S \qquad \text{and}\qquad u_m = -\log f_m \ \text{on } \{f_m > 0\}.
$$
Since the expression \eqref{eq:reciprocal-ball-barrier} of $\beta_{p,R}$ depends continuously on its center and radius, we have $f_m \uparrow f_V$ pointwise on $S$. By Dini's theorem, this convergence is uniform on every compact subset of $V$. In particular, on any given compact subset $A \subset V$, $f_m$ is strictly positive for all sufficiently large $m$. For each point $x_0 \in \{f_m > 0\}$, those functions $\beta_{p_j,\delta_V(p_j)}$ that vanish at $x_0$ can be dropped without changing the maximum $f_m$ in a sufficiently small neighborhood of $x_0$. Thus, on such a neighborhood, $u_m$ is represented as the minimum of finitely many functions $-\log(\beta_{p_j,\delta_V(p_j)})$, each of which satisfies the distributional inequality $\boldsymbol\Delta u\leq\bigl(e^{2u}-1\bigr)\mu_0$ by \cref{lem:ball-barriers}. It then follows from \cref{app:minimum} that
$$
\boldsymbol{\Delta} u_m \leq (e^{2u_m} - 1)\mu_0
$$
in the sense of distributions on $\{f_m > 0\}$. Since $u_m \to u_V$ uniformly on each compact subset of $V$, testing against arbitrary non-negative test functions $\varphi \in C_c^\infty(V)$ and passing to the limit as $m \to \infty$ establishes \eqref{eq:envelope-curvature-inequality}.

Since $u_V\geq0$, it is immediate by definition that $\widehat d_V\geq d_0$. We now show that $\widehat d_V$ is a complete metric. Note that the conformal factor $e^{u_V} = 1/f_V$ is continuous and strictly positive, hence bounded above and away from zero on compact subsets of $V$. If $\overline{B}_{d_0}(q, 3r) \subset V$, then any $d_0$-geodesic connecting points in $B_{d_0}(q,r)$ remains inside $\overline{B}_{d_0}(q, 3r)$. It follows that for all $x, y \in B_{d_0}(q,r)$,
$$
d_0(x,y) \leq \widehat d_V(x,y) \leq \left( \min_{\overline{B}_{d_0}(q,3r)} f_V \right)^{-1} d_0(x,y).
$$
These estimates ensure that $\widehat d_V$ is locally bi-Lipschitz equivalent to $d_0$ on $V$, thereby inducing the original topology on $V$. Since $V$ is connected, we have $\widehat d_V(x,y) < \infty$ for all $x,y \in V$. Furthermore, since $f_V$ is $1$-Lipschitz, along each unit-speed $d_0$-rectifiable curve $\gamma(t)$, we have 
$$
\left|\frac{d}{dt}\log f_V(\gamma(t))\right|= \frac{1}{f_V(\gamma(t))} \left|\frac{d}{dt} f_V(\gamma(t))\right|
\le\frac1{f_V(\gamma(t))}\quad a.e.,
$$
which implies that
\begin{equation}\label{eq:log-envelope-lipschitz}
|\log f_V(x) - \log f_V(y)| \leq \widehat d_V(x,y), \qquad\forall x, y \in V.
\end{equation}
Let $\{x_k\}_{k=1}^\infty \subset V$ be a $\widehat d_V$-Cauchy sequence. Since $\widehat d_V \ge d_0$, $\{x_k\}$ is also a $d_0$-Cauchy sequence and thus converges in $(S,d_0)$ to some limit point $x_\infty$. By \eqref{eq:log-envelope-lipschitz}, $\{\log f_V(x_k)\}_{k=1}^\infty$ is a Cauchy sequence in $\mathbb{R}$. Hence, $f_V(x_k)$ is bounded away from zero. Since $f_V$ is continuous on $S$ and vanishes identically on $S \setminus V$, the limit point $x_\infty$ must belong to $V$. The local equivalence of $\widehat d_V$ and $d_0$ around $x_\infty$ then guarantees that $x_k \to x_\infty$ with respect to $\widehat d_V$. Therefore, $(V, \widehat d_V)$ is a complete metric space. By \cref{app:deformation-criterion}, the distributional inequality \eqref{eq:envelope-curvature-inequality} and the metric completeness imply that $(V, \widehat d_V)$ is a $\CBB (-1)$ surface, which completes the proof.
\end{proof}

\begin{proposition}\label{prop:strict-finite-funnel-approximants}
Let
$$
 V_1\Subset V_2\Subset\cdots\Subset S,
 \qquad \bigcup_nV_n=S,
 \qquad o\in V_1,
$$
be an exhaustion of $S$ by finitely connected Jordan domains, and let $(\lambda_n)_{n=1}^\infty$ be a sequence in $(1,\infty)$ converging to $1$. Then there exists a sequence of locally Lipschitz functions $v_n \ge 0$ on $V_n$ such that the conformal path metrics $d_n$, defined by the length elements $\mathrm{d}s_{d_n} = e^{v_n}\mathrm{d}s_{d_0}$ as in \cref{prop:finite-domain-completion}, satisfy 
\begin{enumerate}
\item[(1)] each $(V_n,d_n)$ is a complete $\CBB (-\lambda_n^{-2})$ surface of finite funnel type;
\item[(2)] $d_n \ge d_0$ on $V_n\times V_n$ for every $n \ge 1$;
\item[(3)] for every compact set $A \subset S$,
\begin{equation}\label{eq:strict-approximants-convergence}
\lim_{n\to\infty} \sup_A v_n = 0, \qquad \lim_{n\to\infty} \sup_{x,y\in A} |d_n(x,y) - d_0(x,y)| = 0.
\end{equation}
\end{enumerate}
\end{proposition}

\begin{proof}
For each $n \ge 1$, we define $v_n = u_{V_n} + \log\lambda_n$. It follows that $d_n = \lambda_n \widehat d_{V_n}$. By \cref{prop:finite-domain-completion}, each metric space $(V_n, \widehat d_{V_n})$ is complete, induces the original conformal structure on $V_n$, and satisfies $\widehat d_{V_n} \ge d_0$ on $V_n\times V_n$. Since $\lambda_n > 1$, scaling the metric by $\lambda_n$ preserves completeness and the conformal structure and yields $d_n \ge \widehat d_{V_n} \ge d_0$. Moreover, $(V_n, d_n)$ is a $\mathrm{CBB}(-\lambda_n^{-2})$ space, where $-\lambda_n^{-2} > -1$. By a similar argument as in \cref{subsubsec:fin-circle}, each $V_n$ is conformally equivalent to a circle domain whose complementary components are all closed disks. Consequently, \cref{cor:finite-funnel-type} implies that $(V_n, \widehat d_{V_n})$ is of finite funnel-type.

It remains to establish the uniform convergence on compact sets. For any compact subset $A \subset S$, let $a_n(A) = \operatorname{dist}_{d_0}(A, S \setminus V_n)$. Since $(S,d_0)$ is a noncompact complete $\CBB(-1)$ surface, it is proper. It follows that $\lim_{n\to\infty} a_n(A) = \infty$. Applying the interior bound \eqref{eq:completion-interior-estimate} on $A$, we obtain
$$
0 \leq \sup_A v_n \leq \log\lambda_n + \log\coth \frac{a_n(A)}{2} \to 0,\qquad n\to\infty.
$$
To establish the uniform convergence of metrics on a given compact set $A \subset S$, choose $R > 0$ such that $A \subset \overline{B}_{d_0}(o,R)$, and set $A^* = \overline{B}_{d_0}(o,3R)$. By the triangle inequality, any $d_0$-geodesic segment connecting two points of $A$ remains entirely within $A^*$. For any $x,y\in A$, evaluating the length functional of $d_n$ along a minimizing $d_0$-geodesic connecting $x$ and $y$ yields
$$
d_0(x,y) \leq d_n(x,y) \leq e^{\sup_{A^*} v_n} d_0(x,y) \qquad \text{for all } x,y \in A.
$$
Since $\sup_{A^*} v_n \to 0$ as $n \to \infty$ and $A$ is bounded in diameter, the uniform convergence asserted in \eqref{eq:strict-approximants-convergence} follows directly.
\end{proof}

\subsection{Proof of the main theorems}
\label{subsec:mian-theorem}

Throughout the remainder of this section, suppose that $(S,d_0)$ is
conformally equivalent to a circle domain.  Fix a conformal
homeomorphism $\psi:S\to\Omega\subset\widehat{\CC}$ and a point
$o\in S$.  Choose three distinct points $a_1,a_2,a_3\in\Omega$ and
normalize coordinates so that $a_1=\infty$.  By
\cite[Theorem~2.1]{NtalampekosRajala2025}, there exists a quasiround
exhaustion $\Omega_1\Subset\Omega_2\Subset\cdots\Omega_n\Subset\cdots$ of $\Omega$ by finitely connected Jordan domains. Without loss of generality, we may assume that all the points $a_1,a_2,a_3,\psi(o)\in\Omega_1$.
For each $n$, let $k_n$ be the conformal homeomorphism from $\Omega_n$ to a circle domain $U_n$ such that $k_n(a_i)=a_i$, $i\in\{1,2,3\}$. Then it is shown in \cite{NtalampekosRajala2025} that after passing to a subsequence, $k_n$ converge locally uniformly to a conformal homeomorphism $k$ from $\Omega$ onto a circle domain $U$.
The inverse maps converge locally uniformly on $U$.  

Now let $E(S)$ be the end space of $S$. For each $e\in E(S)$, let $P_e$ and $D_e$ denote the components of $\widehat{\CC}\setminus\Omega$ and $\widehat{\CC}\setminus U$ corresponding to $e$, respectively. For each $n\in\NN$, the homeomorphism $k_n$ extends to a homeomorphism $\widehat{k_n}$ between the end compactification $\mathcal{E}(\Omega_n)$ and $\mathcal{E}(U_n)$. Let $P_{n,e}$ be the component of
$\widehat{\CC}\setminus\Omega_n$ containing $P_e$; then there is a unique component $D_{n,e}$ of $\widehat{\CC}\setminus U_n$ such that $\widehat{k_n}([P_{n,e}])=[D_{n,e}]$. By \cite[Theorem~3.1]{NtalampekosRajala2025}, we have
\begin{equation}\label{eq:labelled-convergence-of ends}
 D_{n,e}\to D_e
 \quad\text{in spherical Hausdorff distance for every end }e\in E(S).
\end{equation}

Set $V_n=\psi^{-1}(\Omega_n)$. After taking a further subsequence, we assume that $V_n\Subset V_{n+1}$ and the closed ball $\overline B_{d_0}(o,n)\subset V_n$ for every $n$. Applying \cref{prop:strict-finite-funnel-approximants} to $V_n$, we obtain complete metrics
$d_n\geq d_0$ such that $(V_n,d_n)$ is a $\CBB(\kappa_n)$ surface of finite funnel-type for some $\kappa_n>-1$ and
\begin{equation}\label{eq:infinite-distance-convergence}
 \sup_{x,y\in A}|d_n(x,y)-d_0(x,y)|\to0
 \qquad \text{on every compact set }A\subset S.
\end{equation}
By \cref{thm:finite-funnel-end}, there exist full-dimensional
compact convex sets $K_n\subset\overline{\BB^3}$ and isometries
$$
 F_n:(V_n,d_n)\to
 \bigl(\mathcal B_{\HH}(K_n),\delta_{\mathcal B(K_n)}\bigr)
$$
such that $\mathscr{C}_n:=K_n\cap\SS^2$ is a finite union of disjoint round
closed disks. After post-composing $F_n$ with suitable hyperbolic isometries, we may assume that
\begin{equation}\label{eq:infinite-normalization}
 F_n(o)=0,\qquad K_n\subset\{x_3\ge0\}\cap\overline{\BB^3}.
\end{equation}
In particular, every disk component of $\mathscr{C}_n$ has a spherical radius of at most $\pi/2$. By \cref{thm:blashke-rewrite}, after passing to a subsequence, there exists a compact convex set $K\subset\overline{\BB^3}$ such that $K_n\to K$ in Hausdorff distance. It then follows from \cref{ccv:realization-package} that $\dim K\ge2$ and $F_n$ converge locally uniformly in a compact common metric space $(W,d_W)$ to an isometric embedding
\begin{equation}\label{eq:infinite-visible-limit}
 F:(S,d_0)\to
 (\mathcal B_{\HH}(K),\delta_{\mathcal B(K)}).
\end{equation}
The following lemmas are immediate:
\begin{lemma}\label{lem:infinite-inverse-control}
If $x_n\in V_n$ and $F_n(x_n)\to a\in F(S)$ in $W$, then
$x_n\to F^{-1}(a)$ in $S$.
\end{lemma}
\begin{proof}
By \cref{prop:ccv-finite-stability}, we have
$$
 d_0(F^{-1}(a),x_n)\leq d_n(F^{-1}(a),x_n)
 =\delta_{\mathcal B(K_n)}(F_n(F^{-1}(a)),F_n(x_n))
 \to\delta_{\mathcal B(K)}(a,a)=0,
$$
which completes the proof.
\end{proof}

\begin{lemma}\label{lem:infinite-stable-separation}
Suppose compact sets $A_n,B_n\subset\mathcal B(K_n)$ converge in
$W$ with respect to Hausdorff distance to disjoint compact sets $A,B\subset\mathcal B(K)$ of positive
$d_{\mathcal B(K)}$-diameter, respectively. Then, there exists $c,N>0$ such that for every $n\geq N$ and every compact set $T_n\subset \mathcal B(K_n)$ separating $A_n$ from
$B_n$, one has $\operatorname{diam}_{d_{\mathcal B(K_n)}}T_n\ge c$.
\end{lemma}
\begin{proof}
Otherwise, passing to a subsequence, there exist separators $T_n$ that converge in $W$ to a single point $p\in \mathcal B(K)$. Choose $a \in A \setminus \{p\}$ and $b \in B \setminus \{p\}$. Since $\mathcal{B}(K)$ is a topological sphere equipped with its Euclidean intrinsic metric, there exists a rectifiable path $\gamma\subset \mathcal{B}(K)\setminus \{p\}$ connecting $a$ and $b$. By \cref{ccv:curve-recovery}, there exist curves $\gamma_n$ in $\mathcal{B}(K_n)$ joining points in $A_n$ and $B_n$ that converge uniformly to $\gamma$ in $W$. It follows that $\gamma_n \cap T_n=\varnothing$ for all sufficiently large $n$, which contradicts the assumption that $T_n$ separates $A_n$ from $B_n$ in $\mathcal B(K_n)$.
\end{proof}

For every $n$, the end compactification $\mathcal E(\mathcal B_{\HH}(K_n))$ is obtained by collapsing each of the finitely many components of $\mathscr{C}_n$ to a point. As in \cref{sec:fin-source-el}, we let $d_{\EE}'$ be the pullback metric on the Klein ball model of the ambient Euclidean metric in the Poincar\'e ball model. Let $d_{n}'$ be the induced path metric on $\mathcal B_{\HH}(K_n)$ by $d_{\EE}'$. By \cref{prop:luowu-area}, we have $\operatorname{Area}_{d_{n}'}(\mathcal B_{\HH}(K_n))\le16\pi$. For a fixed component $C_*$ of $\mathscr{C}_n$, let $\mm_{n,C_*}=(d_n',\mu_n)$, where $\mu_n$ assigns value $0$ to $C_*$ and value $\diam_{\SS^2}C$ to every other component $C$ of $\mathscr{C}_n$. Then $\mm_{n,C_*}$ is an extended metric on $\mathcal B_{\HH}(K_n)$. Since $\mathscr{C}_n$ is contained in the upper half-sphere, we have
\begin{equation}\label{eq:infinite-extended-area}
 A(\mm_{n,C_*})\le16\pi+\sum_{C}\bigl(\operatorname{diam}_{\SS^2}C\bigr)^2
 \le16\pi+\pi \Area_{\SS^2}(\mathscr{C}_n)\leq 16\pi +4\pi^2.
\end{equation}
Thus, $\mm_{n,C_*}$ is an admissible metric for the restricted transboundary extremal length $\mathrm{EL}_{\mathcal B_{\HH}(K_n);C_*}$.

\begin{lemma}\label{lem:infinite-short-lift}
Let $A_n,B_n\subset\mathcal B_{\HH}(K_n)$ be disjoint compact sets, and let $\gamma\subset\mathcal E(\mathcal B_{\HH}(K_n))$ be a simple closed curve separating them, with $l_{\mm_{n,C_*}}(\gamma)<\infty$.  Then, there exists a simple closed curve $\widetilde \gamma\subset\mathcal B(K_n)$ separating $A_n$ from $B_n$ such that
\begin{equation}\label{eq:infinite-short-lift}
 l_{d_n}(\widetilde \gamma)
 \le(2+\pi)l_{\mm_{n,C_*}}(\gamma),
\end{equation}
where $d_n$ is the intrinsic path metric on $\mathcal B(K_n)$ induced by the Euclidean metric on $\overline{\BB^3}$.
\end{lemma}
\begin{proof}
We only need to consider the case where $l_{\mm_{n,C_*}}(\gamma)<+\infty$. By \cref{prop:luowu-area}, every curve $\alpha \subset \mathcal{B}_{\HH}(K_n)$ satisfies $l_{d_n}(\alpha) \leq 2 l_{d_n'}(\alpha)$. Thus, if $\gamma \subset \mathcal{B}_{\HH}(K_n)$, we can simply choose $\widetilde{\gamma} = \gamma$. Otherwise, $\gamma$ traverses finitely many ends $e_1, \dots, e_k$ of $\mathcal{B}_{\HH}(K_n)$, each visited exactly once, with the indices cyclically ordered in $\ZZ_k$. For each end $e_i$, let $C_i$ be the corresponding disk component of $\mathscr C_n$.
The preimage of $\gamma \setminus \{e_1, \dots, e_k\}$ in $\mathcal{B}(K_n)$ is the disjoint union of open arcs $\alpha_1, \dots, \alpha_k$, where for each $i\in\ZZ_k$, $\alpha_i$ has endpoints $p_i \in \partial C_i$ and $q_{i+1} \in \partial C_{i+1}$. For each $i$, let $\beta_i$ be the shortest spherical geodesic in $C_i$ joining $p_i$ and $q_i$, so that $l_{d_n}(\beta_i) = l_{\SS}(\beta_i) \leq \operatorname{diam}_{\SS^2} C_i$. Setting $\widetilde{\gamma} = \bigcup_{i=1}^k (\alpha_i \cup \beta_i)$, it is clear that $\widetilde{\gamma}$ is a simple closed curve in $\mathcal{B}(K_n)$ separating $A_n$ from $B_n$.

If $\widetilde{\gamma}$ does not traverse the component $C_*$, then
\begin{equation}\label{eq:hat-gamma}
l_{d_n}(\widetilde{\gamma}) = \sum_{i=1}^k l_{d_n}(\alpha_i) + \sum_{i=1}^k l_{d_n}(\beta_i) \leq 2 l_{d_n'}(\gamma) + \sum_{i=1}^k \operatorname{diam}_{\SS^2} C_i \leq 2l_{\mm_{n,C_*}}(\gamma).
\end{equation}
It remains to consider the case where $\widetilde{\gamma}$ traverses the component $C_*$. Without loss of generality, we may assume that $C_1 = C_*$. Let $\widetilde{\alpha} = \widetilde{\gamma} \setminus \beta_1$, which is an open arc joining $p_1$ to $q_1$ in $\mathcal{B}(K_n) \setminus \operatorname{Int}(C_1)$. As in \eqref{eq:hat-gamma}, we have $l_{d_n}(\widetilde{\alpha}) \leq 2 l_{\mm_{n,C_*}}(\gamma)$. Since $\beta_1$ is the shortest spherical geodesic connecting $p_1$ and $q_1$, and the spherical radius of $C_1$ is at most $\pi/2$, we have
$$
l_{d_n}(\beta_1) \leq \frac{\pi}{2} \Vert{}p_1 - q_1\Vert{} \leq \frac{\pi}{2} l_{d_n}(\widetilde{\alpha}) \leq \pi l_{\mm_{n,C_*}}(\gamma).
$$
Therefore,
$$
l_{d_n}(\widetilde{\gamma}) = l_{d_n}(\widetilde{\alpha}) + l_{d_n}(\beta_1) \leq (2 + \pi) l_{\mm_{n,C_*}}(\gamma),
$$
which completes the proof.
\end{proof}

Set
$$
 \Phi_n=F_n\circ\psi^{-1}\circ k_n^{-1}:U_n\to\mathcal B_{\HH}(K_n),
 \qquad G_n=\pi_{K_n}\circ\Phi_n.
$$
The maps $\Phi_n$, $G_n$ are conformal or anticonformal and each $\Phi_n$ induces a homeomorphism between the end compactification spaces $\mathcal{E}(U_n)$ and $\mathcal{E}(\mathcal{B}_{\HH}(K_n))$. By local uniform convergence of $F_n$ and $k_n$, we see that
\begin{equation}\label{eq:infinite-interior-convergence}
 \Phi_n\to
 \Phi:=F\circ\psi^{-1}\circ k^{-1}:U\to F(S)
 \quad\text{locally uniformly in }W.
\end{equation}

\begin{proposition}\label{prop:infinite-forward-equicontinuity}
The maps $G_n:(U_n,d_{\SS^2})\to
(\overline{\BB^3},d_{\EE})$ are equicontinuous. That is, for every $\varepsilon>0$, there exists $\delta>0$ such that for all $n\in\mathbb N$ and all $x,y\in U_n$, we have $\|G_n(x)-G_n(y)\|<\varepsilon$ whenever $d_{\SS}(x,y)<\delta$.
\end{proposition}
\begin{proof}
Suppose, for contradiction, that there exists $\varepsilon_0>0$ and sequences $z_m,z_m'\in U_{n_m}$ converging in $\SS^2$ to the same point $z$ such that $\|G_{n_m}(z_m)-G_{n_m}(z_m')\|\geq\varepsilon_0$ for every $m\in\NN$. 

We first consider the case where $n_m\to\infty$ as $m\to\infty$. Without loss of generality, we may assume that $n_m=m$. By \eqref{eq:infinite-interior-convergence}, the limit $z\notin U$. Let $D_e$ be the component of $\widehat{\mathbb C}\setminus U$ containing $z$. For each $m$, let
$D_{m,e}$ be the component labelled by this end as in \eqref{eq:labelled-convergence-of ends}. Then $D_{m,e}\to D_e$ in spherical Hausdorff distance. Fix a nondegenerate closed round disk $B_*\Subset U_m$ for all sufficiently large $m$. After passing to a subsequence, we may assume that the sets
$$
 A_m:=\Phi_m(B_*)\qquad\text{and}\qquad
 B_m:=\{\Phi_m(z_m),\Phi_m(z_m')\}
$$
converge in $W$ with respect to Hausdorff distance to compact sets $A,B\subset \mathcal{B}(K)$ in $W$ of
positive diameter. Clearly, $A\subset F(S)$. Moreover, we have $B\cap F(S)=\varnothing$ and hence $A\cap B=\varnothing$. 
This is because an accumulation point of $B_m$ in $F(S)$, by
\cref{lem:infinite-inverse-control} and local convergence of
$k_m$, would force $z\in U$.

Let $\Gamma_m$ be the family of simple closed curves in $\mathcal E(U_m)$
separating $\{z_m,z_m'\}$ from $B_*$. By
\eqref{eq:labelled-convergence-of ends} and
\cref{prop:sourse_vanishing}, 
$$
 \mathrm{EL}_{U_m;D_{m,e}}(\Gamma_m)\to0.
$$
Let $C_{m,e}$ be the component of $\mathscr{C}_m$ corresponding to
$D_{m,e}$ under $\Phi_m$, and let $\Gamma_m'$ be the curve
family in $\mathcal E(\mathcal B_{\HH}(K_m))$ which is the image of $\Gamma_m$ under $\Phi_m$.
By \cref{lem:infinite-stable-separation}, there exists $c_0>0$ such that every compact
separator of $A_m$ and $B_m$ has diameter at least
$c_0$ for sufficiently large $m$. It then follows from \cref{lem:infinite-short-lift} that $l_{\mm_{m,C_{m,e}}}(\Gamma_m')\geq\frac{c_0}{2+\pi}$.
Since the restricted transboundary extremal length is invariant under conformal and anti-conformal homeomorphisms, \eqref{eq:infinite-extended-area} yields
$$
 \mathrm{EL}_{U_m;D_{m,e}}(\Gamma_m)
 =
 \mathrm{EL}_{\mathcal B_{\HH}(K_m);C_{m,e}}(\Gamma_m')
 \ge\frac{c_0^2}{(2+\pi)^2(16\pi+4\pi^2)}>0,
$$
which is a contradiction. 

Now we only need to consider the case where the sequence $n_m$ is bounded. In this case we may assume that $n_m=n_0$ for some fixed $n_0\in\NN$. The same argument above and the fact that $\Phi_n$ is a homeomorphism show that $G_{n_0}$ is uniformly continuous, which is again a contradiction. This completes the proof.
\end{proof}

\begin{proof}[Proof of \cref{thm:MainThm0}]
Let the sets $U_n,U,K_n,K,D_{n,e},C_{n,e}$ and maps $\Phi_n,\Phi,G_n$ be defined as above. If $S$ is isometric to the hyperbolic plane $\HH^2$, it is clear that a hyperbolic half-space realizes
$S$ and has a round hemisphere as its ideal set. So in the following we assume that $S$ is not isometric to the hyperbolic plane $\HH^2$.

By \cref{prop:infinite-forward-equicontinuity}, the maps $G_n$ admit equicontinuous extensions $\overline G_n$ to $\overline{U_n}$, satisfying $\overline G_n(\partial D_{n,e})=\partial C_{n,e}$ by compactness and the end correspondence. Since the domains $U_n$ contain a common round disk for sufficiently large $n$, after applying a fixed M\"obius transformation to the source $\widehat{\CC}$, we may assume that all complementary disks $D_{n,e}$ have spherical radii at most $\pi/2$, while the target disks $C_{n,e}$ also have spherical radii at most $\pi/2$ by \eqref{eq:infinite-normalization}. Since $\dim K_n=3$, applying the central extension argument of \cite[Section~7.3]{LuoWu} to these boundary surjections and returning to the original source coordinates yields continuous surjections $\widehat\Phi_n:\widehat{\CC}\to\mathcal B(K_n)$ extending $\Phi_n$, such that the maps $\widehat G_n:=\pi_{K_n}\circ\widehat\Phi_n$ form an equicontinuous family and satisfy
$$
\widehat G_n(D_{n,e})=C_{n,e}
\qquad\text{for all }n\text{ and }e\in E(S).
$$
By the Arzel\`a--Ascoli theorem, after passing to a subsequence, we may assume that $\widehat G_n$ converges uniformly to a map $\widehat G:\widehat{\CC}\to\pi_K(\mathcal{B}(K))$ with $\widehat G|_U=\pi_K\circ\Phi$.
For each end $e\in E(S)$, set $C_e^*=\widehat G(D_e)$. Since
$\widehat G_n(D_{n,e})=C_{n,e}$, uniform convergence and
\eqref{eq:labelled-convergence-of ends} imply that $C_{n,e}\to C_e^*$ in spherical Hausdorff distance. Consequently, $C_e^*\subset \mathscr{C}:=K\cap\SS^2$ is either a round closed
disk or a point. Since each $\widehat{G_n}$ is a surjection onto $\pi_{K_n}(\mathcal{B}(K_n))$, and $\pi_{K_n}(\mathcal{B}(K_n))\to \pi_K(\mathcal B(K))$ in Hausdorff distance, we have
$\widehat G(\SS^2)=\pi_K(\mathcal B(K))$, which implies that $\mathscr{C}=\bigcup_{e\in E(S)}C_e^*$. Moreover, since $\widehat G(U)=\pi_K(F(S))\subset\BB^3$, we have $\pi_K(F(S))=\pi_K(\mathcal B_{\HH}(K))$.

If $\dim K=3$, injectivity of $\pi_K$ gives $F(S)=\mathcal B_{\HH}(K)$.
If $\dim K=2$, since $S$ is not isometric to the hyperbolic plane, we have $\partial_{\mathrm{rel}}K\cap\BB^3\ne\varnothing$, so the metric double $\mathcal{B}_\HH(K)$ is connected.
Invariance of domain and completeness imply $F(S)=\mathcal B_{\HH}(K)$ in this case as well. 
Since $\mathcal B(K)$ is a topological sphere, the
homeomorphism $\Phi:U\to\mathcal B_{\HH}(K)$ induces a
bijection between the complementary components of $U$
and those of $\mathcal B_{\HH}(K)$ in $\mathcal B(K)$;
see the opening paragraphs of
\cite[Section~2]{NtalampekosRajala2025}.
The restriction of $\pi_K$ identifies the latter
components with the components of $\mathscr C$.
Denote by $C_e$ the component corresponding to $D_e$.
It suffices to show that $C_e=C_e^*$. Fix an end $e\in E(S)$. For a point $z\in\partial D_e$, let $z_j\in U$ be a sequence converging to $z$ in $\SS^2$. Since $\Phi$ carries the corresponding ends
to each other, we have
$$
 \widehat G(z)
 =\lim_{j\to\infty}\pi_K(\Phi(z_j))
 \in C_e\cap C_e^*.
$$
Since $C_e^*\subset\mathscr{C}$ is connected and has nonempty intersection with component $C_e$, we have $C_e^*\subset C_e$. Suppose, for contradiction, that there exists a point $z\in C_e\setminus C_e^*$. Since $\mathscr{C}=\bigcup_{e\in E(S)}C_e^*$, $z\in C_{e'}^*\subset C_{e'}$ for some end $e'\neq e$, which contradicts the fact that $C_e\cap C_{e'}=\varnothing$. This shows that $\mathscr{C}$ is a circle-type closed set. Taking $X=K\cap\BB^3$, $X$ is a closed convex set in $\HH^3$ with $\partial_{\infty} X=\mathscr{C}$ and $F:(S,d_0)\to (\mathcal{B}_{\HH}(K),\delta_{\mathcal{B}(K)})$ is an isometry, which completes the proof.
\end{proof}

\begin{proof}[Proof of \cref{thm:MainThm1}]
He--Schramm's uniformization theorem \cite[Theorem~0.2]{HeSchramm1993} provides circle-domain coordinates for a planar surface with
countably many ends. Then the existence of circle-type convex realizations follows directly from \cref{thm:MainThm0}, while the
uniqueness follows from \cite{luo2026rigidity} and \cite{Zhao-rigidity}.
\end{proof}

\appendix

\section{Auxiliary Results on the Convergence of Convex Surfaces}\label{sec:appendix-convergence}

This appendix supplies the auxiliary results used in
\cref{subsec:convergence}. 

The following lemma concerning the convergence of convex surfaces in $\RR^3$ is well known to experts, though rarely stated explicitly in the literature. Its hyperbolic counterpart in arbitrary dimensions was rigorously established in \cite[Proposition~2.5]{AKP2023}, from which the Euclidean version below follows readily.

\begin{lemma}\label{lem:ccv-euclidean-convergence}
Let $K_n\subset\RR^3$ be a sequence of full-dimensional compact convex sets converging in the Hausdorff distance to a nonempty compact convex set $K\subset\RR^3$. Then
$$
 (\mathcal{B}(K_n),d_{\mathcal{B}(K_n)})
 \xrightarrow{\mathrm{GH}}
 (\mathcal{B}(K),d_{\mathcal{B}(K)}).
$$
\end{lemma}

\begin{lemma}\label{lem:ccv-dist-fn}
    Suppose that a sequence of compact metric spaces $(L_n, d_{L_n})$ converges to $(L, d_L)$ in the Gromov--Hausdorff distance, and that there exist surjections $f_n:L_n\to L$ and positive numbers $\varepsilon_n\to 0$ satisfying
    \begin{equation}\label{eq:ccv-dist-f-n}
        d_{L}(f_n(a),f_n(b))\leq d_{L_n}(a,b)+\varepsilon_n, \qquad \forall a,b\in L_n.
    \end{equation}
    Then
    $$
    \dist(f_n):=\sup_{a,b\in L_n}
    \bigl|d_{L}(f_n(a),f_n(b))-d_{L_n}(a,b)\bigr|\to0\quad \text{as} \quad n\to\infty.
    $$
\end{lemma}
\begin{proof}
    By a standard result on Gromov--Hausdorff convergence due to Gromov (see \cite[\S~6, p.~65]{MR623534}), $(L_n,d_{L_n})$ and $(L,d_{L})$ can be isometrically embedded into a common compact metric space $(Z,d_Z)$ such that $d_H^Z(L_n,L)\to0$ as $n\to\infty$. By \cref{thm:blashke-rewrite}, after passing to a subsequence, the closures of the graphs 
    $$
    \Gamma_n=\{(a,f_n(a)):a\in L_n\}\subset Z\times L
    $$ 
    of $f_n$ converge in the Hausdorff distance in $Z\times L$ to a compact subset $\Gamma$. We claim that $\Gamma$ is the graph of a surjective $1$-Lipschitz map $f:L\to L$. Since $d_H^Z(L_n, L)\to 0$ and $f_n$ is surjective, we have $\Gamma\subset L\times L$ and the projection of $\Gamma$ onto the second factor is exactly $L$. If both $(x,y)$ and $(x,y')$ lie in $\Gamma$, we may choose sequences $a_n, b_n \in L_n$ such that
    $$
    (a_n,f_n(a_n))\to(x,y),\qquad
    (b_n,f_n(b_n))\to(x,y').
    $$
    By \eqref{eq:ccv-dist-f-n}, we have $d_L(f_n(a_n),f_n(b_n))\leq d_Z(a_n,b_n)+\varepsilon_n\to0$. Thus, $y=y'$, and $\Gamma$ defines a surjective map $f:L\to L$. The $1$-Lipschitz property of $f$ also follows directly from \eqref{eq:ccv-dist-f-n}. Since any nonexpanding surjective self-map on a compact metric space is an isometry (see, e.g., \cite[Theorem~1.6.15(1)]{BBI}), $f:L\to L$ is an isometry.

    Now, suppose for contradiction that there exist $\delta>0$, a subsequence $\{n_k\}$, and points $a_k,b_k\in L_{n_k}$ such that
    $$
    \left|
    d_L(f_{n_k}(a_k),f_{n_k}(b_k))
    -d_{L_{n_k}}(a_k,b_k)
    \right|\ge\delta.
    $$
    By the above argument, we may assume that the closures of the graphs $f_{n_k}$ converge in the Hausdorff distance in $Z\times L$ to the graph of an isometry $f:L\to L$. Passing to a further subsequence, we assume that $(a_k,f_{n_k}(a_k))$ and $(b_k,f_{n_k}(b_k))$ converge to $(a,f(a))$ and $(b,f(b))$ in $Z\times L$, respectively. Then
    $$
    \lim_{k\to\infty}\left|d_L(f_{n_k}(a_k),f_{n_k}(b_k))-d_{L_{n_k}}(a_k,b_k)\right|= 
    \left| d_L(f(a),f(b))-d_{L}(a,b)\right|=0,
    $$
    which yields a contradiction.
\end{proof}

We are now ready to prove \cref{prop:common-W-space}.
\newtheorem*{proposition1}{\cref{prop:common-W-space}}
\begin{proposition1}
Let $K_n \subset \overline{\BB^3}$ be a sequence of full-dimensional compact convex sets converging in the Hausdorff metric to a nonempty compact convex set $K \subset \overline{\BB^3}$. Then $(\mathcal{B}(K), d_{\mathcal{B}(K)})$ and all $(\mathcal{B}(K_n), d_{\mathcal{B}(K_n)})$ can be isometrically realized as subsets of a common compact metric space $(W, d_W)$ equipped with a $1$-Lipschitz map $\Pi : W \to \overline{\BB^3}$, such that
$$
  \Pi|_{\mathcal{B}(K_n)} = \pi_{K_n}, \qquad \Pi|_{\mathcal{B}(K)} = \pi_K,
$$
and
\begin{equation}\label{app:ccv-common-space}
  d_H^W(\mathcal{B}(K_n), \mathcal{B}(K)) \to 0 \quad \text{as } n \to \infty.
\end{equation}
\end{proposition1}

\begin{proof}
We first claim that there exist continuous surjections $r_n:\mathcal{B}(K_n)\to\mathcal{B}(K)$ such that
\begin{equation}\label{eq:ccv-distortion}
 \dist(r_n)
 :=\sup_{a,b\in\mathcal{B}(K_n)}
 \bigl|d_{\mathcal{B}(K)}(r_n(a),r_n(b))-d_{\mathcal{B}(K_n)}(a,b)\bigr|
 \to0
\end{equation}
and
\begin{equation}\label{eq:ccv-coordinate-error}
 \kappa_n:=\sup_{a\in\mathcal{B}(K_n)}
       \|a-\pi_K(r_n(a))\|\to0.
\end{equation}
Set $h_n:=d_H(K_n,K)$. We divide the proof of the claim into three cases according to the dimension of $K$. 

\noindent\textbf{Case 1.} $\dim K=3$.
By translating $K$ if necessary, we may assume that the origin lies in $\Int K$, and choose $\rho>0$ sufficiently small so that $\rho\overline{\BB^3}\subset K$. Take a sequence $t_n\to0$ with $t_n\rho>h_n$. For any compact convex set $L\subset\RR^3$ and any unit vector $u\in\SS^2$, let $h_L(u)=\sup_{x\in L}\langle x,u\rangle$. Clearly, $|h_{K_n}(u)-h_K(u)|\leq h_n$ for all $u\in\SS^2$. Setting $A_n=(1-t_n)K$, we have
$$
h_{A_n}(u)=(1-t_n)h_K(u)\leq h_K(u)-t_n\rho<h_K(u)-h_n\leq h_{K_n}(u),\quad \forall u\in\SS^2,
$$
which implies that $A_n\subset \Int K_n$ for all sufficiently large $n$. For any compact convex set $L\subset\RR^3$, denote by $P_{L}:\RR^3\to L$ the nearest-point projection onto $L$. Define $r_n(a)=(1-t_n)^{-1}P_{A_n}(a)$ for $a\in\mathcal{B}(K_n)$. Since $A_n\subset\Int K_n$, $r_n$ is a continuous map from $\mathcal{B}(K_n)$ to $\mathcal{B}(K)$. To see that $r_n$ is surjective, take any $y\in\mathcal{B}(K)$ and let $x=(1-t_n)y\in\mathcal{B}(A_n)$. Let $u_x$ be an outward unit normal vector to a supporting plane of $A_n$ at $x$, and let $x'$ be the unique intersection point of $\mathcal{B}(K_n)$ with the ray $x+\RR_{\geq0} u_x$. Then $P_{A_n}(x')=x$ and thus $r_n(x')=y$, establishing surjectivity. For any pair of points $a,b\in\mathcal{B}(K_n)$, projecting shortest paths on $\mathcal{B}(K_n)$ onto $\mathcal{B}(A_n)$ yields
$$
 (1-t_n)d_{\mathcal{B}(K)}(r_n(a),r_n(b))= d_{\mathcal{B}(A_n)}(P_{A_n}(a),P_{A_n}(b))\leq d_{\mathcal{B}(K_n)}(a,b).
$$
Consequently,
$$
 d_{\mathcal{B}(K)}(r_n(a),r_n(b))
 \leq d_{\mathcal{B}(K_n)}(a,b)
       +t_n\sup_{x,y\in \mathcal{B}(K)}d_{\mathcal{B}(K)}(x,y),
$$
which implies $\dist(r_n)\to 0$ by \cref{lem:ccv-euclidean-convergence} and \cref{lem:ccv-dist-fn}. Next, we estimate $\kappa_n$. Set $R:=\max_{x\in K}\|x\|$. For every $a\in \mathcal{B}(K_n)$, there exists $b\in K$ such that $\|a-b\|\leq h_n$. Since $(1-t_n)b\in A_n$, the definition of $P_{A_n}$ gives
$$
\|a-P_{A_n}(a)\|\leq \|a-(1-t_n)b\|
\leq h_n+t_nR.
$$
On the other hand,
$$
\|P_{A_n}(a)-r_n(a)\|=\|(1-t_n)r_n(a)-r_n(a)\|
=t_n\|r_n(a)\|
\leq t_nR.
$$
Therefore, $\kappa_n\leq h_n+2t_nR\to0$, which proves \eqref{eq:ccv-coordinate-error}.

\smallskip
\noindent\textbf{Case 2.} $\dim K=2$.
After applying a rigid motion, we may assume that $K\subset H:=\RR^2\times\{0\}$. Let $P_H$ be the orthogonal projection onto $H$. Since the Hausdorff distance is nonincreasing under $P_H$, a radial scaling argument analogous to Case~1 yields boundary-preserving homeomorphisms $\varphi_n:P_H(K_n)\to K$ satisfying 
$$
\theta_n:=\sup_{x\in P_H(K_n)}\|\varphi_n(x)-x\|\to0.
$$
For each point $p$ in the relative interior of $P_H(K_n)$, the preimage $P_H^{-1}(p)\cap\mathcal{B}(K_n)$ consists of exactly two points (an upper and a lower point), which we map via $\varphi_n$ to the corresponding copies $K^+$ and $K^-$ in $\mathcal{B}(K)$, respectively. For each point $p$ in the relative boundary of $P_H(K_n)$, we map the entire set $P_H^{-1}(p)\cap \mathcal{B}(K_n)$ to the common boundary point $\varphi_n(p)\in\partial_{\mathrm{rel}}K\subset\mathcal{B}(K)$. It is straightforward to verify that this defines a continuous surjection $r_n:\mathcal{B}(K_n)\to\mathcal{B}(K)$ satisfying $\kappa_n\leq h_n+\theta_n\to0$ and
$$
 d_{\mathcal{B}(K)}(r_n(a),r_n(b))\leq d_{\mathcal{B}(K_n)}(a,b)+4\theta_n.
$$
Applying \cref{lem:ccv-euclidean-convergence} and \cref{lem:ccv-dist-fn} once more gives $\dist(r_n)\to 0$.

\smallskip
\noindent\textbf{Case 3.} $\dim K\leq1$.
If $K$ is a nondegenerate line segment, let $l$ be the line spanned by $K$ and $P_{l}$ be the orthogonal projection onto $l$. We define $r_n: \mathcal{B}(K_n)\to \mathcal{B}(K)=K$ as the restriction to $\mathcal{B}(K_n)$ of the composition of $P_l$ with the orientation-preserving affine homeomorphism from $P_l(K_n)$ onto $K$. If $K$ is a singleton $\{q\}$, we simply take $r_n$ to be the constant map with value $q$. Clearly, $r_n$ is a continuous surjection satisfying
$$
 \kappa_n=
 \sup_{a\in \mathcal{B}(K_n)}\|a-r_n(a)\|\to0
$$
and
$$
 d_{\mathcal{B}(K)}(r_n(a),r_n(b))
 =\|r_n(a)-r_n(b)\|\leq \|a-b\|+2\kappa_n\leq d_{\mathcal{B}(K_n)}(a,b)+2\kappa_n.
$$
An application of \cref{lem:ccv-euclidean-convergence} and \cref{lem:ccv-dist-fn} then completes the proof of the claim.

\medskip
Now, choose a sequence $\tau_n > 0$ converging to zero such that $2\tau_n \ge \dist(r_n)$ and $\tau_n \ge \kappa_n$. Let $W_0$ be the disjoint union of $\mathcal{B}(K)$ and all $\mathcal{B}(K_n)$. For any two points $x, y \in W_0$, a \emph{chain} from $x$ to $y$ is a finite sequence of points
$$
x = p_0, p_1, \dots, p_m = y, \quad m \ge 1,
$$
such that for each step $(p_{i-1}, p_i)$, one of the following holds:
\begin{enumerate}
  \item[(1)] $p_{i-1}$ and $p_i$ both lie in $\mathcal{B}(K)$ or both lie in $\mathcal{B}(K_n)$ for some $n$, in which case the cost of the step is $d_{\mathcal{B}(K)}(p_{i-1}, p_i)$ or $d_{\mathcal{B}(K_n)}(p_{i-1}, p_i)$, respectively;
  \item[(2)] $\{p_{i-1}, p_i\} = \{a, r_n(a)\}$ for some $n$ and some $a \in \mathcal{B}(K_n)$, in which case the cost of the step is $\tau_n$. Such a pair $\{p_{i-1}, p_i\}$ is called a \emph{bridge}.
\end{enumerate}
The cost of a chain is the sum of the costs of its steps, and we define $d_{W_0}(x,y)$ as the infimum of the costs over all chains joining $x$ and $y$. It is readily checked that $(W_0,d_{W_0})$ is a metric space.
By our choice of $\tau_n$, for all $a,b\in\mathcal{B}(K_n)$ we have
$$
d_{\mathcal{B}(K_n)}(a,b)
\leq d_{\mathcal{B}(K)}(r_n(a),r_n(b))+2\tau_n\quad\text{and}\quad d_{\mathcal{B}(K)}(r_n(a),r_n(b))
\leq d_{\mathcal{B}(K_n)}(a,b)+2\tau_n.
$$
This ensures that the distance between any two points within $\mathcal{B}(K)$ (resp.\ $\mathcal{B}(K_n)$) cannot be shortened by taking a detour through bridges, so $\mathcal{B}(K)$ (resp.\ $\mathcal{B}(K_n)$) embeds isometrically into $W_0$. Moreover, it follows immediately from the surjectivity of $r_n$ that
$$
  d_H^{W_0}(\mathcal{B}(K_n), \mathcal{B}(K)) \leq \tau_n.
$$
Notice that $W_0$ is totally bounded, since any finite $\varepsilon$-net of $\mathcal{B}(K)$ also serves as a finite $2\varepsilon$-net for $\mathcal{B}(K_n)$ in $W_0$ whenever $\tau_n < \varepsilon$. Let $(W,d_W)$ be the metric completion of $(W_0, d_{W_0})$; then $(W,d_W)$ is a compact metric space. Finally, note that the coordinate maps $\pi_K$ and $\pi_{K_n}$ are $1$-Lipschitz on each piece $\mathcal{B}(K)$ and $\mathcal{B}(K_n)$, as well as across every bridge since
$$
\|\pi_{K_n}(a)-\pi_K(r_n(a))\|=\|a-\pi_K(r_n(a))\|\leq \kappa_n \leq \tau_n.
$$ 
Consequently, they extend uniquely to a $1$-Lipschitz map $\Pi : W \to \overline{\BB^3}$. This completes the proof.
\end{proof}

Now we consider the approximation of curves and the convergence of intrinsic hyperbolic distances under the Euclidean Hausdorff convergence of compact convex sets. Here, $K_n$, $K$, and $(W, d_W)$ are as in \cref{prop:common-W-space}.

\begin{lemma}\label{ccv:curve-recovery}
Let $c : [0,1] \to \mathcal{B}(K)$ be a rectifiable curve parametrized with constant speed $\ell$ with respect to $d_{\mathcal{B}(K)}$. Suppose that sequences $p_n, q_n \in \mathcal{B}(K_n)$ converge in $(W, d_W)$ to $c(0)$ and $c(1)$, respectively. Then there exist curves $c_n : [0,1] \to \mathcal{B}(K_n)$ with $c_n(0) = p_n$ and $c_n(1) = q_n$ such that $c_n \to c$ uniformly in $(W, d_W)$ and
\begin{equation}\label{ccv:strict-recovery}
  \Lip_{d_{\mathcal{B}(K_n)}}(c_n) \leq \ell + o(1), \qquad l_{d_{\mathcal{B}(K_n)}}(c_n) \to \ell\qquad \text{as}\quad n\to\infty.
\end{equation}
If, in addition, $c([0,1]) \subset \mathcal{B}_{\HH}(K)$, then for all sufficiently large $n$, the curve $c_n$ is contained in $\mathcal{B}_{\HH}(K_n)$ and its hyperbolic length $l_{\HH}(\pi_{K_n} \circ c_n)$ in $\HH_K^3$ converges to $l_{\HH}(\pi_K \circ c)$ as $n\to\infty$.
\end{lemma}

\begin{proof}
Choose a sequence $\alpha_n > 0$ converging to zero such that
$$
\alpha_n \ge \max\bigl\{ d_H^W(\mathcal{B}(K_n), \mathcal{B}(K)),\, d_W(p_n, c(0)),\, d_W(q_n, c(1)) \bigr\},
$$
and select an integer sequence $m_n \to \infty$ satisfying $m_n \alpha_n \to 0$. For each $i \in \{0, \dots, m_n\}$, set $t_i := i / m_n$. Since $d_H^W(\mathcal{B}(K_n), \mathcal{B}(K)) \leq \alpha_n$, we may choose points $p_n^i \in \mathcal{B}(K_n)$ such that $d_W(p_n^i, c(t_i)) \leq \alpha_n$ for each $0 \leq i \leq m_n$, with the endpoints fixed as $p_n^0 := p_n$ and $p_n^{m_n} := q_n$. Define $c_n : [0,1] \to \mathcal{B}(K_n)$ to be the curve whose restriction to each subinterval $[t_{i-1}, t_i]$ is a constant-speed shortest geodesic in $\mathcal{B}(K_n)$ connecting $p_n^{i-1}$ to $p_n^i$.
Since $\mathcal{B}(K)$ and $\mathcal{B}(K_n)$ embed isometrically into $(W, d_W)$, the triangle inequality yields
$$
d_{\mathcal{B}(K_n)}(p_n^{i-1}, p_n^i) = d_W(p_n^{i-1}, p_n^i) \leq d_W(c(t_{i-1}), c(t_i)) + 2\alpha_n \leq \frac{\ell}{m_n} + 2\alpha_n.
$$
It follows that the Lipschitz constant of $c_n$ on each subinterval $[t_{i-1}, t_i]$ is bounded from above by
$$
\frac{\ell/m_n + 2\alpha_n}{1/m_n} = \ell + 2m_n\alpha_n = \ell + o(1),
$$
and therefore $\Lip_{d_{\mathcal{B}(K_n)}}(c_n)\leq \ell + o(1)$.
Moreover, for any $t \in [t_{i-1}, t_i]$, we have
$$
d_W(c_n(t), c(t)) \leq d_W(c_n(t), p_n^{i-1}) + d_W(p_n^{i-1}, c(t_{i-1})) + d_W(c(t_{i-1}), c(t)) \leq \frac{2\ell}{m_n} + 3\alpha_n,
$$
which implies that $c_n \to c$ uniformly in $(W, d_W)$. Regarding the curve lengths, we have
$$
l_{d_{\mathcal{B}(K_n)}}(c_n) = \sum_{i=1}^{m_n} d_{\mathcal{B}(K_n)}(p_n^{i-1}, p_n^i) \leq \ell + 2m_n\alpha_n=\ell+o(1).
$$
Since $\ell =l_{d_{\mathcal{B}(K)}}(c)\leq \liminf_{n \to \infty} l_{d_{\mathcal{B}(K_n)}}(c_n)$ by the lower semicontinuity of curve length under uniform convergence in $(W, d_W)$, we conclude that $l_{d_{\mathcal{B}(K_n)}}(c_n)\to\ell$ as $n\to\infty$.

Next, we prove the $L^2$-convergence of the derivatives of the curves $\pi_{K_n}\circ c_n$. Since the projection $\Pi : W \to \overline{\BB^3}$ is $1$-Lipschitz and preserves the lengths of $c$ and $c_n$, we have $\|(\pi_K \circ c)'\| = \ell$ and $\|(\pi_{K_n} \circ c_n)'\| \leq \ell + o(1)$ almost everywhere on $[0,1]$. The uniform convergence $\pi_{K_n} \circ c_n \to \pi_K \circ c$ combined with this uniform $L^\infty$-bound ensures that
$$
(\pi_{K_n} \circ c_n)' \rightharpoonup (\pi_K \circ c)' \quad \text{weakly in } L^2([0,1]; \RR^3).
$$
Since
$$
\limsup_{n \to \infty} \|(\pi_{K_n} \circ c_n)'\|_{L^2}^2 \leq \ell^2 = \|(\pi_K \circ c)'\|_{L^2}^2,
$$
it follows that $(\pi_{K_n} \circ c_n)' \to (\pi_K \circ c)'$ strongly in $L^2([0,1]; \RR^3)$. 

Finally, suppose that $c([0,1]) \subset \mathcal{B}_{\HH}(K)$. Then $(\pi_K \circ c)([0,1])$ is a compact subset of the open unit ball $\BB^3$. By uniform convergence, there exists a compact ball $B \subset \BB^3$ containing $(\pi_K \circ c)([0,1])$ and $(\pi_{K_n} \circ c_n)([0,1])$ for all sufficiently large $n$; in particular, $c_n([0,1])\subset\mathcal{B}_{\HH}(K_n)$. Let $F(x, v) := \sqrt{g_{\HH}(x)(v,v)}$ denote the Riemannian norm in the Klein model $\HH^3_K$. We claim that there exists a constant $C > 0$ such that for all $x, y \in B$ and $v, w \in \RR^3$,
\begin{equation}\label{eq:Klein-norm}
    |F(x, v) - F(y, w)| \leq |F(x,v)-F(x,w)|+ |F(x,w)-F(y,w)|\leq C\|v - w\| + C\|x - y\|\|w\|.
\end{equation}
Indeed, the triangle inequality for the norm $F(x,\cdot)$ gives $|F(x,v)-F(x,w)|\leq F(x,v-w)\leq C\|v-w\|$. Moreover, the smoothness and uniform positive definiteness of $g_{\HH}$ on the compact set $B$ imply that, for $w\neq0$,
$$
 |F(x,w)-F(y,w)|
 =\frac{|(g_{\HH}(x)-g_{\HH}(y))(w,w)|}
       {F(x,w)+F(y,w)}
 \leq \frac{C\|x-y\|\|w\|^2}{\|w\|}= C\|x-y\|\|w\|,
$$
which proves the claim. Let $\gamma:=\pi_K\circ c$ and $\gamma_n:=\pi_{K_n}\circ c_n$. Applying \eqref{eq:Klein-norm} to the integrands of the hyperbolic lengths yields
$$
\begin{aligned}
\bigl| l_{\HH}(\gamma_n) - l_{\HH}(\gamma) \bigr|
&\leq \int_0^1 |F(\gamma_n(t), \gamma_n'(t)) - F(\gamma(t), \gamma'(t))| \, dt \\
&\leq C \int_0^1 \|\gamma_n'(t) - \gamma'(t)\| \, dt + C \int_0^1 \|\gamma_n(t) - \gamma(t)\| \|\gamma'(t)\| \, dt \\
&\leq C \|\gamma_n' - \gamma'\|_{L^1} + C \|\gamma_n - \gamma\|_\infty \|\gamma'\|_{L^1}.
\end{aligned}
$$
By H\"older's inequality, $\|\gamma_n' - \gamma'\|_{L^1} \leq \|\gamma_n' - \gamma'\|_{L^2} \to 0$. Since $\|\gamma_n - \gamma\|_\infty \to 0$ and $\|\gamma'\|_{L^1} = \ell < \infty$, the right-hand side vanishes as $n \to \infty$. This completes the proof.
\end{proof}

\newtheorem*{proposition2}{\cref{prop:ccv-finite-stability}}
\begin{proposition2}
Suppose that sequences $p_n, q_n \in \mathcal{B}_{\HH}(K_n)$ converge in $(W, d_W)$ to $p, q \in \mathcal{B}_{\HH}(K)$, respectively. Then the hyperbolic intrinsic distances satisfy
\begin{equation}\label{app:ccv-finite-stability-equation}
  \delta_{\mathcal{B}(K_n)}(p_n, q_n) \to \delta_{\mathcal{B}(K)}(p, q) \quad \text{in } [0, +\infty] \quad \text{as } n \to \infty.
\end{equation}
In particular, $\delta_{\mathcal{B}(K_n)}(p_n, q_n) \to +\infty$ whenever $p$ and $q$ lie in different connected components of $\mathcal{B}_{\HH}(K)$.
\end{proposition2}

\begin{proof}
First, suppose that $\delta_{\mathcal{B}(K)}(p, q) < \infty$. For any $\eta > 0$, choose a rectifiable curve $c : [0,1] \to \mathcal{B}_{\HH}(K)$ connecting $p$ to $q$, parametrized with constant speed with respect to $d_{\mathcal{B}(K)}$, such that $l_{\HH}(\pi_K \circ c) \leq \delta_{\mathcal{B}(K)}(p, q) + \eta$. By \cref{ccv:curve-recovery}, there exist curves $c_n : [0,1] \to \mathcal{B}_{\HH}(K_n)$ connecting $p_n$ to $q_n$ such that $l_{\HH}(\pi_{K_n} \circ c_n) \to l_{\HH}(\pi_K \circ c)$. Consequently,
$$
\limsup_{n \to \infty} \delta_{\mathcal{B}(K_n)}(p_n, q_n) \leq \lim_{n \to \infty} l_{\HH}(\pi_{K_n} \circ c_n) = l_{\HH}(\pi_K \circ c) \leq \delta_{\mathcal{B}(K)}(p, q) + \eta.
$$
Letting $\eta \to 0$, we obtain $\limsup_{n \to \infty} \delta_{\mathcal{B}(K_n)}(p_n, q_n) \leq \delta_{\mathcal{B}(K)}(p, q)$. This inequality holds trivially when $\delta_{\mathcal{B}(K)}(p, q) = +\infty$.

It remains to establish the reverse inequality $\liminf_{n \to \infty} \delta_{\mathcal{B}(K_n)}(p_n, q_n) \ge \delta_{\mathcal{B}(K)}(p, q)$. This holds vacuously if the limit inferior is $+\infty$. Otherwise, after passing to a subsequence, we may assume that
$$
L := \liminf_{n \to \infty} \delta_{\mathcal{B}(K_n)}(p_n, q_n) < \infty
$$
is actually the limit of $\delta_{\mathcal{B}(K_n)}(p_n, q_n)$. For each $n$, choose a curve $\gamma_n : [0,1] \to \mathcal{B}_{\HH}(K_n)$ joining $p_n$ and $q_n$ such that
$$
l_{\HH}(\pi_{K_n} \circ \gamma_n) \leq \delta_{\mathcal{B}(K_n)}(p_n, q_n) + \frac{1}{n}.
$$
Then the hyperbolic lengths $l_{\HH}(\pi_{K_n} \circ \gamma_n)$ are uniformly bounded. Since $p_n \to p \in \mathcal{B}_{\HH}(K)$, it follows that the images $(\pi_{K_n} \circ \gamma_n)([0,1])$ are all contained in some compact ball $B \subset \BB^3$ for sufficiently large $n$. 
Parametrize each $\gamma_n$ with constant speed with respect to the intrinsic Euclidean metric $d_{\mathcal{B}(K_n)}$. Since the Euclidean metric is bounded from above by the hyperbolic metric on $\HH^3_K$ (see \eqref{eq:H-bound-E}), we have
$$
l_{d_{\mathcal{B}(K_n)}}(\gamma_n) \leq l_{\HH}(\pi_{K_n} \circ \gamma_n) \leq C < \infty.
$$
Hence, the curves $\gamma_n : [0,1] \to (W, d_W)$ are equicontinuous and take values in the compact metric space $(W, d_W)$. By the Arzel\`a--Ascoli theorem, up to a further subsequence, $\gamma_n$ converges uniformly in $(W, d_W)$ to a continuous curve $\gamma : [0,1] \to \mathcal{B}(K)$ joining $p$ and $q$. Since $\pi_{K_n} \circ \gamma_n \to \pi_K \circ \gamma$ uniformly in $B \subset \BB^3$, the limit curve satisfies $\gamma([0,1]) \subset \mathcal{B}_{\HH}(K)$.
Finally, the lower semicontinuity of hyperbolic length with respect to uniform convergence in $B\subset\HH^3_K$ implies that
$$
\delta_{\mathcal{B}(K)}(p, q) \leq l_{\HH}(\pi_K \circ \gamma) \leq \liminf_{n \to \infty} l_{\HH}(\pi_{K_n} \circ \gamma_n) = \liminf_{n \to \infty} \delta_{\mathcal{B}(K_n)}(p_n, q_n).
$$
This completes the proof.
\end{proof}

\section{Conformal Deformations and Curvature Bounds}
\label{sec:appendix}

In this appendix, we record the results used in \cref{subsec:complete-blowup}. We use the canonical
conformal structure recalled in
\cref{sec:prelim-finite-type}.

\subsection{Conformal deformation}
A surface $(S,d_0)$ has \emph{locally bounded integral curvature} (BIC) if its length metric induces the surface topology and locally is a uniform limit of smooth Riemannian distances with uniformly bounded total absolute curvature. Its intrinsic curvature is a signed Radon measure $\omega$, and its area measure $\mu$ is the normalized two-dimensional Hausdorff measure.

Away from the cusps of $S$, i.e., points $p\in S$ such that $\omega(\{p\})\geq 2\pi$, Reshetnyak's representation theorem gives conformal charts in which
\begin{equation}\label{app:eq:representation}
 \dd s=e^\sigma|\dd z|,\qquad
 \mu=e^{2\sigma}\,\dd x\dd y,\qquad
 \omega=-\Delta_0\sigma.
\end{equation}
Here $\sigma\in L^1_{\mathrm{loc}}$ is a difference of subharmonic functions, and $\Delta_0=\partial^2_x+\partial^2_y$ is the Laplacian in the conformal coordinate. The length element is interpreted using the canonical logarithmic-potential representative of $\sigma$, i.e., locally we choose
\begin{equation}\label{app:eq:potential}
 \sigma(z)=\frac1{2\pi}\int\log\frac1{|z-\zeta|}\,\dd\omega(\zeta)+h(z),
\end{equation}
where $h$ is harmonic. We refer the reader to
\cite{Reshetnyak,FillastreSubharmonic} for further details.

\begin{definition}
Let $S$ be a Riemannian surface. A \emph{subharmonic metric} $d$ on $S$ is a metric with length element given in local conformal coordinates $z_\alpha$ by compatible expressions
$$
 ds = e^{\sigma_\alpha}|\mathrm{d}z_\alpha|,
$$
where each $\sigma_\alpha$ can be written as a difference of subharmonic functions, and hence is well-defined up to a set of zero one-dimensional Hausdorff measure. 
More precisely, for any arc $\gamma \subset S$ that is piecewise rectifiable in conformal coordinates, its length $l_d(\gamma)$ is defined by
$$
 l_d(\gamma) = \sum_{i} \int_{\gamma_i} e^{\sigma_{\alpha_i}} |\mathrm{d}z_{\alpha_i}|,
$$
where $\gamma$ is subdivided into subarcs $\gamma_i$ each contained in a single coordinate chart. For any $x,y\in S$, we then have $d(x,y) = \inf_{\gamma} l_d(\gamma)$,
where the infimum is taken over all piecewise rectifiable arcs $\gamma \subset S$ joining $x$ and $y$.
\end{definition}

It is clear that the metric on any BIC surface is a subharmonic metric based on the above facts. Conversely, by \cite[Theorem 4.66]{FillastreSubharmonic} we have
\begin{theorem}\label{app:reshetnyak-converse}
Any subharmonic metric is a metric of bounded curvature.
\end{theorem}

Recall the definition of the distributional Laplacian $\boldsymbol{\Delta}$ in \cref{subsec:complete-blowup}. We have the following proposition:

\begin{proposition}
\label{app:continuous-change}
Let $(S,d)$ be a BIC surface without cusps. Suppose that $a$ is
continuous on a connected open set $V\subset S$ and
$\boldsymbol\Delta a$ is a signed Radon measure. The path metric
$\widetilde d$ with length element $e^a\dd s_d$ is a BIC metric,
has the same conformal structure as $d$, and satisfies
\begin{equation}\label{app:eq:change}
 \mu_{\widetilde d}=e^{2a}\mu_d,\qquad
 \omega_{\widetilde d}=\omega_d-\boldsymbol\Delta a,\qquad
 \omega_{\widetilde d}(\{p\})=\omega_d(\{p\}).
\end{equation}
\end{proposition}
\begin{proof}
Since $e^{a}$ is a continuous positive function, $\widetilde{d}$ is locally bi-Lipschitz equivalent to $d$. The length element of the metric $\tilde{d}$ is given in conformal
coordinates by compatible expressions
$$
 ds=e^{a_\alpha+\sigma_\alpha}|\mathrm dz_\alpha|.
$$
Since $a$ is continuous and $\Delta_0a$ is a signed Radon
measure, the local logarithmic-potential representation
and Weyl's lemma imply that outside a set of zero one-dimensional Hausdorff measure, $a$ is locally a difference
of two subharmonic functions; see \cite[Section~3]{FillastreSubharmonic}. It then follows from \cref{app:reshetnyak-converse} that $(V,\widetilde{d})$ is a BIC surface. 

The first two identities in \eqref{app:eq:change} are immediate from the definition. We now show that $\Delta_0a (\{p\})=0$ for each $p\in V$. Fix $p\in V$ and choose a conformal coordinate with $z(p)=0$. Let $\chi$ be a smooth cutoff function in the unit disk with $\chi(0)=1$ and set $\chi_r(z)=\chi(z/r)$. Then in a conformal chart near $p$, we have
$$
 |\langle\Delta_0a,\chi_r\rangle|=\left|\int a\Delta_0\chi_r\,\dd x\dd y\right|
 =\left|\int(a-a(0))\Delta_0\chi_r\,\dd x\dd y\right|
 \le\|\Delta_0\chi\|_{L^1}\sup_{B(0,r)}|a-a(0)|\to0.
$$
Since $\Delta_0a$ is a Radon measure, by the dominated convergence theorem we have $\Delta_0a (\{p\})=0$. This proves the last identity.
\end{proof}

By \cite[Theorem~3.2 and the proof of Theorem~1.1]{AmbrosioBertrand}, every $\CBB(-1)$ surface satisfies
\begin{equation}\label{app:eq:cbb-measures}
 \omega_d\ge-\mu_d,
 \qquad
 \omega_d(\{p\})<2\pi
 \quad\text{for every }p\in S.
\end{equation}
The following result then gives a characterization of $\CBB(-1)$ surfaces with the language of BIC surfaces. 

\begin{theorem}[(Marot, {\cite[Theorem~1]{MarotCurvatureBounds}})]
\label{app:cbb-criterion}
A complete BIC surface $(S,d)$ satisfying
$$
 \omega_d\ge-\mu_d,\qquad
 \omega_d(\{p\})<2\pi\quad\text{for every }p\in S
$$
is a complete $\CBB(-1)$ surface.
\end{theorem}

\begin{corollary}\label{app:deformation-criterion}
Let $(S,d)$ be a $\CBB(-1)$ surface and let $a$ be a continuous function on a connected open set $V\subset S$. If $a$ satisfies the distributional inequality
$$
 \boldsymbol\Delta a\le(e^{2a}-1)\mu_d,
$$
and the path metric $\widetilde d$ induced by the length element $e^a\dd s_d$ is complete, then $(V,\widetilde d)$ is a $\CBB(-1)$ surface.
\end{corollary}
\begin{proof}
The inequality ensures that $\boldsymbol\Delta a$ is a signed Radon measure. Proposition~\ref{app:continuous-change} and
\eqref{app:eq:cbb-measures} then imply
$$
 \omega_{\widetilde d}+\mu_{\widetilde d}
 =\omega_d+\mu_d+(e^{2a}-1)\mu_d-\boldsymbol\Delta a\ge0
$$
and $\omega_{\widetilde d}(\{p\})=\omega_d(\{p\})<2\pi$. The result then follows from \cref{app:cbb-criterion}.
\end{proof}

\subsection{Weak inequalities used in \cref{subsec:complete-blowup}}
We first identify the distributional Laplacian in conformal
coordinates with the weak Laplacian defined by the intrinsic Dirichlet form.
Let $(S,d)$ be a $\CBB(-1)$ surface without boundary, endowed with
the area measure $\mu$ and curvature measure $\omega$.
For every locally $d$-Lipschitz function $f$, the distributional Laplacian
$\boldsymbol\Delta f$ admits the following Dirichlet representation:
in any isothermal chart $z=x+iy:\mathcal O\to\Omega\subset\CC$, after identifying
functions with their coordinate expressions, we have
$$
\langle\boldsymbol\Delta f,\varphi\rangle
=
-\int_{\Omega}\nabla_0f\cdot\nabla_0\varphi\,dx\,dy
=
-\int_{\mathcal O}
\langle\nabla f,\nabla\varphi\rangle\,d\mu,
\qquad
\varphi\in C_c^\infty(\mathcal O),
$$
where $\nabla_0$ and $\nabla$ denote the Euclidean gradient in the
$z$-coordinates and the intrinsic gradient with respect to $d$, respectively.

To verify this, let $\omega=\omega^+-\omega^-$ be the Jordan decomposition
of $\omega$, and let $e^{2\sigma}|dz|^2$ be the local conformal representation of the metric. Since $\omega^+(\{p\})<2\pi$ at every point $p\in S$, we may shrink the chart if necessary and choose $q>1$ such that
$\omega^+(\mathcal O)<2\pi/q$. The local potential representation \eqref{app:eq:potential} together with \cite[Lemma~2.10]{AmbrosioBertrand}
then yields $e^{2\sigma}\in L^q_{\mathrm{loc}}(\Omega)$. Moreover, since $\omega^-\leq\mu$ by \cite[Theorem~3.2]{AmbrosioBertrand},
the Radon--Nikodym density of $\omega^-$ with respect to $dx\,dy$ also belongs to $L^q_{\mathrm{loc}}(\Omega)$, and hence its logarithmic potential is locally bounded by H\"older's inequality. Combined with the fact that the potential of $\omega^+$ is locally bounded from below, \eqref{app:eq:potential} implies that $e^\sigma$ is locally bounded away from zero on $\Omega$. Consequently, every test function $\varphi\in C_c^\infty(\mathcal O)$ is $d$-Lipschitz. Now, on a sufficiently small neighborhood where $f$ is $L$-Lipschitz with respect to $d$, we see that $f$ is absolutely continuous on almost every coordinate line and satisfies
$$
|\partial_x f|,\ |\partial_y f|
\leq L e^\sigma
\qquad\text{a.e. in }\Omega.
$$
Since $e^\sigma\in L^{2q}_{\mathrm{loc}}(\Omega)\subset L^2_{\mathrm{loc}}(\Omega)$, it follows that $f\in W_{\mathrm{loc}}^{1,2}(\Omega)$.
Integrating by parts against the Euclidean Laplacian $\Delta_0$, we obtain
$$
\langle\boldsymbol\Delta f,\varphi\rangle
= \int_{\Omega}f\Delta_0\varphi\,dx\,dy
= -\int_{\Omega}\nabla_0f\cdot\nabla_0\varphi\,dx\,dy.
$$
Since $\langle\nabla f,\nabla\varphi\rangle = e^{-2\sigma}\nabla_0 f\cdot\nabla_0\varphi$ almost everywhere and $d\mu=e^{2\sigma}\,dx\,dy$, the cancellation of the conformal factors yields the desired Dirichlet representation. As a consequence, the local Laplacian estimates in \cite[Claim~1.1]{PetruninCD} apply to $\boldsymbol\Delta$. Furthermore, for any locally $d$-Lipschitz function $f$ whose distributional
Laplacian is a signed Radon measure, the standard weak chain rule gives
$$
\boldsymbol\Delta(\Psi\circ f)
=
\Psi'(f)\,\boldsymbol\Delta f
+
\Psi''(f)|\nabla f|^2\,\mu
$$
for every $\Psi\in C^2(I)$, where $I\subset\RR$ is an open interval
containing the image of $f$.

\begin{lemma}[(Radial comparison)]\label{app:radial-comparison}
Let $(S,d)$ be a complete $\CBB(-1)$ surface and let $p\in S$.
Then the locally Lipschitz function $h=\cosh d(p,\cdot)$ satisfies
$$
 \boldsymbol\Delta h\le2h\mu_d,\qquad
 |\nabla h|^2=h^2-1\quad\mu_d\text{-almost everywhere}.
$$
\end{lemma}
\begin{proof}
Set $r := d(p,\cdot)$. Since $r$ is $1$-Lipschitz, $h = \cosh r$ is locally Lipschitz on $S$. Applying \cite[Theorem~8.23\textup{(b)}]{AKP24} to the modified distance function $\operatorname{md}_{-1}(r) = \cosh r - 1 = h - 1$, we obtain 
\begin{equation}\label{eq:h}
    (h \circ \gamma)'' \le h \circ \gamma
\end{equation}
in the distributional sense along every unit-speed geodesic $\gamma$ in $S$, including at $p$. Fix $\varepsilon > 0$. By continuity, every point in $S$ has a small open neighborhood $V$ on which $\sup_V h - \inf_V h < \varepsilon$. Set $\lambda_V := \sup_V h$. By \eqref{eq:h}, the function $h$ is $\lambda_V$-concave on $V$, i.e., the map $t \longmapsto h(\gamma(t)) - \frac{\lambda_V}{2}t^2$ is concave for every unit-speed geodesic $\gamma$ contained in $V$. By the local Laplacian estimate on $2$-dimensional Alexandrov spaces \cite[Claim~1.1]{PetruninCD}, we have on $V$:
$$
\boldsymbol\Delta h \le 2\lambda_V \mu_d \le 2(h + \varepsilon)\mu_d.
$$
By a partition of unity argument, this distributional inequality holds globally on $S$. Letting $\varepsilon \to 0$ yields $\boldsymbol\Delta h \le 2h\mu_d$.

Finally, the distance function $r$ satisfies $|\nabla r| = 1$ $\mu_d$-almost everywhere on $S \setminus \{p\}$. Since $\mu_d(\{p\}) = 0$, the chain rule for locally Lipschitz functions gives $\nabla h = (\sinh r)\nabla r$, and hence
$$
|\nabla h|^2 = \sinh^2 r\, |\nabla r|^2 = h^2 - 1 \qquad \mu_d\text{-almost everywhere}.
$$
\end{proof}

\begin{lemma}\label{app:minimum}
Let $\mu$ be a measure on an open set $U\subset\CC$ defined by $\alpha\,\dd x\dd y$ for some positive function $\alpha\in L^1_{\mathrm{loc}}(U)$, and let $F$ be a continuous function.
If continuous functions $b,c$ satisfy the distributional inequality $\Delta_0b\leq F(b)\mu$ and $\Delta_0c\leq F(c)\mu$, then we have
$$
 \Delta_0\min\{b,c\}\leq F(\min\{b,c\})\mu.
$$
\end{lemma}
\begin{proof}
Let $\varphi \in C_c^\infty(U)$ with $\varphi \ge 0$, and denote $K := \operatorname{supp}\varphi \Subset U$. It suffices to show that
$$
\int_U \min\{b,c\} \Delta_0\varphi \,\dd x\dd y \leq \int_U \varphi F(\min\{b,c\}) \alpha \,\dd x\dd y.
$$
For any parameter $\eta > 0$, consider the smooth approximation of the minimum function $m_\eta \colon \mathbb{R}^2 \to \mathbb{R}$ defined by
$$
m_\eta(s,t) := -\eta \log \bigl( e^{-s/\eta} + e^{-t/\eta} \bigr).
$$
A direct computation shows that $m_\eta \to \min\{s,t\}$ uniformly on $\mathbb{R}^2$ as $\eta \to 0^+$. Its gradient $(\theta_\eta, 1-\theta_\eta)$ has coefficients in $(0,1)$ with $\theta_\eta(s,t) = (1+e^{(s-t)/\eta})^{-1}$, and its Hessian is negative semi-definite everywhere.

Let $\rho_\varepsilon$ be a standard mollifier. For sufficiently small $\varepsilon > 0$, the smooth functions $b_\varepsilon := b * \rho_\varepsilon$ and $c_\varepsilon := c * \rho_\varepsilon$ satisfy $b_\varepsilon \to b$ and $c_\varepsilon \to c$ uniformly on $K$. Since $\Delta_0 b \leq F(b)\mu$ in the sense of distributions and $\rho_\varepsilon \ge 0$, we have $\Delta_0 b_\varepsilon \leq (F(b)\alpha) * \rho_\varepsilon$ on a neighborhood of $K$, and similarly for $c_\varepsilon$.
Setting $u_{\eta,\varepsilon} := m_\eta(b_\varepsilon, c_\varepsilon)$, the chain rule yields
$$
\begin{aligned}
    \Delta_0 u_{\eta,\varepsilon} 
    &\leq \theta_\eta(b_\varepsilon, c_\varepsilon) \Delta_0 b_\varepsilon + \bigl(1 - \theta_\eta(b_\varepsilon, c_\varepsilon)\bigr) \Delta_0 c_\varepsilon,\\
    &\leq \theta_\eta(b_\varepsilon, c_\varepsilon) \bigl( (F(b)\alpha) * \rho_\varepsilon \bigr) + \bigl(1 - \theta_\eta(b_\varepsilon, c_\varepsilon)\bigr) \bigl( (F(c)\alpha) * \rho_\varepsilon \bigr)
\end{aligned}
$$
pointwise on $K$, where the non-positive Hessian quadratic form is dropped. Testing against $\varphi$ and using integration by parts on the left-hand side, we have
$$
\int_U u_{\eta,\varepsilon} \Delta_0\varphi \,\dd x\dd y \leq \int_U \varphi \Bigl[ \theta_\eta(b_\varepsilon, c_\varepsilon) \bigl( (F(b)\alpha) * \rho_\varepsilon \bigr) + \bigl(1 - \theta_\eta(b_\varepsilon, c_\varepsilon)\bigr) \bigl( (F(c)\alpha) * \rho_\varepsilon \bigr) \Bigr] \,\dd x\dd y.
$$

Since $\theta_\eta(b_\varepsilon, c_\varepsilon) \to \theta_\eta(b,c)$ uniformly on $K$ and the mollifications converge in $L^1(K)$, letting $\varepsilon \to 0$ gives
$$
\int_U m_\eta(b,c) \Delta_0\varphi \,\dd x\dd y \le \int_U \varphi \Bigl[ \theta_\eta(b,c) F(b) + \bigl(1 - \theta_\eta(b,c)\bigr) F(c) \Bigr] \alpha \,\dd x\dd y.
$$
Finally, as $\eta \to 0^+$, $m_\eta(b,c) \to \min\{b,c\}$ uniformly on $K$. On the other hand, since $\theta_\eta(b,c) \to 1_{\{b<c\}} + \frac{1}{2}1_{\{b=c\}}$ pointwise on $K$,
$$
\theta_\eta(b,c) F(b) + \bigl(1 - \theta_\eta(b,c)\bigr) F(c) \to F(\min\{b,c\}).
$$
Since $F(b)$ and $F(c)$ are bounded on $K$ and $\varphi\alpha \in L^1(K)$, the dominated convergence theorem completes the proof.
\end{proof}




\bibliographystyle{amsplain}
\bibliography{main}

\noindent Xinrong Zhao, xrzhao24@m.fudan.edu.cn\\
\emph{School of Mathematical Sciences, Fudan University, Shanghai, 200433, P.R. China}\\[-8pt]

\end{document}